\documentclass[12pt,a4paper]{amsart}

\usepackage{amsmath, amsfonts, xifthen, latexsym, amssymb, amsthm, amscd}

\usepackage[utf8]{inputenc}
\usepackage{graphicx}
\usepackage{url}

\usepackage{hyperref}
\hypersetup{colorlinks=true,citecolor=blue,filecolor=blue,linkcolor=blue,urlcolor=blue}
\usepackage[margin=1.4in]{geometry}

\newtheorem{theorem}{Theorem}
\newtheorem{lemma}{Lemma}[section]

\newtheorem{proposition}{Proposition}[section]

\newtheorem{conjecture}{Conjecture}[section]

\newtheorem{corollary}{Corollary}[section]

\theoremstyle{definition}
\newtheorem{definition}{Definition}[section]

\newtheorem{remark}{Remark}

\newtheorem{example}{Example}

\newcommand{\conf}{\mbox{Conf}}

\newcommand{\eps}{{\varepsilon}}

\renewcommand{\phi}{{\varphi}}
\newcommand{\grd}{{Gr_d}}

\newcommand{\Stab}{\operatorname{Stab}}

\newcommand\diam{\operatorname{diam}}

\newcommand{\cB}{\mathcal{B}}\newcommand{\cM}{\mathcal{M}}
\newcommand{\cCB}{\mathcal{CB}}
\newcommand\R{\mathbb R}
\newcommand\Q{\mathbb Q}
\newcommand\Z{\mathbb Z}

\newcommand\N{\mathbb N}
\newcommand\cG{\mathcal{G}}
\newcommand{\cA}{\mathcal{A}}
\newcommand\rg{\mathcal{RG}_d}
\newcommand\crg{\mathcal{CRG}_d}
\newcommand\drg{d_{\mathcal{RG}}}
\newcommand\gsg{\Gamma_\Sigma \cG}
\newcommand\tbkr{\cB_r(\Gamma,\Sigma,k)}
\newcommand\tblr{\cB_r(\Gamma,\Sigma,l)}
\newcommand\tbls{\cB_s(\Gamma,\Sigma,l)}
\newcommand\tbks{\cB_s(\Gamma,\Sigma,k)}
\newcommand\tbrp{\cB_{r+1}(\Gamma,\Sigma,r)}\newcommand\tbr{\cB_{r}(\Gamma,\Sigma,r)}

\newcommand{\Free}{\operatorname{Free}}
\newcommand{\ha}{\hat{A}}

\newcommand{\actson}{\curvearrowright}

\newcommand{\cC}{\mathcal{C}} \newcommand{\cl}{\mathcal{L}}\newcommand{\cN}{\mathcal{N}}\newcommand{\cH}{\mathcal{H}}

\newcommand{\cal}[1]{{\mathcal #1}}
\newcommand{\ov}{\overline}
\newcommand{\tart}{\stackrel{n}{\rightarrow}}
\newcommand{\vi}{\vskip 0.1in \noindent}
\usepackage{etoolbox}
\newtoggle{no_cases}\toggletrue{no_cases}
\newcommand{\case}[2][]{\iftoggle{no_cases}{\left\{\begin{array}{ll}#2 & #1}{\\#2 & #1}\togglefalse{no_cases}}
\newcommand{\esac}{\end{array}\right.\toggletrue{no_cases}}

\newcommand{\Prob}{\operatorname{Prob}}
\newcommand{\stab}{\operatorname{Stab}}\newcommand{\Sub}{\operatorname{Sub}}
\newcommand{\Spec}{\operatorname{Spec}}

\newcommand{\Inv}{{\operatorname{Inv}}}

\newcommand{\supp}{\operatorname{Supp}}
\newcommand{\Cay}{\operatorname{Cay}}

\begin{document}
\title[Qualitative graph limit theory]{Qualitative graph limit theory. \\ Cantor Dynamical Systems and Constant-Time Distributed Algorithms}
\author{Gábor Elek}
\address{Department of Mathematics And Statistics, Fylde College, Lancaster University, Lancaster, LA1 4YF, United Kingdom}

\email{g.elek@lancaster.ac.uk}  

\thanks{The author was partially supported
by the ERC Consolidator Grant "Asymptotic invariants of discrete groups,
sparse graphs and locally symmetric spaces" No. 648017. }

\begin{abstract}
The goal of the paper is to lay the foundation for
the qualitative analogue of the classical, quantitative sparse graph limit 
theory. In the first part of the paper we introduce the qualitative analogues
of the Benjamini-Schramm and local-global graph limit theories for sparse graphs.
The natural limit objects are continuous actions of finitely generated groups on totally disconnected compact metric spaces.
We prove that the space of weak equivalent classes of free Cantor actions is compact and contains a smallest element, as
in the measurable case. We will introduce and study various notions of almost finiteness, the qualitative
analogue of hyperfiniteness, for classes of bounded degree graphs. We prove the almost finiteness of a new class of \'etale groupoids associated
to Cantor actions and construct an example of a nonamenable, almost finite totally disconnected \'etale groupoid, answering a query
of Suzuki.  Motivated by the notions and results on qualitative graph limits, in the second part of our paper we
give a precise definition of constant-time distributed algorithms on sparse graphs. We construct such
constant-time algorithms for various approximation problems for hyperfinite
and almost finite graph classes. We also prove the Hausdorff convergence of
the spectra of convergent graph sequences in the strongly almost finite
category.

\end{abstract}\maketitle
\newpage

\setcounter{tocdepth}{2}
\tableofcontents

\newpage
\section{Introduction and basic definitions}
Sparse graph limit theory deals with very large finite graphs of small vertex degrees (see \cite{Lovasz1} for a recent survey). 
One can study these large graphs by their subgraph statistics. Informally speaking, two large graphs are close to each other in 
the statistical sense if
the frequencies of their small subgraphs do not differ too much. Benjamini and Schramm \cite{BS} defined limit objects for growing
sequences of graphs that are closer and closer to each other in the statistical sense. Another kind of limit objects were defined by
Aldous and Lyons \cite{AL}, and soon turned out that sparse graph limits are intimately related to finitely generated infinite groups and
their measure preserving actions. 
Goldreich and Ron \cite{GR} introduced the notion of property testing and parameter estimation
of bounded degree graphs. A graph parameter can be estimated using finite samplings if and only if for
two graphs that are close in the statistical sense the parameters are close to each other, as well.
In other words, sparse graph limit theory is about the information one can extract from large graphs via
sampling. In his seminal monograph: Metric Structures for Riemannian and Non-Riemannian Spaces 
\cite{Gromov}, Gromov studied sampling (quantitative) and observable (qualitative) convergence of metric measure spaces.
In our paper we develop a study of bounded degree graphs via observables: qualitative graph limit
theory. The table below is intended to briefly summarize the analogies between various notions
of the quantitative and the qualitative approaches.
 
\begin{table}[h!]
  \begin{center}
    \begin{tabular}{l|l} 
      \textbf{Quantitative graph limit theory} & \textbf{Qualitative graph limit theory} \\
      \hline
      Benjamini-Schramm  convergence  & Naive convergence\\
      \hline
      Invariant Random Subgroups &  Uniformly Recurrent Subgroups\\
      \hline
      Measure preserving actions & Stable actions \\
      \hline 
      Hyperfiniteness & Almost finiteness \\
      \hline
      Sofic groups & LEF-groups \\
      \hline
      Local-global convergence & Qualitative convergence \\
      \hline
      Constant-time \\  randomized algorithms & 
      Constant-time  distributed algorithms \\
      \hline
      Property testing & Distributed property testing \\
      \hline
      Weak convergence of the spectra  & Hausdorff convergence of the spectra \\
      \hline
      von Neumann algebras & $C^*$-algebras\\
      \hline
    \end{tabular}
  \end{center}
\end{table}
\subsection{Qualitative graph limits}
Let $d>0$ be a natural number and denote by $\grd$ the set of all countable graphs (up to isomorphisms) of
vertex degree bound $d$.
For $G\in\grd$ and $r>0$, let $\cB_r(G)$ be the set of all rooted balls of radius $r$ contained in $G$. 
\begin{definition}[Naive convergence] \label{defnaive}
A sequence of countable graphs \\ $\{G_n\}^\infty_{n=1}\subset \grd$ is {\bf convergent in the naive sense}
if for any $r>0$ and rooted ball $B$ of radius $r$, there exists $N_B>0$ such that
\begin{itemize}
\item either $B\in \cB_r(G_n)$ for any $n\geq N_B$,
\item or $B\notin \cB_r(G_n)$ for any $n\geq N_B$.
\end{itemize}
\end{definition}
\noindent
One can immediately observe that naive convergence is the qualitative analogue of the Benjamini-Schramm graph convergence (defined for finite graphs).
It is not hard to see that for any graph sequence $\{G_n\}^\infty_{n=1}$ that converges in the naive sense, there
exists a graph $G\in\grd$ such that $G_n\stackrel{n}{\to} G$, that is, $B\in \cB_r(G)$ if and only if  the rooted ball
$B$ is contained in all but finitely many of the graphs $G_n$.
A much finer notion of convergence,  qualitative graph convergence can be defined as the analogue of
the local-global convergence notion introduced by Hatami, Lov\'asz and Szegedy \cite{HLSZ}.
Let $H\in\grd$ and let $Q$ be a finite set. Let $\phi:V(H)\to Q$ be a labeling function.
We can define the $r$-configuration set of $\phi$ in the following way.
Let $U_d^{r,Q}$ denote the finite set of all rooted, $Q$-labeled balls of radius $r$ and vertex degree bound $d$
(up to rooted, labeled isomorphisms). Then $\mbox{Conf}_{r,H}(\phi)\subseteq U^{r,Q}_d$ is the set
of all rooted, $Q$-labeled balls that occur in the labeled graph $(H,\phi)$.
\begin{definition}[Qualitative convergence] \label{defqual}
A sequence of countable graphs $\{G_n\}^\infty_{n=1}\subset \grd$ is {\bf qualitatively convergent}
if for any finite set $Q$, integer $r>1$ and subset $S\subseteq U^{r,Q}_d$, 
\begin{itemize}
\item either there exists $\phi^n:V(G_n)\to Q$, such that $\mbox{Conf}_{r,G_n}(\phi^n)=S$ for any $n\geq N_S$,
\item or there exists no $\phi^n:V(G_n)\to Q$, such that $\mbox{Conf}_{r,G_n}(\phi^n)=S$ for any $n\geq N_S$.
\end{itemize}
\end{definition}
\noindent
Clearly, if $\{G_n\}^\infty_{n=1}$ qualitatively converges then it converges in the naive sense, as well.
On the other hand, let $\{C_n\}^\infty_{n=1}$ be a sequence of cyclic graphs, where $|V(C_n)|=n$.
Then,  $\{C_n\}^\infty_{n=1}$ converges in the naive sense and it does not converge qualitatively (see Proposition \ref{kedd1}).
In Section \ref{qual}, we will define qualitative convergence for $\Gamma$-Schreier graphs as well, where $\Gamma$ is
a finitely generated group. It turns out that the natural limit objects for qualitative convergence of $\Gamma$-Schreier graphs
are stable actions of the group $\Gamma$ on a totally disconnected compact metric space. 
\subsection{The space of free Cantor actions} Qualitative convergence leads to the notion of
weak equivalence of free Cantor actions.  For minimal actions of the integers, such notion has already been introduced by
Lin and Matui \cite{LM}. We will prove (Theorem \ref{compact}) that the space of weak equivalence classes for
free $\Gamma$-actions is compact, similarly to the space of weak equivalence class of measurable actions \cite{AE},\cite{TDR}.
We also prove the qualitative analogue of the Abert-Weiss Theorem \cite{AW} showing that the space above contains
a smallest element, the continuous free analogue of the standard Bernoulli shifts. We will define the
qualitative graph limits for sequences of simple graphs as well and study the associated \'etale groupoids.
\subsection{Almost finiteness}
Recall \cite{Elekcost} that a class of finite graphs $\cG\subset \grd$ is {\bf hyperfinite} if
for any $\eps>0$ there exists $L_\eps>0$ such that 
\begin{itemize}
\item for any $G\in \cG$, we can remove $\eps|V(G)|$ edges such that the resulting graph $G'$ consists of
components of size at most $L_\eps$. 
\end{itemize}
\noindent
The key notion of our paper is almost finiteness, the qualitative analogue of hyperfiniteness.
Almost finiteness was originally introduced by Matui \cite{Matui} for totally disconnected \'etale groupoids, thus,
in particular, for free and even stable Cantor actions. Before getting further, let us recall some definitions from
graph theory.
Let $G\in\grd$ be a countable graph and $H\subset V(G)$ be a finite subset. Then,
\begin{itemize}
\item $\diam_G(H)=\max_{x,y\in H} d_H(x,y)$,
\item the boundary of $H$, $\partial(H)$, is the set of vertices $x\in H$ such that
there exists a vertex $y\in V(G)\backslash H$ adjacent to $x$,
\item the isoperimetric constant of $H$, $i_G(H):=\frac{|\partial(H)|}{|H|}\,.$
\end{itemize}
\begin{definition}
A class of graphs $\cG\subset \grd$ is {\bf almost finite} if for any $\eps>0$, there
exists $K_\eps>0$ such that for each $G\in\cG$ we have a partition $V(G)=\{H_1,H_2,\dots\}$ satisfying the
following properties.
\begin{enumerate}
\item For any $j\geq 1$, $\diam_G(H_j)\leq K_\eps$.
\item For any $j\geq 1$, $i_G(H)\leq \eps\,.$
\end{enumerate}
We call such a partition an $(\epsilon,K_\eps)$-partition. An infinite graph $G\in\grd$ is called an almost finite graph
if the class consisting of the single element $G$ is almost finite. 
\end{definition}
\noindent
Clearly, any almost finite graph class is hyperfinite. On the other hand,
the class of finite trees is hyperfinite, but not almost finite.
 We say that a class $\cG\subset \grd$ is {\bf distributed almost finite} if
for any $\eps>0$ there exists a constant-time distributed algorithm that computes an $(\eps,K_\eps)$-partition for
any $G\in\cG$ (see Section \ref{CTDA} for further details). 
\begin{conjecture} A class $\cG$ is almost finite if and only if it is distributed almost finite.
\end{conjecture}
\begin{remark}
Let $\Gamma$ be a finitely generated amenable group and let $C_\Gamma$ be a Cayley graph of $\Gamma$. 
Using the language of our paper, the breakthrough Tiling Theorem of Downarowicz, Huczek and Zhang \cite{DHZ} states
that the graph $C_\Gamma$ is  almost finite.  On the other hand, the
distributed almost finiteness of $C_\Gamma$ would imply that all free $\Gamma$-actions on the Cantor set are
almost finite in the sense of Matui. Such result would have very important consequences for the reduced $C^*$-algebras of the actions
\cite{Kerr}, \cite{Suzuki}.
\end{remark}
\noindent 
The following two definitions are motivated by the notion of fractional hyperfiniteness introduced by Lov\'asz \cite{Lovhyp}.
\begin{definition}
A class $\cG\subset \grd$ is called {\bf strongly almost finite} if for any $\eps>0$, there exists
$T_\eps>0$ and $K_\eps>0$ such that for any $G\in\cG$ we have $T_\eps$ pieces of $(\eps,K_\eps)$-partitions
$\{H^i_1,H^i_2,\dots\}_{i=1}^{T_\eps}$ and
for any $x\in V(G)$
$$\frac{|\{i\,\mid \, x\in\partial(H^i_j)\,\mbox{for some $1\leq j$}\}|}{T_\eps}<\eps\,.$$
\end{definition}
\noindent
\begin{definition} \label{weaklyalm}
A countable graph $G$ is {\bf fractionally almost finite} if for any $\eps>0$, there exists
$T_\eps>0$ and $K_\eps>0$  such that we have $T_\eps$ amount of partitions of $V(G)$,
$\{H^i_1,H^i_2,\dots\}_{i=1}^{T_\eps}$ into subsets of diameter at most $K_\eps$, so that
for any $x\in V(G)$ \begin{equation} \label{este1}
\frac{|\{i\,\mid \, x\in\partial(H^i_j)\,\mbox{for some $1\leq j$}\}|}{T_\eps}<\eps\,.\end{equation}
\end{definition}
\noindent
We will see that regular trees (that are clearly not almost finite graphs) are actually fractionally almost finite.
The main technical result of our paper is that $D$-doubling graphs are, in fact, strongly almost finite and even
distributed strongly almost finite, see Section \ref{CTDA}  (Theorem \ref{fotetel}).
This theorem implies that stable actions with doubling graph structure are almost finite in the sense of Matui. Many minimal
stable actions with the doubling property was constructed in \cite{Elekurs}, hence by the Main Result of \cite{Suzuki} they all amount
to new examples of simple $C^*$-algebras of stable rank one. 
We will also construct an example of an almost finite minimal \'etale groupoid (of a minimal stable action) which is not topologically
amenable, answering a query of Suzuki \cite{Suzuki}.

\subsection{Constant-time distributed algorithms} The second part of our paper is an application in computer science, and it
is strongly related to the notions and results of the first part.  
Let us recall the classical $\mathcal{LOCAL}$-model of distributed graph algorithms (see e.g. \cite{BE}). 
Let us fix a constant $d>0$ for the rest of the paper. Let $G=(V,E)$ be
a simple graph of vertex degree bound $d$, such that the vertices have a unique ID from the
set $\{1,2,\dots,|V|\}$. The vertices are identified with processors and the edges between adjacent
vertices are identified with communication ports. In each round each vertex $x\in V$
\begin{itemize}
\item  sends some message to each of its adjacent vertices,
\item  receives some message from each of its adjacent vertices,
\item  performs some calculation based on all the received messages.
\end{itemize}
After a certain amount of rounds (the running time) the procedure halts and each of the vertices produces 
some output, e.g.
an element of a given finite set $F$. Note that the local calculation performed by the
vertices can be unbounded and are not taken into consideration in the calculation of the running time.
Also, we do not bound
the length of the individual messages. The process starts with the communication of the ID's.
The simplest and most basic distributed graph algorithm problem is the vertex colouring problem that
are used for breaking the symmetries of the graphs. 
The goal is to produce a legal colouring of the vertices of the graph $G$ by $(d+1)$-colours (that is,
adjacent vertices must have different colours). Linial \cite{Linial} proved that one needs $O(\log^*(n))$ rounds
for such vertex colouring. On the other hand, such $(d+1)$-colouring can be computed within
$(O(d^2)+\log^*(n))$-time \cite{GPS}. Note that $\log^*(n)$ denoted the iterated logarithm, where
\begin{itemize}
\item $\log^*(n):= 0$ if $n\leq 1$.
\item $\log^*(n):=1 +\log^*(\log_2(n))$, otherwise.
\end{itemize}
\noindent
Our goal is to study distributed algorithms that can be performed in
constant-time provided that  a {\bf symmetry breaking vertex colouring} is already given. 
Let $r>0$ be an integer, $F$ be
a finite set and 
$\cG\subseteq\grd$ be a class of graphs.
Let us assume that the vertices of the graphs $G\in\cG$ are labeled by the finite set $Q$ in such a way
that if $0<d_G(x,y)\leq r$, then the labels of $x$ and $y$ are different. A {\bf constant-time distributed algorithm} starts with such
 $(r,Q)$-labelings and computes a labeling of the vertices of the graphs $G$ by the finite set $F$ in at most $r$ rounds. 
Note that for graphs of size $n$ the required
precolouring can be computed within $O(d^{2r})+\log^*(n)$-time.  
\vskip 0.1in
\noindent
Our definition of a contant-time distributed algorithm is motivated by Theorem \ref{AWe}. We use a local
symmetry breaking mechanism that gives away the minimum amount of information about the graphs on which
our algorithm will be performed.
So, our approach can be viewed as the qualitative analogue of the randomized local framework introduced by
Goldreich and Ron \cite{GR}, which uses a uniformly random labeling of the vertices to break most of the symmetries with
high probability. We will show that for hyperfinite graph classes (see Section \ref{CTDA}), for any $\eps>0$ we
have a constant-time distributed algorithm that produces an $(1-\eps)$-approximation of the maximum independent
set problem or the minimum vertex cover problem (Proposition \ref{hypind}). For a smaller class of graphs we have
a constant-time distributed algorithm even for the weighted unrestricted independent set problem in the
deterministic-randomness framework (Proposition \ref{detran}). For arbitrary graphs, we can show the existence
of a constant-time distributed algorithm that produces an $(1+\eps)$-approximation of the maximum matching
problem (Theorem \ref{matching}). This is a typical infinite-to-finite proof that uses the notion of 
qualitative graph limits. Finally, we consider distributed parameter testing, that is strongly related to naive graph convergence.
We prove a general spectral convergence result (Theorem \ref{spectralconv}) and show that for classes of $D$-doubling graphs the spectrum of the graph can be tested in a distributed fashion.

\section{Graph convergence in the naive sense}\label{tgc}
\noindent The goal of this section is to introduce the notion of naive
convergence and to define the relevant compact spaces of countable graphs.
\subsection{The space of rooted connected graphs}
Until the end of this paper we fix an integer $d>0$. First of all, let us recall the notion
of rooted graph convergence (see e.g. \cite{Elekhyper}). Let $\rg$ be the set of all connected graphs $G$ with vertex degree bound
$d$ and a distinguished vertex (root) $x\in V(G)$.
We can define a metric $\drg$ on the set $\rg$ in the following way. Let $(G,x),(H,y)\in \rg$. Then
$$\drg((G,x),(H,y))=2^{-n}\,,$$
\noindent
where $n$ is the largest integer for which the rooted balls
$B_n(G,x)$ and $B_n(H,y)$ are rooted-isomorphic.
Then, $\rg$ is a compact metric space with respect to $\drg$. We will also consider
the space of Cantor-labeled rooted graphs $\crg$.
An element of $\crg$ is a rooted, connected graph $(G,x)$ equipped with
a vertex labeling $\phi:V(G)\to \{0,1\}^{\N}$ by the Cantor set.
For $s>0$ let $\phi_s:V(G)\to  \{0,1\}^{[s]}$ denote the projection of $\phi$ onto the first $s$ coordinates, where
$[s]=\{1,2,\dots,s\}$. Hence, if $(G,x,\phi)$ is a $\{0,1\}^{\N}$-labeled rooted graph then $(G,x,\phi_s)$ is
a $\{0,1\}^{[s]}$-labeled rooted graph. Again, we can define a metric on $\crg$ by
$$d_{\crg}((G,x,\phi), (H,y,\psi))=2^{-n}\,,$$
\noindent
where $n$ is
the largest integer for which the $\{0,1\}^{[n]}$-labeled rooted balls
$B_n(G,x,\phi_n)$ and $B_n(H,y,\psi_n)$ are
rooted-labeled isomorphic. Again, $\crg$ is a compact metric space with respect to
$d_{\crg}$. Similarly, we can define a compact metric structure on the space $\rg^Q$ of rooted countable graphs
of vertices labeled by the finite set $Q$.
\subsection{The space of countable graphs}
Now let $Gr_d$ be the set of all (not necessarily connected) countable graphs of vertex degree bound $d$.
For $G\in Gr_d$ denote by $\cB(G)$ the set of rooted balls $B$ for which there is an $x\in V(G)$ and $k\geq 1$
such that $B_k(G,x)$ is rooted isomorphic to $B$.
Let $G,H\in Gr_d$. We say that
$G$ and $H$ is equivalent, if $\cB(G)=\cB(H)$. 
We will denote by $\overline{Gr_d}$ the set of equivalence classes of $Gr_d$.
We can define a metric on $\overline{Gr_d}$ as follows.
Let $G,H\in Gr_d$ representing the elements $[G],[H]\in \overline{Gr_d}$. Then
$d_{Gr}([G],[H])= 2^{-n}$ if
\begin{itemize}
\item
For any $1\leq i \leq n$ and ball $B\in U^i_d$ (the set of
all rooted balls of radius $i$ and vertex degree bound $d$), either $B\in \cB(G)$ and $B\in \cB(H)$, 
or $B\notin \cB(G)$ and $B\notin \cB(H)$.
\item There exists $B\in U^{n+1}_d$ such that $B$ is a rooted ball in exactly one of the two graphs.
\end{itemize}
\noindent
Again, we consider the Cantor labeled graphs. An element of the set
$CGr_d$ is a countable graph $G\in Gr_d$ equipped with a vertex labeling
$\phi:V(G)\to \{0,1\}^\N$.
For $k\geq 1$ we denote by $CU^k_d$ the set of rooted balls $B$ of radius $k$ equipped with a vertex labeling
$\rho:V(B)\to \{0,1\}^{[k]}$. So, if $(G,\phi)\in CGr_d$ and $k\geq 1$, then for any $x\in V(G)$ we assign
an element of $CU^k_d$. Now we can proceed exactly the same way as in the unlabeled case. For $G\in CGr_d$
and $B\in CU^k_d$, $B\in \cCB(G)$ if and only if there exists $x\in V(G)$ such that the ball $B_k(G,x,\phi_k)$ is rooted-labeled
isomorphic to $B$.
We say that $(G,\phi)$ and $(H,\psi)$ are equivalent if $\cCB(G)=\cCB(H)$. The set of
equivalence classes will be denoted by $\overline{CGr_d}$. The metric
on $\overline{CGr_d}$ is defined as follows. Let $(G,\phi),(H,\psi)\in CGr_d$ representing
the classes $[(G,\phi)]$ and $[(H,\psi)]$. Then,
$$d_{CGr}((G,\phi),(H,\psi))=2^{-n}\,,$$
if
\begin{itemize}
\item For any $1\leq i \leq n$ and $B\in CU^i_d$, either $B\in \cCB(G)$ and $B\in \cCB(H)$, 
or $B\notin \cCB(G)$ and $B\notin \cCB(H)$.
\item There exists $B\in CU^{n+1}_d$ such that
$B$ is rooted-labeled isomorphic to a $\{0,1\}^{[n+1]}$-labeled ball of exactly one of the two graphs.
\end{itemize}
\noindent
Similarly, we can define a metric on the equivalence classes of countable graphs
$\overline{Gr}^Q_d$ with vertices labeled by a finite set $Q$.

\subsection{Orbit invariant subspaces}
Now, let $G\in Gr_d$. Then we can consider the set $O(G)\subset \rg$, the orbit of $G$. The elements
of $O(G)$ are all the rooted graphs $(G^x,x)$, where $x\in V(G)$ and $G^x$ is the component
of $G$ containing $x$. The orbit closure of $G$, $\overline{O(G)}$ is
the closure of $O(G)$ in the compact metric space $\rg$.
We say that a closed set $M\subseteq \rg$ is {\bf orbit invariant} if for
any $(G,x)\in M$ and $y\in V(G)$, $(G,y)\in M$ as well.
We denote the set of all orbit invariant closed sets by $\Inv(\rg)$.

\begin{proposition}\label{compactpro}
The metric space $\overline{Gr}_d$ is compact and $\overline{O}: \overline{Gr}_d \to \Inv(\rg)$ is a homeomorphism,
where the topology on $ \Inv(\rg)$ is given by the Hausdorff metric and $\overline{O}$ assigns to the graph $G$ its orbit closure.
\end{proposition}
\proof
Let $\{G_n\}^\infty_{n=1}$ be a Cauchy-sequence in $Gr_d$. Let $\cB$ be the set of rooted balls that
are eventually contained in the graphs $\{G_n\}^\infty_{n=1}$.
In order to prove compactness, it is enough to show that there exists $G\in Gr_d$ such that
$\cB(G)=\cB$.
If $B\in \cB\cap U^i_d$, then there exists $n_B>0$ and for any $n\geq n_B$ a vertex $x^B_n\in V(G_n)$
such that the rooted ball of radius $i$ around $x^B_n$ is rooted-isomorphic to $B$.
Let $(G^{x^B_{n_k}}_{n_k}, x^B_{n_k})$ be a convergent sequence in $\rg$ and let $(G_B,x_B)$ be its limit.
Then, $\cB(G_B)\subset \cB$ and $B\in G_B$. Let $G=\cup_B G_B$ be the graph that consists of the disjoint
copies of the graphs $G_B$. Clearly, $\cB(G)=\cB$. Hence, $Gr_d$ is compact.

\noindent
Suppose that $\ov{O(G)}=\ov{O(H)}$ for some $G,H\in Gr_d$.
Then $\cup_{J\in \ov{O(G)}} \cB(J)=\cB(G)$ and also, $\cup_{J\in \ov{O(H)}} \cB(J)=\cB(H)$, therefore $G$ and $H$
are equivalent.
Thus, the map $\ov{O}: \overline{Gr}_d \to \Inv(\rg)$ is injective. Also, if $M\in \Inv(\rg)$ then
it is easy to see that $M=\ov{O(G)}$, where $G$ is the disjoint union of graphs $\{G_n\}^\infty_{n=1}$, where
$(G_n,x_n)$ is a dense set in $M$.
Hence, the map $\ov{O(G)}$ is surjective as well.
Finally, we prove that if $G_n\to G$ in the metric space $Gr_d$, then $\ov{O(G_n)}\to \ov{O(G)}$ in the Hausdorff metric
of $\Inv(\rg)$. Let $(H,x)\in\ov{O(G)}$. We need to show that there exists $y_k\in V(G_k)$ so that
$(G_k,y_k)\to (H,x)$ in the metric space $\rg$. For $k\geq 1$, let $s_k$ be the
largest integer such that the following two conditions are satisfied.
\begin{enumerate}
\item $s_k\leq k$.
\item There exists $y_k\in G_k$ so that
$B_{s_k}(G_k,y_k)$ and $B_{s_k}(H,x)$ are rooted-isomorphic.
\end{enumerate}
\noindent
Since $\{G_n\}^\infty_{n=1}$ converges to $G$ and $(H,x)\in\overline{Gr}_d$ we have that $s_k\to\infty$ as $k\to\infty$.
Thus, $\{(G_k,y_k)\}^\infty_{k=1}$ tends to $(H,x)$ in the metric space $\rg$.
In order to finish our proof we need to show that if $(H,x)\notin \ov{O(G)}$, then there is no subsequence
$(H_{n_k},x_{n_k})\in \ov{O(G_{n_k})}$ such that $(H_{n_k},x_{n_k})\to (H,x)$.
That is, there exists $k\geq 1$ such that if $n$ is large enough, then $B_k(H,x)\notin \cB(H_n)$ if
$(H_n,z)\in\ov{O(G_n)}$.
Since $(H,x)\notin \ov{O(G)}$, there exists $k\geq 1$ so that $B_k(H,x)\notin \cB(G)$. 
Consequently, $B_k(H,x)\notin \cB(G_n)$ if $n$ is large enough. Therefore, if
$n$ is large enough, then $B_k(H,x)\notin \cB(H_n)$ provided that
$(H_n,z)\in\ov{O(G_n)}$. \qed
\vskip 0.1in
\noindent
Similarly, one can prove the following proposition.
\begin{proposition}\label{labeledversion}
The metric spaces $\overline{C\grd}$ and $\grd^Q$ are compact.
Also, $\overline{O}:\overline{C\grd}\to Inv(\crg)$ and 
$\overline{O}:\overline{\grd^Q}\to Inv(\rg^Q)$ are
homeomorphisms.
\end{proposition}
\begin{remark}
Let $\{G_n\}^\infty_{n=1}\subset \grd$ be a sequence of graphs that is convergent in the naive sense (see Introduction). 
By Proposition \ref{compactpro}, there exists $G\in\grd$ such that $G_n\tart G$.
\end{remark}
\noindent
\subsection{Schreier graphs} \label{schreier} 
In this subsection we define naive convergence for Schreier graphs.
Let $\Gamma$ be a finitely generated group and $\Sigma$ be a finite, symmetric generating set for $\Gamma$.
Let $H\subset \Gamma$ be a subgroup of $\Gamma$. Recall that the Schreier graph $S(\Gamma, H)$ is defined in the following
way.
\begin{itemize}
\item $V(S(\Gamma,H))=\Gamma/H$, the set of right cosets of $H$,
\item there is a directed edge from $Ha$ to $Hb$ labeled by $\sigma\in\Sigma$ if $Ha\sigma=Hb$.
\end{itemize}
\noindent
The root of $S(\Gamma, H)$ is the subgroup $H$ itself. Note that by definition any Schreier graph is a rooted, connected
graph. The set of all Schreier graphs are denoted by $\gsg$. Again,
we can define a metric $d_{\gsg}$ of the set $\gsg$.
Let $S(\Gamma,H), S(\Gamma,K)\in\gsg$. Then,
$$d_{\gsg}(S(\Gamma,H), S(\Gamma,K))=2^{-n}\,,$$
\noindent
where $n$ is the largest integer for which the balls $B_n(S(\Gamma,H),H)$ and \\ $B_n(S(\Gamma,K),K)$ are
rooted $\Sigma$-edge labeled isomorphic. Again, $\gsg$ is a compact metric space with respect to the 
metric $d_{\gsg}$.
Note that the map $H\to S(\Gamma,H)$ provides a homeomorphism $\tau:\Sub(\Gamma)\to \gsg$ (\cite{Elekurs}).
One can also define the compact space of Cantor-vertex labeled, rooted Schreier graphs $\cC\gsg$ and
the compact space of $Q$-labeled rooted Schreier graphs in a similar fashion.
Observe that we have a natural $\Gamma$-action
$\alpha:\Gamma\actson \gsg$, where
$\alpha(\gamma)(S(\Gamma,H))=S(\Gamma,g H g^{-1})$.
The action extends to the spaces $\cC \gsg$ and $\gsg^Q$ as well 
(\cite{Elekmin}).
If $\phi:\Gamma/H\to \cC$ is a map, then
$$\alpha(g)(S(\Gamma,H),\phi)=(S(\Gamma, gHg^{-1}),\psi)\,,$$
\noindent
where $\psi(gHg^{-1} a)=\phi(Hag)$.

\noindent
A {\bf generalized Schreier graph $S$} is the countable union of some Schreier graphs. Generalized Schreier graphs are
associated to not necessarily transitive actions of $\Gamma$ on countable sets. We denote the space of all
generalized Schreier graphs of $\Gamma$ with respect to $\Sigma$ by $G\Gamma_\Sigma \cG$. Let $U^r_{\Gamma,\Sigma}$ be the set of rooted,
$\Sigma$-edge labeled balls of radius $r$ (up to isomorphisms) that occur in some Schreier graph of $\Gamma$.
For $S\in G\gsg$, we denote by $\cB_{\Gamma,\Sigma}(S)$ the set of all rooted, $\Sigma$-edge labeled balls $B$ for which there
exists $x\in V(S)$ and $k\geq 1$ such that $B_k(S,x)$ is rooted-labeled isomorpic to $B$.
Let $S,T\in G\gsg$. We say that $S$ and $T$ are equivalent if  $\cB_{\Gamma,\Sigma}(S)= \cB_{\Gamma,\Sigma}(T)$. Again, we denote by
$\overline{G\gsg}$ the set of equivalence classes of generalized Schreier graphs.
We can define a metric on $\overline{G\gsg}$ in the same way as we did for $\overline{\grd}$. We can also consider
the set of Cantor labeled generalized Schreier graphs $\cC G \gsg$ and the set $G\gsg^Q$ as well, together with
the metric spaces  $\overline{\cC G\gsg}$ and $\overline{G\gsg^Q}$. Then, we have the following straightforward generalization of
Proposition \ref{compactpro} and Proposition \ref{labeledversion}.
\begin{proposition}
The metric spaces $\overline{G\gsg}$,  $\overline{\cC G\gsg}$ and $\overline{G\gsg^Q}$ are compact.
Also, $\overline{O}:\overline{G\gsg}\to Inv(\gsg)$, $\overline{O}:\overline{\cC G\gsg}\to Inv(\cC \gsg)$ and 
$\overline{O}:\overline{G\gsg^Q}\to Inv(\gsg^Q)$ are
homeomorphisms.
\end{proposition}

\subsection{Infinite graphs that cannot be approximated by finite graphs}
It is known that all graphons are limits of a convergent sequence
of dense finite graphs \cite{LSZ}. On the other hand, it is
not known whether every measured graphing is the limit of sparse graphs.
For naive convergence, we have the following proposition.
\begin{proposition}
There exists a countable, connected infinite graph of \\bounded vertex degree, which is
not the naive limit of a sequence of finite graphs.
\end{proposition}
\proof Let us recall the notion of a LEF-group \cite{GV}.
Let $\Gamma$ be a finitely generated group with a symmetric generating
set $\Sigma$. Let $\Cay(\Gamma,\Sigma)$ be the corresponding right
Cayley graph.That is,
\begin{itemize}
\item
$V(\Cay(\Gamma,\Sigma))=\Gamma$,
\item
$\overrightarrow{(\gamma,\delta)}$ is a directed edge of
$\Cay(\Gamma,\Sigma)$ labeled by $\sigma\in\Sigma$, if $\sigma\gamma=\delta$.
\end{itemize}
\noindent
A $\Sigma$-graph is a finite graph $G$ such that
\begin{itemize}
\item for each vertex $x\in V(G)$, $\deg(x)=|\Sigma|$,
\item and each edge $\overrightarrow{(x,y)}$ is labeled uniquely
by some element of $\Sigma$, where
\item the label of $\overrightarrow{(x,y)}$ is the inverse
of the label of $\overrightarrow{(y,x)}$.
\end{itemize}
\noindent
We say that a $\Sigma$-graph $G$ is an $n$-approximant of
$\Cay(\Gamma,\Sigma)$ if all the $n$-balls in $G$ are edge-labeled-isomorphic
to the $n$-ball of the Cayley graph $\Cay(\Gamma,\Sigma)$.
The group $\Gamma$ is a LEF-group if for any $n\geq 1$ $\Cay(\Gamma,\Sigma)$ 
 possesses
an $n$-approximant (it is not hard to see that being LEF is independent
of the choice of the generating system).
Gordon and Vershik \cite{GV} showed that not all groups $\Gamma$ are LEF.
Now, let $\Gamma$ be
a finitely generated group and let $\Sigma=\{\sigma_i\}^k_{i=1}$
be a symmetric generating system of $\Gamma$.
We encode $\Cay(\Gamma,\Sigma)$ by an undirected, unlabeled graph in
the following way.
For any $1\leq i \leq k$, connect $x\in\Gamma$ with $\sigma_i x$ by
a path of length $3k$. So, from any $x\in\Gamma$ there are $k$ ongoing
path. If the path corresponding to the element $\sigma_i$, then
glue a hanging edge to the $i$-th vertex of the path. Notice that
the path between $x$ and $\sigma_i x$ will receive two hanging edges, one
from the direction of $x$ and one from the direction of $\sigma_i x$.
Let $G$ be the resulting infinite graph. Clearly, there exists
a sequence of finite graphs
$\{G_n\}^\infty_{n=1}\subset Gr_d$, $G_n\stackrel{n}{\rightarrow} G$
if and only if $\Gamma$ is a LEF-group. Hence our proposition follows
from the existence of a non-LEF group. \qed

\section{Qualitative weak equivalence of free Cantor actions}
The goal of this section is to study the qualitative analogue of the weak equivalence of essentially free probability measure preserving actions (see
\cite{Burton} for a survey).
The section also serves as a preparation for our qualitative graph limit theory.
\subsection{Weak containment and weak equivalence}
Let $\Gamma$ be a finitely generated group and $\Sigma$ be a finite, symmetric generating 
set for $\Gamma$. Let $r>0$ be an integer and $\phi:\cC\to Q$ be a continuous map, where $Q$ is
a finite set and $\alpha:\Gamma\actson \cC$ is a free Cantor-action.  For $x\in \cC$ let the map $\tau^r_{x,\alpha,\phi}:B_r(\Gamma,\Sigma,e_\Gamma)\to Q$ defined by setting
$$\tau^r_{x,\alpha,\phi}(\gamma):=\phi(\alpha(\gamma)(x))\,.$$
\noindent
Also, \begin{enumerate}
\item Let $\conf_{r,\alpha}(\phi)=\cup_{x\in \cC} \tau^r_{x,\alpha,\phi}\,.$
\item If  $Q=\{0,1\}^k$ we set $\conf_{r,k}(\alpha)=\cup_{\phi:\cC\to \{0,1\}^k} \conf_{r,\alpha}(\phi)\,,$
\item Let $\conf(\alpha)=\cup_{r,k} \conf_{r,k}(\alpha).$
\end{enumerate}
\begin{definition}
Let $\alpha:\Gamma\actson \cC$ and $\beta:\Gamma\actson \cC$ be free Cantor-actions. We
say that $\alpha$ {\bf qualitatively weakly contains} $\beta$, $\alpha\succeq \beta$,
if for any finite set $Q$, continuous map $\psi:\cC\to Q$ and $r>0$,
there exists a continuous map $\phi:\cC\to Q$ such that
$$ \conf_{r,\alpha}(\phi)=\conf_{r,\beta}(\psi)\,.$$
\noindent
We say that $\alpha$ is {\bf qualitatively weakly equivalent} to $\beta$, $\alpha \simeq \beta$, if
$\alpha\succeq \beta$ and $\beta\succeq\alpha$. 
\end{definition}
\noindent
Clearly, qualitative weak containment does not depend on the choice of the generating system $\Sigma$.
The set of all quantitative weak equivalent classes will be denoted by $\Free(\Gamma)$.
\begin{remark} In case of free $\Z$-actions, our notion of qualitative weak equivalence coincide with the notion
 of weak approximate conjugacy introduced by Lin and Matui \cite{LM}.
\end{remark}
\noindent
Recall that the cost and strong ergodicity are invariants with respect to weak equivalence in the measurable
framework \cite{Burton}. The next theorem shows that certain important properties of Cantor actions are, in fact, invariants
of the qualitative weak equivalence class.
\begin{theorem}
Let  $\alpha:\Gamma\actson \cC$ and  $\beta:\Gamma\actson \cC$ be Cantor actions. 
\begin{enumerate}
\item If $\alpha$ admits an invariant probability measure, then $\beta$ admits an invariant probability measure
as well.
\item If $\beta$ is an amenable action, then $\alpha$ is an amenable action as well.
\end{enumerate}
\noindent
Hence, both amenability and admitting an invariant probability measure are qualitative weak invariants.
\end{theorem}
\proof
Let $\{\{0,1\}^k\}^{B_r(\Gamma,\Sigma,e_\Gamma)}=\tbkr$ be
the set of all $\{0,1\}^k$-valued functions on the ball $B_r(\Gamma,\Sigma,e_\Gamma) $ of radius
$r$ centered around the unit element of $\Gamma$ in the Cayley graph of $\Gamma$ with respect
to $\Sigma$. If $A\in \tbkr$ then we can define the clopen set 
$$U_A:=\{x\in\cC\,\mid\,(\beta(\gamma)(x))_{[k]}=
A(\gamma)\,\mbox{for any}\, \gamma\in B_r(\Gamma,\Sigma,e_\Gamma)\}\,,$$
\noindent
where $(x)_{[k]}$ is the first $k$ coordinates of $x$.
Clearly, $\cup_{A\in \tbkr} U_A$ is a clopen partition of $\cC$.
For $\sigma\in\Sigma$ and $A,B\in\tbkr$, we define two further clopen sets
$$U^{\sigma}_{A\to B}:=\{x\in \cC\,\mid\, x\in U_A, \beta(\sigma)(x)\in U_B\}\,.$$
\noindent
and
$$U^{\sigma}_{B\leftarrow A}:=\{y\in \cC\,\mid\, y\in U_B, \beta(\sigma^{-1})(y)\in U_A\}\,.$$
\noindent
Let $r,k\geq 1$ and $\phi:\cC\to \{0,1\}^h$ be a continuous map, where $h\geq k$.
For $A\in \tbkr$ we set
$$V_A^\phi=\{x\in \cC\,\mid\, \phi_k(\alpha(\gamma))(x)=A(\gamma)\,\mbox{for any}\,\gamma
\in B_r(\Gamma,\Sigma,e_\Gamma)\}\,.$$
\noindent
(recall that $\phi_k(x)=(\phi(x))_{[k]}$) 
Also, if $A\in \tbkr$,  where $1\leq l\leq k$,  let
$A_{[l]}\in \tblr$ is the projection of the values of $A$ onto the first $l$ coordinates.
If $1\leq s\leq r$, let $A_{\mid s}\in \tbks$ is the restriction of $A$ onto the
ball $B_s(\Gamma,\Sigma,e_\Gamma)$. Finally, we use the notation $A\sqsubset B$ if $B=(A_{\mid s})_{[l]}$.
Note that $A\sqsubset B$ implies that $U_A\subset U_B$.
So, if  $B\in\tbls$, then
$$ V^\phi_B=\cup_{A\in\tbkr,A\sqsubset B} V^\phi_A\,.$$
\noindent
Since $\alpha\succeq\beta$, for all $r\geq 1$ there exists a continuous map  $\phi^r:\cC\to \{0,1\}^r$ such that
for all  $A\in \tbr$, $V^{\phi^r}_A$ is nonempty if and only if $U_A$ is nonempty.
Let $\{r_n\}_{n=1}^\infty$ be an increasing sequence of integers such that
for all $s,l\geq 1$ and $B\in \tbls$,
$$\lim_{n\to\infty} \mu(V_B^{\phi^{r_n}})=l(B)$$
\noindent
exists, where $\mu$ is a probability measure on $\cC$ invariant under the action $\alpha$.
Set $\nu(U_B)=l(B)$.
\begin{lemma}\label{lf5}
The function $\nu$ extends to a Borel probability measure on $\cC$ in a unique way.
\end{lemma}
\proof
If $B\in \tbls$, $s\leq r$, $l\leq k$ and $n$ is a large enough integer, then
$\mu(V_B^{\phi^{r_n}})$ is well-defined and equals to $\sum_{A\in\tbkr,A\sqsubset B}V^{\phi^{r_n}}_A$
and also
$$\sum_{B\in \tbls}\mu (V^{\phi^{r_n}}_B)=1\,.$$
\noindent
Also, if for some $B\in\tbls$, the set $U_B$ is empty, then $\mu(V_B^{\phi^{r_n}})=0.$
Hence, if $1\leq l \leq k$ and $1\leq s\leq r$, then 
$$\nu(U_B)=\sum_{A\in\tbkr,A\sqsubset B} \nu (U_A)\,,$$
\noindent
and
$$\sum_{B\in \tbls}\nu (U_B)=1\,.$$
\noindent
Also, $\nu(U_B)$ is defined being zero, if $U_B$ is the empty set.
Hence, $\nu$ is a premeasure on the basic sets of $\cC$, therefore, $\nu$ extends to a Borel probability
measure in a unique way. \qed
\begin{proposition}\label{pf6}
The measure $\nu$ is invariant under the action $\beta$.
\end{proposition}
\proof
It is enough to show that for any $r\geq 1$, $\sigma\in\Sigma$ and pair $A,B\in \tbr$, we have that
$$ \nu( U^{\sigma}_{A\rightarrow B})=  \nu( U^{\sigma}_{B\leftarrow A})\,.$$
\noindent
For large enough $n$, we can define the clopen sets $V^{\sigma,\phi^{r_n}}_{A\rightarrow B}$ and  $V^{\sigma,\phi^{r_n}}_{B\leftarrow A}$ by
$$V^{\sigma,\phi^{r_n}}_{A\rightarrow B}:=\{x\in V_A^{\phi^{r_n}}\,\mid\, \alpha(\sigma)(x)\in V_B^{\phi^{r_n}}\}\,,$$
$$V^{\sigma,\phi^{r_n}}_{A\leftarrow B}:=\{z\in V_B^{\phi^{r_n}}\,\mid\, \alpha(\sigma^{-1})(z)\in V_A^{\phi^{r_n}}\}\,.$$
\begin{lemma}\label{lf6}
$$\nu( U^{\sigma}_{A\rightarrow B})= \lim_{n\to\infty} \mu(V_{A\rightarrow B}^{\phi^{r_n}})\,, $$
$$\nu( U^{\sigma}_{B\leftarrow A})= \lim_{n\to\infty} \mu(V_{B\leftarrow A}^{\phi^{r_n}})\,. $$    
\end{lemma}
\proof
For $D\in\tbrp$, let us use the notation $D\sqsubset A\to B$ if
$U_D \subset  U^{\sigma}_{A\rightarrow B}$.
Also, for $E\in\tbrp$, let $E\sqsubset B\leftarrow A$ if
$U_E \subset  U^{\sigma}_{B\leftarrow A}$.
Then we have that
\begin{itemize}
\item $U^{\sigma}_{A\rightarrow B}=\cup_{D\in\tbrp\,,D\sqsubset A\to B} U_D\,,$
\item $U^{\sigma}_{B\leftarrow A}=\cup_{E\in\tbrp\,,E\sqsubset B\leftarrow A } U_\cC\,,$
\item for large enough  $n$, 
$$V_{A\rightarrow B}^{\phi^{r_n}}=\cup_{D\in\tbrp\,,D\sqsubset A\to B} V_{D}^{\phi^{r_n}}\,,$$
\item for large enough $n$, 
$$V_{B\leftarrow A}^{\phi^{r_n}}=\cup_{E\in\tbrp\,,E\sqsubset B\leftarrow A} V_{E}^{\phi^{r_n}}\,.$$ 
\end{itemize}
\noindent
Therefore our lemma follows. \qed
\vskip 0.1in
\noindent
Since $\mu$ is invariant under the action $\alpha$, for all $n\geq 1$, we have that
$$\mu(V_{A\rightarrow B}^{\phi^{r_n}})= \mu(V_{B\leftarrow A}^{\phi^{r_n}})\,.$$
\noindent
Therefore,
$\nu(U^{\sigma}_{A\rightarrow B})=\nu(U^{\sigma}_{B\leftarrow A})\,.$ Thus our proposition follows. \qed
\vskip 0.1in
\noindent
Now, the first part of our theorem follows from Lemma \ref{lf5} and Proposition \ref{pf6}. 
Let us turn to proof of the second part of our theorem.
First, let us recall the notion of an amenable action from \cite{HR}.
The free action $\beta:\Gamma\actson\cC$ is {\bf an amenable action} if
there exists a sequence of weak$^*$-continuous maps $b_n:X\to \Prob(\Gamma)$ such that
for every $\gamma\in \Gamma$,
\begin{equation}\label{ef8}
\lim_{n\to\infty} \sup_{x\in\cC} \| S(\gamma)(b_n(x))-b_n(\beta(\gamma)(x))\|_1=0\,,
\end{equation}
\noindent
where $S$ is the natural action of the group $\Gamma$ on $\Prob(\Gamma)$.
Since $\Gamma$ is finitely generated, it is enough to assume \eqref{ef8} for the generators $\sigma\in\Sigma$.
\begin{lemma}\label{lf17}
Let $b:\cC\to\Prob(\Gamma)$ be a weak$^*$-continuous function. Then for any $\eps>0$, there exists $R>0$ 
and a  weak$^*$-continuous function $b':\cC\to\Prob(\Gamma)$ such that for all $x\in\cC$,
$\supp(b'(x))\subseteq B_R(\Gamma,\Sigma,e_\Gamma)$ and \\ $\|b(x)-b'(x)\|_1<\eps.$
\end{lemma}
\proof
Suppose that there exists a sequence $\{x_n\}^\infty_{n=1}$ such that for all $n\geq 1$
$$ b(x_n)(B_n(\Gamma,\Sigma,e_\Gamma))\leq 1-\frac{\eps}{3}\,.$$
\noindent
Then, for any limit point $x$ of the sequence, $b(x)(\Gamma)\leq 1-\frac{\eps}{3}$ leading to a contradiction.
That is, there must exist some $R>0$ , such that for all $x\in\cC$,
$$b(x)(B_R(\Gamma,\Sigma,e_\Gamma))>1-\frac{\eps}{3}\,.$$
\noindent
We define the function $b'$ in the following way.
\begin{itemize}
\item Let $b'(x)(\gamma)=b(x)(\gamma)$, if $e_\Gamma\neq \gamma \in B_R(\Gamma,\Sigma,e_\Gamma)$.
\item Let $b'(x)(\gamma)=0$, if $\gamma\notin B_R(\Gamma,\Sigma,e_\Gamma)$.
\item Let  $b'(x)(e_\Gamma)=b(x)(e_\Gamma)+b(x)(\Gamma\backslash B_R(\Gamma,\Sigma,e_\Gamma)$.
\end{itemize}
\noindent
Clearly, $b':\cC\to\Prob(\Gamma)$ is  weak$^*$-continuous, for any $x\in\cC$ we have that $\supp(b'(x))\subseteq
 B_R(\Gamma,\Sigma,e_\Gamma)$ and $\|b(x)-b'(x)\|_1<\eps.$ \qed
\vskip 0.1in
\noindent
So, from now on we can assume that for any $n\geq 1$, there exists some $R_n>0$ such that for all $x\in\cC$,
\begin{equation}\label{ef18}
\supp(b_n(x))\subseteq B_{R_n}(\Gamma,\Sigma,e_\Gamma).
\end{equation}
\noindent
The following lemma is trivial.
\begin{lemma}\label{1lf18}
Let $f\in \Prob(\Gamma)$ such that $\supp(f)\subseteq B_R(\Gamma,\Sigma,e_\Gamma)$
for some $R>0$. Let $k\geq 1$ be an integer. Then, there exists $g\in \Prob(\Gamma)$,
$\supp(g)\subseteq B_R(\Gamma,\Sigma,e_\Gamma)$, such that
\begin{itemize}
\item for any $\gamma\in  B_R(\Gamma,\Sigma,e_\Gamma)$, $g(\gamma)=\frac{i}{k}$ for some integer
$i\geq 0$,
\item $\|f-g\|_1=\sum_{\gamma\in B_R(\Gamma,\Sigma,e_\Gamma)} |f(\gamma)-g(\gamma)|\leq \frac{| B_R(\Gamma,\Sigma,e_\Gamma)|}
{k}\,.$
\end{itemize}
\end{lemma}
\begin{lemma} \label{2lf18}
Let $\{b_n\}^\infty_{n=1}$ as in \eqref{ef18}. Then, for any $n\geq 1$ and $\eps>0$, there exists $k>1$ and 
continuous function $c_n:\cC\to\Prob(\Gamma)$ such that
for all $\cC$,
\begin{itemize}
\item $\supp(c_n(x))\subseteq B_{R_n}(\Gamma,\Sigma,e_\Gamma),$
\item for all $\gamma\in B_{R_n}(\Gamma,\Sigma,e_\Gamma), c_n(x)(\gamma)=\frac{i}{k}$ for some integer
$i\geq 0$,
\item $\|b_n(x)-c_n(x)\|_1\leq \eps.$
\end{itemize}
\end{lemma}
\proof
Let $\{U_\alpha\}_{\alpha\in I}$ be a finite clopen partition of the Cantor set $\cC$ such that
if $x,y\in U_\alpha$ for some $\alpha$, then $\|b_n(x)-b_n(y)\|_1\leq \eps/3.$ Let $k\geq 1$ be an integer
such that $$\frac{|B_{R_n}(\Gamma,\Sigma,e_\Gamma)|}{k}< \frac{\eps}{3}\,.$$
\noindent

For each $\alpha$  pick
an element $x_\alpha\in U_\alpha $. Then for each $n\geq 1$ choose  $c^\alpha_n \in \Prob(\Gamma)$ in such a way that
\begin{itemize}
\item $\|c^\alpha_n-b_n(x_\alpha)\|\leq \frac{\eps}{3}\,,$
\item $\supp(c^\alpha_n)\subseteq B_{R_n}(\Gamma,\Sigma,e_\Gamma)\,,$
\item for any $\gamma\in B_{R_n}(\Gamma,\Sigma,e_\Gamma), c^\alpha_n(\gamma)=\frac{i}{k}$, where
$i\geq 0$ is an integer.
\end{itemize}
\noindent
Finally, define $c_n:\cC\to \Prob(\Gamma)$ by setting $c_n(y):=c^\alpha_n$ if $y\in U_\alpha$.
Then the functions $\{c_n\}^\infty_{n=1}$ satisfy the conditions of our lemma. \qed
\vskip 0.1in
\noindent
Now, let $\alpha,\beta:\Gamma\actson\cC$ be free Cantor actions such that
$\alpha\succeq\beta$ and $\beta$ is amenable.
Let $b:\cC\to\Prob(\Gamma)$ be a weak$^*$-continuous map, $k\geq 1$ and $R\geq 1$ be integers, $\eps>0$
such that for all $x\in\cC$ and $\sigma\in\Sigma$.
\begin{itemize}
\item $\supp(b(x))\subseteq  B_R(\Gamma,\Sigma,e_\Gamma)\,,$
\item  for any $\gamma\in B_R(\Gamma,\Sigma,e_\Gamma), b(x)(\gamma)=\frac{i}{k}$, where
$i\geq 0$ is an integer,
\item $\|S(\sigma)(b(x))-b(\beta(\sigma)(x)\|_1 <\eps.$
\end{itemize}
\noindent
Let $Q$ be the finite set of elements $c$ in $\Prob(\Gamma)$, \\such
that $\supp(c)\subseteq  B_R(\Gamma,\Sigma,e_\Gamma)$ and $c(\gamma)=\frac{i}{k}$, where
$i\geq 0$ is an integer. Define the continuous map $\pi:\cC\to Q$ by setting
$\pi(x):=b(x)\in Q$. Let $\hat{\pi}:\cC\to Q$ be a  continuous map  such that
$$\conf_{1,\alpha}(\hat{\pi})=\conf_{1,\beta}(\pi)\,.$$
\noindent
Now, let $\hat{b}:\cC\to\Prob(\Gamma)$ be defined by
$\hat{b}(y)=\hat{\pi}(y).$
Then, for any $y\in\cC$ and $\sigma\in\Sigma$,
$$\| S(\sigma)(\hat{b}(y))-\hat{b}(\alpha(\sigma)(y))\|_1<\eps\,.$$
\noindent
Hence, $\alpha$ is amenable as well. This finishes the proof of our theorem. \qed
\subsection{Cantor subshifts}
\begin{definition}
Let $Z\subset \cC^\Gamma$ be a closed subset invariant under the right shift action, $(R(\delta)(z))(\gamma)=z(\gamma\delta)$.
 We call $Z$
a {\bf Cantor subshift} if $Z$ is homeomorphic to $\cC$ and the action on $Z$ is free.
\end{definition}
\noindent
Let $\alpha:\Gamma\actson\cC$ be a free Cantor action. Let the equivariant map $\kappa_\alpha:\cC\to \cC^\Gamma$ be 
defined by setting
$\kappa_\alpha(x)(\gamma):=\alpha(\gamma)(x).$  We call $\kappa_\alpha(\cC)$ the Cantor subshift of $\alpha$.
Clearly, every action is conjugate isomorphic (hence, qualitatively weakly equivalent) to its own Cantor subshift.
Note that not all Cantor subshift $Z\subset \cC^\Gamma$ are in the form of  $\kappa_\alpha(\cC)$ for some action
$\alpha$. Let $Y\subset \cC^\Gamma$ be a Cantor subshift and $r\geq 1$. Then 
$\conf_{r,Y}([t])\subseteq \{\{0,1\}^t\}^{B_r(\Gamma,\Sigma,e_\Gamma)}$ is defined in
the following way. The map $\phi:B_r(\Gamma,\Sigma,e_\Gamma)\to \{0,1\}^t$ is an element of $\conf_{r,Y}([t])$
if there exists $y\in Y$ such that for all  $\gamma\in B_r(\Gamma,\Sigma,e_\Gamma)$, $\phi(\gamma)=(y(\gamma))_{[t]}.$
The following two lemmas will  be used in our proofs.
\begin{lemma}\label{lf14}
Let  $Y\subset \cC^\Gamma$ be a Cantor subshift and $\alpha:\Gamma\actson\cC$ be a free Cantor action. Then
$\alpha\succeq Y$ if and only for any $r,t\geq 1$, there exists a continuous map $\phi:\cC\to \{0,1\}^t$ such that
$\conf_{r,\alpha}(\phi)=\conf_{r,Y}([t])$.
\end{lemma}
\proof
Clearly, the condition of the lemma is necessary for $Y$ begin qualitatively weakly contained by $\alpha$.
Let us show that the condition is sufficient as well. Let $Q$ be a finite set, $r\geq 1$ and $\psi:Y\to Q$
be a continuous map.
We need to show that there exists $\hat{\psi}:\cC\to Q$ such that
$$\conf_{r,\alpha}(\hat{\psi})=\conf_{r,Y}(\psi).$$
\noindent
For $A\in \conf_{n,Y}{[t]}$ set
$$W_A:=\{y\in Y\,\mid\, y(\gamma)_{[t]}=A(\gamma),\,\mbox{for any}\, \gamma\in  B_n(\Gamma,\Sigma,e_\Gamma)\}\,.$$
\noindent
Since $\psi$ is continuous, there exists $n,t\geq 1$ such that $\psi$ is constant on the clopen
sets $W_A$, for all $A\in\conf_{n,Y}{[t]}$.
Let $\pi:\cC\to\{0,1\}^t$ such that 
$$\conf_{n+r,\alpha}(\pi)=\conf_{n+r,Y}([t])\,.$$
\noindent
Again, if $A\in\conf_{n,\alpha}(\pi)=\conf_{n,Y}([t]),$
let
$$V_A:=\{x\in \cC\,\mid\, \pi(\alpha(\gamma))(x))=A(\gamma),\,\mbox{for any}\, \gamma\in 
 B_n(\Gamma,\Sigma,e_\Gamma)\}\,.$$
\noindent
Now, we define $\hat{\psi}:\cC\to Q$ in the following way.
If $x\in V_A$, then $\hat{\psi}(x)=q$ where $q$ is the value of $\psi$ on the clopen set $W_A$.
Then, it is easy to see that
$$\conf_{r,\alpha}(\hat{\psi})=\conf_{r,Y}(\psi)\,,$$
\noindent
hence the lemma follows. \qed
\vskip 0.1in
\noindent
For $z\in\{0,1\}^t$ and $x\in\cC=\{0,1\}^\N$, let $x\downarrow z\in \cC$ be defined in the following way.
\begin{itemize}
\item If $1\leq i \leq t$, $x\downarrow z(i)=z(i),$
\item if $t+1\leq i$, $x\downarrow z(i)=x(i-t)\,.$
\end{itemize}
\noindent
Now, let $\alpha:\Gamma\actson\cC$ be a free Cantor action and $\pi:\cC\to \{0,1\}^t$ be 
a continuous map.
Then we define an equivariant homeomorphism $\kappa^\pi_\alpha:\cC\to\cC^\Gamma$ by setting
$$\kappa^\pi_\alpha(x)=\alpha(\gamma)(x)\downarrow \pi(\alpha(\gamma)(x))\,.$$
\noindent
The following lemma is trivial.
\begin{lemma} \label{lf16}
The space $\kappa^\pi_\alpha(\cC)$ is a Cantor subshift.
\end{lemma}

\subsection{The smallest element}
Our next goal is to prove a qualitative analogue of a result of Ab\'ert and Weiss \cite{AW}. Their theorem
states that all essentially free p.m.p. action of a countable group weakly contains the Bernoulli actions.
\begin{theorem} \label{AWe}
Let $\Gamma$ be a finitely generated group. Then, there exists a free Cantor action $\beta:\Gamma\actson\cC$
such that any free Cantor action $\alpha:\Gamma\actson\cC$ qualitatively weakly contains $\beta$.
\end{theorem}
\proof
First, we explicitly construct the action $\beta$. For $n\geq 1$, let $r_n$ be a positive integer such
that $2^{r_n}\geq | B_n(\Gamma,\Sigma,e_\Gamma)|$. So, one can label the elements of $\Gamma$ by the set
$\{0,1\}^{r_n}$ in such a way that if $0<d_{Cay(\Gamma,\Sigma)}(x,y)\leq n$, then the labels of $x$ and $y$
are different.
\begin{definition}
An {\bf $n$-block } is a labeling $$\phi:\prod^n_{j=1} \phi^j: B_n(\Gamma,\Sigma,e_\Gamma)\to 
\prod^n_{j=1}\{0,1\}^{r_j}\,,$$
\noindent
where $\phi^j: B_n(\Gamma,\Sigma,e_\Gamma)\to \{0,1\}^{r_j}$, such that
if $1\leq j \leq n$ and \\ $0< d_{Cay(\Gamma,\Sigma)}(x,y)\leq j$, then $\phi^j(x)\neq \phi^j(y)$.
\end{definition}
\noindent
Now, let $y=\prod^\infty_{j=1} y_j: \Gamma\to \prod^\infty_{j=1} \{0,1\}^{r_j}$
be a labeling such that
\begin{itemize}
\item for any $j\geq 1$, if $0< d_{Cay(\Gamma,\Sigma)}(x,y)\leq j$, then $y_j(x)\neq y_j(y)$,
\item for any $n\geq 1$ and $n$-block
$\phi: B_n(\Gamma,\Sigma,e_\Gamma)\to 
\prod^n_{j=1}\{0,1\}^{r_j}$, there exists $\delta\in\Gamma$ such that for any $1\leq j\leq n$ and $\gamma\in  B_n(\Gamma,\Sigma,e_\Gamma)$
we have 
$\phi^j(\gamma)=y_j(\delta\gamma)$. 
\end{itemize}
\noindent
We call such a labeling $y$ a full labeling. Clearly, full labelings exist.
\begin{lemma} \label{lf10}
Let $y$ be a full labeling and $Y\subset \cC^\Gamma$ be the orbit closure of $y$
in $\cC^\Gamma$, where $\cC=\{0,1\}^\N$ is identified with $\prod^\infty_{j=1} \{0,1\}^{r_j}$.
Then the restricted shift action $\beta:\Gamma\actson Y$ is free.
\end{lemma}
\proof
Observe that for any $r\geq 1$, there exists $s_r\geq 1$ such that if $0< d_{Cay(\Gamma,\Sigma)}(\gamma,\delta)\leq r$
then $(y(\gamma))_{[s_r]}\neq (y(\delta))_{[s_r]}$.
Hence, if $\hat{y}\in Y$ and $0< d_{Cay(\Gamma,\Sigma)}(\gamma,\delta)\leq r$,
then $(\hat{y}(\gamma))_{[s_r]}\neq (\hat{y}(\delta))_{[s_r]}$. That is,
$\beta(\gamma)(\hat{y})\neq \hat{y}$ provided that $\gamma\neq e_\Gamma$. Hence the action $\beta$ is free on $Y$.
\qed
\begin{lemma}\label{lf11}
The space $Y$ is homeomorphic to the Cantor set.
\end{lemma}
\noindent
\proof
It is enough to show that $Y$ does not have an isolated point.
Let $\hat{y}\in Y$. Observe that if we restrict $\hat{y}=\prod^\infty_{j=1}\hat{y}_j$ onto the ball
$ B_n(\Gamma,\Sigma,e_\Gamma)$ we obtain an $n$-block $\phi$. On the other hand, there are more than one $n+1$-blocks
$\psi$ such that if we restrict $\psi$ onto $ B_n(\Gamma,\Sigma,e_\Gamma)$ taking only the first $\sum_{j=1}^nr_j$
coordinates, the resulting $n$-block coincides with $\phi$. Hence, for any $\eps>0$, there exists $z\in Y$ such
that $0<d_Y(y,\hat{y})<\eps.$ Therefore, the space $Y$ is homeomorphic to $\cC$. \qed
\vskip 0.1in
\noindent
In order to finish the proof of our theorem. It is enough to prove the following proposition.
\begin{proposition} \label{pf12}
Let $\alpha:\Gamma\actson \cC$ be a free Cantor action. Then, $\alpha\succeq \beta$.
\end{proposition}
\proof
Let $\{r_j\}^\infty_{j=1}$ be as above and set $k_n:=\sum^n_{j=1} r_j$. Our goal is to construct
a continuous map $\pi_n:\cC\to \{0,1\}^{k_n}$ such that the set $\conf_{n,\alpha}(\pi_n)$ equals
to the set of all $n$-blocks. Then, we have that
$$\conf_{n,k_n}(\kappa^{\pi_n}_\alpha)=\conf_{n,k_n}(\beta)\,.$$
\noindent
Then by Lemma \ref{lf14}, it follows that $\alpha\succeq\beta$.
First, we make sure that $\conf_{n,\alpha}(\pi_n)$ contains all the $n$-blocks.
We use an insertion process to define $\pi_n$ partially.
Let $M_n$ be a positive integer, which is so large that \\ $\{\{0,1\}^{k_n}\}^{B_{M_n}(\Gamma,\Sigma,e_\Gamma)}$
contains an element
$$\psi^n=\prod^n_{j=1} \psi^n_j: B_{M_n}(\Gamma,\Sigma,e_\Gamma)\to \{0,1\}^{k_n}$$
\noindent
such that
\begin{itemize}
\item for any $1\leq j \leq n$, $\psi^n_j(\gamma_1)\neq\psi^n_j(\gamma_2)$,
whenever $0<d_{Cay(\Gamma,\Sigma)}(\gamma_1,\gamma_2)\leq j$,
\item for any $n$-block $\phi:B_n(\Gamma,\Sigma,e_\Gamma)\to\{0,1\}^{k_n}$,
there exists \\$\zeta\in B_{M_n}(\Gamma,\Sigma,e_\Gamma)\to \{0,1\}^{k_n}$ such
that $B_n(\Gamma,\Sigma,\zeta)$ is contained in the ball $B_{M_n}(\Gamma,\Sigma,e_\Gamma)$
and $\psi^n(\zeta\gamma)=\phi(\gamma)$ for any $\gamma\in B_n(\Gamma,\Sigma,\zeta)$.
\end{itemize}
\noindent
Informally speaking, the labeling $\psi^n$ contains an inserted copy of all the $n$-blocks.
Now let $U_n$ be a nonempty clopen subset of $\cC$ such that if $x\neq y\in U_n$ then
$$ \bigcup_{\gamma\in B_{2M_n}(\Gamma,\Sigma,e_\Gamma)} \alpha(\gamma)(x) \cap
\bigcup_{\gamma\in B_{2M_n}(\Gamma,\Sigma,e_\Gamma)} \alpha(\gamma)(y) =\emptyset\,.$$
\noindent
We define the map $\pi_n$ on the clopen set $\bigcup_{x\in U_n} \bigcup_{B_{M_n}(\Gamma,\Sigma,e_\Gamma)}
\alpha(\gamma)(x)$ in the following way. If $y=\alpha(\gamma)(x)$, where $x\in U_n$ and
$\gamma\in B_{M_n}(\Gamma,\Sigma,e_\Gamma)$, then let $\pi_n(y)=\psi^n(\gamma)$.
It follows from our construction that $\conf_{n,\alpha} (\pi_n)$ will contain all the $n$-blocks.
We only need to make sure that $\conf_{n,\alpha} (\pi_n)$ contains only $n$-blocks.
Thus, our proposition follows from the lemma below.
\begin{lemma}\label{lf23}
Let $r\geq 1$ be a positive integer and $Q$ be a finite set such that $|Q|\geq 
|B_r(\Gamma,\Sigma,e_\Gamma)|$. Also, let $W\subset \cC$ be a clopen set and $\eta:W\to Q$
be a continuous function.
Suppose that for all pairs $x,y\in W$ such that $y=\alpha(\gamma)(x)$ for some $e_\Gamma\neq \gamma\in
B_r(\Gamma,\Sigma,e_\Gamma)$, we have that $\eta(x)\neq \eta(y)$.
Then there exists $\hat{\eta}:\cC\to Q$ such that $\hat{\eta}_{\mid W}=\eta$ and
for all pairs $x,y\in \cC$ such that $y=\alpha(\gamma)(x)$ for some $e_\Gamma\neq \gamma\in
B_r(\Gamma,\Sigma,e_\Gamma)$, we have that $\hat{\eta}(x)\neq \hat{\eta}(y)$.
\end{lemma}
\proof Let $\cup^s_{i=1} V_i=\cC$ be a clopen partition such that for any
$1\leq i \leq s$ and $x\neq y\in V_i$,
$$ \bigcup_{\gamma\in B_{2r}(\Gamma,\Sigma,e_\Gamma)} \alpha(\gamma)(x) \cap
\bigcup_{\gamma\in B_{2r}(\Gamma,\Sigma,e_\Gamma)} \alpha(\gamma)(y) =\emptyset\,.$$
\noindent
By induction, it is enough to show that there exists a continuous map $\eta': V_1\cup W\to Q$
such that $\eta'_{\mid W}=\eta$ and 
for all pairs $x,y\in V_1\cup W$ such that $y=\alpha(\gamma)(x)$ for some $e_\Gamma\neq \gamma\in
B_r(\Gamma,\Sigma,e_\Gamma)$, we have that $\eta'(x)\neq \eta'(y)$.
So, let us construct the map $\eta'$.
Let $z\in V_1\backslash W$ and $Q=\{1,2,\dots,k\}$. Let $\tau(z)$ be the smallest positive integer
such that there exists no  $\gamma\in
B_r(\Gamma,\Sigma,e_\Gamma)$ satisfying both of the two conditions below.
\begin{itemize}
\item $\alpha(\gamma)(z)\in W.$
\item $\eta(\alpha(\gamma)(z))=i.$
\end{itemize}
\noindent
Observe that $\tau^{-1}(i)$ is a clopen set for each $1\leq i \leq k$ and
$$V_1\cup W= W\cup \bigcup^k_{i=1} \tau^{-1}(i)\,.$$
\noindent
Define $\eta': V_1\cup W\to Q$ by setting
\begin{itemize}
\item $\eta'(z)=\eta (z)$ for $z\in W$,
\item  $\eta'(z)=i$, for $z\in  \tau^{-1}(i)$.
\end{itemize}
\noindent
It is easy to see that the map $\eta'$ satisfies the condition of our lemma. \qed
\vskip 0.1in
\noindent
Since the previous lemma implies Proposition \ref{pf12}, our theorem follows. \qed
\subsection{The space of the weak equivalence classes is compact}
We can introduce a natural metric on the set $\Free(\Gamma)$. Let $CF(n,k)$ be the set of
all subsets of $\{\{0,1\}^k\}^{B_n(\Gamma,\Sigma,e_\Gamma)}$.
Let $\alpha,\beta:\Gamma\actson\cC$ be free Cantor actions.
The weak distance of the classes $[\alpha]$ and $[\beta]$ is defined in the following way.
$$d_w([\alpha],[\beta])=2^{-n}$$
\noindent
if
\begin{itemize}
\item there exists $\cA\in CF(n,n)$ such that
either $\cA\in\conf_{n,\alpha}(\phi)$ for some continuous map $\phi:\cC\to \{0,1\}^n$
and  $\cA\notin\conf_{n,\beta}(\psi)$ for any continuous map $\psi:\cC\to \{0,1\}^n$, or
$\cA\in\conf_{n,\beta}(\psi)$ for some continuous map $\psi:\cC\to \{0,1\}^n$
and  $\cA\notin\conf_{n,\alpha}(\phi)$ for any continuous map $\phi:\cC\to \{0,1\}^n$,
\item for any $1\leq i \leq n-1$ and $\cB\in CF(i,i),$
$\cB\in\conf_{i,\alpha}(\phi)$ for some continuous map $\phi:\cC\to \{0,1\}^i$ if and only if
$\cB\in\conf_{i,\beta}(\psi)$ for some continuous map $\psi:\cC\to \{0,1\}^i$.
\end{itemize}
\noindent
The following theorem is the qualitative analogue of the main result of \cite{AE} 
(see also \cite{TDR}) on the compactness of the weak equivalent classes of p.m.p. actions.
\begin{theorem}\label{compact}
The space $\Free(\Gamma)$ is compact with respect to the metric $d_w$.
\end{theorem}
\proof
First, we introduce some operations on configuration sets.
For $ 1\leq n \leq m$ and $k\geq 1$, let
$$\rho^{m,n}: \{\{0,1\}^k\}^{B_m(\Gamma,\Sigma,e_\Gamma)}\to 
\{\{0,1\}^k\}^{B_n(\Gamma,\Sigma,e_\Gamma)}\,,$$
\noindent
be the restriction map. So, if $\cA\in CF(m,k)$, then $\rho^{m,n}(\cA)\in CF(n,k)$.
Also, if $1\leq a\leq b\leq k$ and $l=b-a+1$, let the map $\pi^{[a,b]}:\{0,1\}^k\to \{0,1\}^l$ be
defined by setting for $1\leq i \leq k$,
$$\pi^{[a,b]}(c)(i)=c(a-i+1)\,.$$
\noindent
That is, if $\cA\in CF(n,k)$, then $\pi^{[a,b]}\circ\cA\in CF(n,l).$
Now, let $\{\alpha_j:\Gamma\actson\cC\}^\infty_{j=1}$ be a sequence
of free actions such that the classes $\{[\alpha_j]\}^\infty_{j=1}$ form a Cauchy sequence
in the $d_w$-metric. We need to show that there exists $\alpha:\Gamma\actson\cC$ such
that $\lim_{j\to\infty} [\alpha_j]=[\alpha]$.
For $n\geq 1$, let us call $\cA\in CF(n,k)$ a {\bf surviving configuration}
if for all but finitely many $j$'s, there exists a continuous map
$\phi^j:\cC\to\{0,1\}^k$ such that
$\conf_{n,\alpha_j}(\phi^j)=\cA$, or in shorthand notation, $\cA\in\conf(\alpha_j).$
Our goal is to construct a free
action $\alpha:\Gamma\actson \cC$ such that $\cA\in \conf(\alpha)$ 
if and only if $\cA$ is a surviving configuration.
Let $\{\cA_i\}^\infty_{i=1}$ be an enumeration of the surviving configurations, where
$\cA_i\in\conf_{n_i,l_i}$.
Let $m_1\leq m_2 \leq \dots $ be a sequence of integers such that if
$1\leq j \leq i$, then $m_i\geq n_j$. For all pairs $ 1\leq j \leq i$, let $\hat{A}^i_j\in CF(m_i,l_j)$
be a surviving configuration such that
$\rho^{m_i,n_j} (\ha^i_j)=\cA_j$.
For $i\geq 1$, let $k_i=\sum^i_{j=1} l_j$.
Let $a_1\leq b_1< a_2\leq b_2<\dots$ be an infinite sequence of integers such
that for any $i\geq 1$, $a_{i+1}=b_i$ and $b_i-a_i+1=l_i$.
Let $CO(i)\subset CF(m_i,k_i)$ be the set of all surviving
configurations $\cB$ such that for any $1\leq j\leq i$, $\pi^{[a_j,b_j]}\circ \cB=\ha^i_j.$
Note that if $\cB\in CO(i)$, then
$$\pi^{[1, k_{i-1}]}\circ (\rho^{m_i,m_{i-1}}(\cB))\in CO(i-1)\,.$$
\noindent
Therefore by Konig's Lemma, we have a sequence of surviving configurations
$\{ \cB_i\in CO(i)\}^\infty_{i=1}$ such that for any $i\geq 1$,
$$\pi^{[1, k_{i-1}]}\circ (\rho^{m_i,m_{i-1}}(\cB_i))= \cB_{i-1}\,.$$
\noindent
For $n\geq 1$, let $1\leq i_n\leq n$ be the largest integer such that $\cB_{i_n}\in \conf(\alpha_n)$ (we can 
assume that $\cB_1 \in \conf(\alpha_n)$, for every $n\geq 1$).
Clearly $i_n\to \infty$ as $n\to\infty$. Let $\phi_n:\cC\to \{0,1\}^{k_{i_n}}$ be a continuous map
such that $\conf_{m_{i_n},\alpha_n}(\phi_n)=\cB_{i_n}$ and let $\kappa^{\phi_n}_{\alpha_n}(\cC)=Y_n\subset \cC^\Gamma$
be the corresponding Cantor subshift. Note 
that by Lemma \ref{lf16}, $Y_n$ is conjugate isomorphic to $\alpha_n$.
Let $Y\subset \cC^\Gamma$ be defined in the following
way. Let $y\in Y$, if there exists a sequence $\{y_{n_k}\in Y_{n_k}\}^\infty_{k=1}$
such that $\lim_{k\to\infty} y_{n_k}=y$.
Our theorem immediately follows from the following proposition.
\begin{proposition}
\label{pf31}
The space $Y$ is a Cantor subshift and $[Y_n]\to [Y]$ in \\ $\Free(\Gamma)$.
\end{proposition}
\proof
Clearly, $Y$ is closed and $\Gamma$-invariant.
\begin{lemma} \label{lf32}
The shift action of $\Gamma$ on $Y$ is free.
\end{lemma}
\proof
Let $\beta:\Gamma\actson\cC$ be the smallest
element in $\Free(\Gamma)$ as in Theorem \ref{AWe}.
So, any $\cA\in\conf(\beta)$ is a surviving configuration.
Let $m\geq 1$ and let $\cA\in \conf_{m,r_m}(\beta)$ such
a configuration that
$c(\gamma_1)\neq c(\gamma_2)$, provided that
$c: B_m(\Gamma,\Sigma,e_\Gamma)\to \{0,1\}^{r_m}\in \cA$ and $\gamma_1\neq\gamma_2\in B_m(\Gamma,\Sigma,e_\Gamma)$.
Therefore, there exists some constant $s_m>0$ such that
for large enough $i$,
$$(y_i(\gamma_1))_{[s_m]}\neq (y_i(\gamma_2))_{[s_m]}\,,$$
\noindent
provided that $d_{Cay(\Gamma,\Sigma)}(\gamma_1,\gamma_2)\leq m.$
Therefore, for such $\gamma_1$ and $\gamma_2$,
$$(y(\gamma_1))_{[s_m]}\neq (y(\gamma_2))_{[s_m]}\,.$$
Consequently, if $e_\Gamma\neq \delta\in B_m(\Gamma,\Sigma,e_\Gamma)$, then 
\begin{equation}\label{ef32}
R(\delta)(y)\neq y,\end{equation}\noindent
where $L$ is translation action on $\cC^\Gamma$.
 Since \eqref{ef32} holds for all $m\geq 1$, the action of $\Gamma$ on $Y$ is free. \qed
\vskip 0.1in
\noindent
We need to show that $\conf(Y)$ equals to the set of all surviving configurations.
Let $\pi_i:Y\to \{0,1\}^{k_i}$ be the projection onto the first $k_i$ coordinates and $n\geq 1$.
\begin{lemma}
\label{l2f1}
The configuration $\conf_{n,Y}(\pi_i)$ is a surviving configuration.
\end{lemma}
\proof By our construction, there exists $\cA\in CF(n,k_i)$ such that
for all large enough $k\geq 1$, $\cA=\conf_{n,Y_k}(\pi_i)$.
We need to show that $\cA=\conf_{n,Y}(\pi_i)$.
Clearly, $\conf_{n,Y}(\pi_i)\subseteq \cA$.
Now, let $c: B_n(\Gamma,\Sigma,e_\Gamma)\to \{0,1\}^{k_i}\in\cA$.
Then, we have a sequence $\{y_k\in Y_k\}^\infty_{k=1}$ so that for
any $\gamma\in B_n(\Gamma,\Sigma,e_\Gamma)$, $\pi_i(y_k(\gamma))=c(\gamma)$.
So, $\pi_i(y(\gamma))=c(\gamma)$, whenever $y$ is a limit point of the sequence $\{y_k\}^\infty_{k=1}$.
Therefore, $c\in \conf(Y).$ \qed
\vskip 0.1in
\noindent
Let $\cA\in CF(n,k)$ be a surviving configuration. Then, there exists $i\geq 1$ and some integers
$a,b$ such that
$$\pi^{[a,b]}\circ(\rho^{i,n} (\cB_i))=\cA\,.$$
\noindent
Since by the previous lemma, $\cB_j\in\conf(Y)$, we can immediately see that $\cA\in\conf(Y)$ as well.
Now, let $\cA\in\conf_{n,k}(Y)$ and let $\phi:Y\to \{0,1\}^k$ be a continuous map such 
that $\conf_{n,Y}(\phi)=\cA$.
By continuity, there exists $i\geq 1$, such that for any Cantor subshift $Z\subset \cC^\Gamma$ for
which $\conf_{i,Z}(\pi_i)=\conf_{i,Y}(\pi_i)$, there exists some continuous map $\phi_Z$ so that
$$\conf_{n,Z}(\phi_Z)=\conf_{n,Y}(\phi)=\cA\,.$$
\noindent
Since for large enough $k\geq 1$,
$$\conf_{i,Y_k}(\pi_i)=\conf_{i,Y}(\pi_i)\,,$$
\noindent
we can see that $\cA\in\conf(Y_k)$ for large enough $k$. Therefore, $\cA$ is a surviving configuration and
so, our proposition follows. \qed
\vskip 0.1in
\noindent
This concludes the proof of our theorem. \qed
\begin{remark} \label{r2f3}
In \cite{Elekmin}, for any free Cantor action 
$\alpha:\Gamma\actson\cC$ we constructed a specific minimal Cantor action $m_\alpha:\Gamma\actson\cC$.
By the construction it is clear that $\alpha\succeq m_\alpha$ holds for any free action $\alpha$.
Therefore, we can conclude that the smallest element of $\Free(\Gamma)$ can be represented by a minimal
action.
\end{remark}
\begin{remark}\label{2r2f3}
Let $\{[\beta_i]\}^\infty_{i=1}$ be a dense subset of $\Free(\Gamma)$.
Then, we can consider the Cantor action $\beta=\prod_{i=1} \beta_i\actson \prod_{i=1} \cC$.
Clearly, for any $i\geq 1$, $\beta\succeq \beta_i$. Hence, $\beta\succeq\alpha$ for any free action
$\alpha:\Gamma\actson\cC$. That is, $[\beta]$ is the largest element of $\Free(\Gamma)$.
\end{remark}
\section{Qualitative convergence and  limits}\label{qual}
The goal of this section is to introduce and study the main notion of our
paper: qualitative convergence of countable graphs.
\subsection{Convergence of generalized Schreier graphs}
First, we extend Definition \ref{defqual} for Schreier graphs.
Let $G\gsg$ be the set of all generalized Schreier graphs as in Section \ref{tgc}.
Let $Q$ be a finite set and $U^{r,Q}_{\Gamma,\Sigma}$ be the set of all
rooted, $\Sigma$-edge labeled, $Q$-vertex labeled balls of radius $r$ that occurs
in some $Q$-labeled Schreier graphs.
Now, let $S\in G\gsg$ and $\phi: V(S)\to Q$ be a continuous map.
Then, let $\conf_{r,S}(\phi)\subseteq U^{r,Q}_{\Gamma,\Sigma}$, the
configuration 
of $\phi$, be
the set of all rooted, $\Sigma$-edge labeled,  $Q$-vertex labeled balls 
that occur in in the labeled graph $(S,\phi)$. Also, we denote by $\conf(S)$ the set of all configurations of $S$.
\begin{definition}\label{d2f4}
A sequence of generalized Schreier graphs $\{S_n\}^\infty_{n=1}\subset G\gsg$ is 
qualitatively convergent, if for any finite set $Q$, integer $r\geq 1$ and
$\cA\subseteq  U^{r,Q}_{\Gamma,\Sigma}$, there exists
$N_\cA>0$ such that
\begin{itemize}
\item either for all $n\geq N_\cA,$ there exists $\phi_n:V(S_n)\to Q$ \\ such that
$\conf_{r,S_n}(\phi)=\cA$,
\item or for all $n\geq N_\cA,$  there exists no $\phi_n:V(S_n)\to Q$ \\ such that
$\conf_{r,S_n}(\phi)=\cA$.
\end{itemize}
\end{definition}
\noindent
Again, qualitative convergence does not depend on the choice of $\Sigma$.
\subsection{The basic example}
Let $\Gamma=\Z$, $\Sigma=\{1,-1\}$ and
 $C_n$ be the cyclic Schreier graph on the set $\Z / n\Z$. 
That is, $a\to b$ is a $1$-labeled directed edge if $b=(a+1)(\mbox{mod}\, n)$ and
$a\to b$ is a $-1$-labeled directed edge if $b=(a-1)(\mbox{mod}\, n)$. Our goal is
to explicitly describe the qualitatively convergent subsequences of $\{C_n\}^\infty_{n=1}$.
A subset $W\subseteq \N$ is called division-closed if $n\in W$ and $m\mid n$ always
 implies that $m\in W$ as well.
Division-closed subsets can be classified in the following way. Let $\{p_s\}_{s=1}^\infty$ be the
sequence of primes $2,3,\dots$. 
The family of division-closed subsets are in a one-to-one correspondence $\kappa\to W_\kappa$ with the
sequences $\kappa=\{\tau_s\}^\infty_{s=1}$, where $\tau_s\in \{\N\cup \{\infty\}\}$, and
$n\in W_\kappa$ if $n=p_{s_1}^{q_1} p_{s_2}^{q_2}\dots$ with $q_i<\tau_i$ for all $s\geq 1$.
\begin{example} If $\kappa=\{1,1,1,\dots\}$, then $W_\kappa=\{1\}$.
If $\kappa=\{\infty,\infty,\dots\}$, then $W_\kappa=\N$. If $\kappa=\{\infty,1,1,1,\dots\}$, then
$W_\kappa=\{1,2,4,8,\dots\}$.
\end{example}
\noindent
Let $\rho=\{n_k\}^\infty_{k=1}$ be a strictly increasing sequence of integers. We say that $\rho$ is arithmetically convergent,
if for any $m>1$
\begin{enumerate}
\item either $m\mid n_k$ for all but finitely many $k$'s,
\item or $m\nmid n_k$ for all but finitely many $k$'s.
\end{enumerate}
\noindent
The set of integers $m$ satisfying condition (1) form the arithmetic spectrum of $\rho$, $\mbox{Ar}(\rho)$.
Clearly, $\mbox{Ar}(\rho)$ is a division-closed set for any arithmetically convergent sequence $\rho$.
Also, for any division-closed set $W$, there exists an arithmetically convergent sequence $\rho$ such that
$\mbox{Ar}(\rho)=W$.
\begin{proposition} \label{kedd1} Let $\{n_k\}^\infty_{k=1}$ be a strictly increasing sequence of integers Then, $\{C_{n_k}\}$ is
qualitatively convergent, if and only if $\{n_k\}^\infty_{k=1}$ is arithmetically convergent.
\end{proposition}
\proof 
Let $I^Q_r$ be the set of all $Q$-labelings of the directed path
$I_r=[-r,-r+1,\dots,r]$. We say that a labeling
$b\in I^Q_r$ is a continuation
of $a\in I_r^Q$, if there exists a labeling $c$ of the
path $[-r,-r+1,\dots,r+1]$ such that the restriction of $c$ onto
the segment $[-r,-r+1,\dots,r]$ equals to $a$ and the restriction
of $b$ onto the segment $[-r+1,-r+2,\dots,r+1]$ equals to $b$.
 For an example, let $Q=\{1,2,3\}$. The element $(1,2,3)\in I^Q_1$ is the continuation
of $(1,1,2),(2,1,2)$ and $(3,1,2)$. Also, the element $(1,1,1)\in I^Q_1$ is
the continuation of itself. So, let $\cA\subset U^{r,Q}_{\Z,\Sigma}=I^Q_r$ be as above, a set
of $Q$-labeled rooted $r$-balls in cyclic Schreier graphs. 
We construct a directed graph $G_\cA$ using the structure of $\cA$ in the following
way.
\begin{itemize}
\item The vertex set of $G_\cA$ consists of the elements of $\cA$,
\item $a\to b$ is a directed edge of $G_\cA$ if $b$ is a continuation of $a$.
\end{itemize}
\noindent
Observe, that there exists a map $\phi_n:V(C_n)\to Q$ such that $\conf_{r,C_n}(\phi_n)=\cA$,
if there
exists a closed directed walk of length $n$ in the graph $G_\cA$ such that the walk passes through 
all the elements of $\cA$ at least once.  We call such a walk
an $\cA$-walk. Our proposition immediately follows from the lemma below.
\begin{lemma}\label{walk}
Let $\cA\subseteq I^Q_r$ be as above.  Suppose that $A$-walks exist in the graph
$G_\cA$.
Then, there exists a positive integer $n_\cA$ with the following property. For
any large enough integer $m$, there exists an $\cA$-walk of length $m$, if and only
if $m$ is divisible by $n_\cA$.
\end{lemma}
\proof
Observe that if $w_1$ is an $\cA$-walk of length $m_1$ and $w_2$ is an
$\cA$-walk of length $m_2$ then the concatenation of $w_1$ and $w_2$ is an
$\cA$-walk of length $m_1+m_2$. Let $n_\cA$ be the
greatest common divisor of the $\cA$-walk lengths. Then, $n_\cA=\sum^k_{i=1} p_im_i$,
where $\{p_i\}^k_{i=1}$ are integers and $\{m_i\}^k_{i=1}$ are lengths of $\cA$-walks.
So, if $n$ is a large enough integer divisible by $n_A$, then $n=\sum^k_{i=1} q_i m_i$, where
$\{q_i\}^k_{i=1}$ are non-negative integers. Consequently, there exists
an $\cA$-walk of length $n$. \qed
\subsection{Stable actions as qualitative limits}
Let $X$ be a totally disconnected compact metric space of continuum size. Then $X=A\cup \cC$,
where $A$ is a countable set.
Let $\alpha:\Gamma\actson X$ be a continuous action. Following
Glasner and Weiss \cite{GW}, we call $\alpha$ a {\bf stable action}
if for any $\gamma\in\Gamma$ and $x\in X$ such that $\alpha(\gamma)(x)=x$, there exists
a neighbourhood $x\in U$ such that $\alpha(\gamma)(y)=y$ for all $y\in U$. Clearly,
any free action is stable.
Let $\Sub(\Gamma)$ be the compact metric space of all the subgroups of $\Gamma$.
Recall that we regard $\Sub(\Gamma)$ as the closed subset of the space of subsets 
$\{0,1\}^\Gamma$, with the natural conjugation action of $\Gamma$.
If $\beta:\Gamma\actson X$ is a continuous action, we have a natural map
$\stab_\beta: X\to \Sub(\Gamma)$ so that for $x\in X$, $\stab_\beta(x)$ is the stabilizer
subgroup of $x$. The action $\beta$ is stable if and only if $\stab_\beta$ is continuous
(see \cite{GW}). In this case, $\stab_\beta(X)$ is a $\Gamma$-invariant closed subset of
$\Sub(\Gamma)$.
Now, let $Q$ be a finite set and $\phi: X\to Q$ be a continuous map. 
Let $\conf_{r,\beta}(\phi)\subseteq U^{r,Q}_{\Gamma,\Sigma}$ be the
set of all rooted $Q$-labeled balls $\eta$ of
radius $r$ such that there exists $x\in X$ so that the $Q$-labeled
ball around $x$ is isomorphic to $\eta$.
\begin{definition}
Let $\{S_n\}^\infty_{n=1}\subset G\gsg$ be 
a qualtitatively convergent sequence of generalized Schreier graphs.
We say that $\{S_n\}^\infty_{n=1}$ is  qualitatively converges to the stable
action $\beta:\Gamma\actson X$, $S_n\stackrel{q}{\to} \beta$, if for any finite set $Q$ and $\cA\subseteq U^{r,Q}_{\Gamma,\Sigma}$,
there exists a continuous map $\psi:X\to Q$ such that $\conf_{r,\beta}(\psi)=\cA$ if and only
if for all large enough $n\geq 1$, there exists $\phi_n:V(S_n)\to Q$ such that
$\conf_{r,S_n}(\phi_n)=\cA$.
\end{definition}
\begin{remark}\label{r2f7}
Let $\beta:\Gamma\actson \cC$ be a free Cantor action such that
for some sequence of generalized Schreier graphs
$S_n\stackrel{q}{\to} \beta$. Suppose that the free Cantor action
$\alpha:\Gamma\actson \cC$ is qualitatively weakly equivalent to $\beta$.
Then $S_n\stackrel{q}{\to} \alpha$, as well.
\end{remark}
\noindent
The main result of this section is the following qualitative analogue of
the theorem of Hatami, Lov\'asz and Szegedy \cite{HLSZ} on the existence of local-global limits.
\begin{theorem} \label{t2f8}
For any qualtitatively convergent sequence $\{S_n\}^\infty_{n=1}\subset G\gsg$, there exists
a totally disconnected continuum $X$ and a stable action $\beta:\Gamma\actson X$ such that
$S_n\stackrel{q}{\to} \beta$.
\end{theorem}
\proof
The proof of our theorem will be very similar to the one of Theorem \ref{compact}.
Again, it is enough to check convergence for finite sets $Q$, where
$Q=\{0,1\}^k$ for some $k\geq 1$.
As in the previous section, we call $\cA\subseteq U_{\Gamma,\Sigma}^{r,\{0,1\}^k}$ a surviving configuration,
if for large enough $n\geq 1$, $\cA\in \conf(S_n)$. Let $\{\cA_i\}^\infty_{i=1}$ be
an enumeration of the surviving configurations, where $\cA_i\subseteq U _{\Gamma,\Sigma}^{r_i,\{0,1\}^{k_i}}$.
Now, we construct a sequence of Cantor labelings 
$\{\psi^n:V(S_n)\to \cC\}^\infty_{n=1}$.
Let $\prod^\infty_{j=1} \psi^n_j: V(S_n)\to \prod^\infty_{j=1} \{0,1\}^{k_j}$ be defined in
the following way.
\begin{itemize}
\item If $\cA_j\in\conf(S_n)$, then let $\psi^n_j:V(S_n)\to\{0,1\}^{k_j}$ be a
function such that $\conf_{r_j,S_n} (\psi^n_j)=\cA_j$,
\item if $\cA_j\notin\conf(S_n)$, then let $\psi^n_j(v)=\{0,0,\dots,0\}$ for all $v\in V(S_n)$.
\end{itemize}
\noindent
Consider the totally disconnected compact space $\cC \gsg$
of rooted $\cC$-labeled Schreier graphs.
Let $(S,\phi)\in \cC G\gsg$. Then its orbit closure, $\overline{O}((S,\phi))$ is defined as in
Section \ref{tgc}.
Now, let us define $X\subset\cC \gsg$ as follows.
Let $y\in X$, if $y=\lim_{k\to\infty} y_{n_k}$ for some convergent sequence
$\{y_{n_k}\}^\infty_{k=1}\subset \cC \gsg$, where $y_{n_k}\in \overline{O}((S_{n_k},\phi_{n,k})).$
Clearly, $X$ is closed and invariant under $\Gamma$.
\begin{lemma} \label{l2f10}
The action of $\Gamma$ on the space $X$ is stable.
\end{lemma}
\proof
Before starting the proof let us make a remark on stability. 
Let $z=(T,H,\rho)\in \cC\gsg$, where $T$ is a Schreier graph, $H\in \Sub(\Gamma)$ and
$\rho:V(T)\to \cC$. Then, we have a natural rooted Schreier graph structure $z'=(T',H',\rho')$
on the orbit set of $z$ in $ \cC\gsg$.
One can observe that $z$ is not always equal to $z'$. For example, if
$T=\mbox{Cay}(\Gamma,\Sigma)$, $H=e_\Gamma$ and $\rho:\Gamma\to \cC$ is a constant-valued
function, then $V(T')$ consists of a singleton.
Now suppose that for any $r\geq 1$, there exists some $s_r\geq 1$
such that
$$(\rho(x))_{[s_r]} \neq (\rho(y))_{[s_r]}\,,$$
\noindent
provided that $0<d_T(x,y)\leq r$. Then, $(T,H,\rho)=(T',H',\rho')$ and
also, the $\Gamma$-action on the orbit closure $\overline{O}(z)$ is stable
(see e.g. the proof of Corollary 3.1. in \cite{Elekmin}).
If $Y\subset \cC \gsg$ is a closed $\Gamma$-invariant subset, $r\geq 1$, $Q$ is a finite set
and $\phi:Y\to Q$ is a continuous map, then one can consider the set
$$\conf^1_{r,Y}(\phi)=\cup_{y\in Y} \conf_{r,T_y}(\phi)\,,$$
\noindent
where $y=(T_y,H_y,\rho_y)$.
In general, it is possible that $\conf_{r,Y}(\phi)\neq \conf^1_{r,Y}(\phi)$.
However, it is clear from the discussion above that
\begin{equation} \label{e2f12}
\conf_{r,Y}(\phi)= \conf^1_{r,Y}(\phi)
\end{equation}
\noindent
holds, whenever the action of $\Gamma$ on $Y$ is stable.
Now, let us turn back to the proof of our lemma.
Observe that for any $r\geq 1$, there exists $i\geq 1$ and a surviving configuration
$\cA_i\subseteq U^{r_i,\{0,1\}^{k_i}}_{\Gamma,\Sigma}$ such that $r_i>r$ and
$\cA_i$ is $r$-separating. That is, for any $c\in \cA_i$, 
$c:B\to \{0,1\}^{k_i}$, where $B\in U^{r_i}_{\Gamma,\Sigma}$,
we have that
$c(u)\neq c(v)$ provided that $0<d_B(u,v) \leq r$.
Hence, for large enough $n\geq 1$, 
\begin{equation}
\label{e2f13}
(\hat{y}(a))_{[s_r]}\neq (\hat{y}(b))_{[s_r]},
\end{equation}
\noindent
for any $\hat{y}=(T_{\hat{y}}, H_{\hat{y}}, \rho_{\hat{y}})\in\overline{O}(S_n,\phi_n)$, where $0< d_{T_{\hat{y}}}(a,b)\leq r$ and $s_r=\sum_{j=1}^ik_j$.
Therefore by definition, $\eqref{e2f13}$ also holds if $\hat{y}\in Y$.
Hence, the action of $\Gamma$ on $Y$ is stable.

\noindent
Again, we need to prove that $\conf(Y)$ equals to the set of surviving
configurations. We can repeat the proof of Theorem \ref{compact}, to show
that if $\cA$ is a surviving configuration, then $\cA\in\conf(Y).$
It is also clear from our construction, that $\conf_{r,Y}(\pi_i)$ is
a surviving configuration for any $r\geq 1$ and $i\geq 1$.
So, the continuity argument after Lemma \ref{l2f1} can be applied to
immediately show that all elements of $\conf(Y)$ are, in fact, surviving
configurations. \qed
\begin{remark} \label{r27}
Let us suppose that for a qualitatively convergent sequence of generalized Schreier graphs
$\{S_n\}^\infty_{n=1}$, there exists an integer $q\geq 1$ such that each graph
$S_n$ contains exactly one vertex of degree $q$. Then if
$S_n\stackrel{n}{\to}\alpha $, where
$\alpha:\Gamma\actson X$ is a stable action, $X$ always contains an isolated 
point. 
\end{remark}
\begin{remark} \label{r272}
Let $S\in G\gsg$ be a generalized Schreier graph. Then, the vertices of $S$ form a
conjugacy invariant subset $A(S)$ in $\Sub(\Gamma)$. Let $\overline{A(S)}$ be
the closure of $A(G)$. If $S_n\stackrel{n}{\to} S$, then it is not hard to to see
that $\overline{A(S)}$ is the Hausdorff limit of the sequence of compact
subsets $\{\overline{A(S_n)}\}^\infty_{n=1}$. Moreover, if
$S_n\stackrel{q}{\to}\beta$, where $\beta:\Gamma\actson X$ is a stable
action, then the closed subset $\Stab_\beta(X)$ coincides with $\overline{A(S)}$.
\end{remark}

\subsection{The limits 
 of cyclic Schreier graphs}
 Let $\alpha:\Z\actson \cC$ be a free minimal Cantor action of the
integers. Following Lin and Matui \cite{LM}, we say that the integer
$n\in\N$ is an element of the {\bf periodic spectrum} of $\alpha$, 
$PS(\alpha)$, if there exists a clopen set $U\subset \cC$ such
that $\cup^{n-1}_{i=0} \alpha(i)(U)$ is a clopen partition of $\cC$.
Clearly, $PS(\alpha)$ is a division-closed set.
In \cite{Shi} (Lemma 3.14), Shimomura constructed
a free minimal Cantor action $\alpha_W:\Gamma\actson \cC$ for every
 division-closed set $W$ such that $PS(\alpha_W)=W$.
Now, we characterize the qualitative limits of cyclic Schreier
graphs using the periodic spectrum.
\begin{theorem} \label{t2f15}
Let $W$ be a  division-closed set and
$\{C_{n_k}\}^\infty_{k=1}$ be a qualitatively convergent
sequence of cyclic Schreier graphs such that
$Ar(\{n_k\}^\infty_{k=1})=W$.
Then, $\alpha_W$ is a qualitative limit for $\{C_{n_k}\}^\infty_{k=1}$.
\end{theorem}
\proof
Let $\cA\subset U^Q_r$ be a set as in Lemma \ref{walk}, for which
$\cA$-walks exist and let $n_\cA$ be the
smallest positive integer for which
 $\cA$-walks of length $m$ exist for large enough integers $m$.
First, suppose that $n_\cA\in W$. That is, $\cA$ is a surviving
configuration for $\{C_{n_k}\}^\infty_{k=1}$.
\begin{lemma} \label{l2f15}
$\cA\in \conf(\alpha_W)$.
\end{lemma}
\proof
Since $n_\cA\in PS(\alpha_W)$, we have 
a clopen set $U\subset \cC$ such that \\
$\cup^{n_\cA-1}_{i=0}\alpha_W(i)(U)$ is a clopen partition
of $C$. Let $l\geq 1$ be an integer such that if $l\leq m$ and
$n_\cA\mid m$, then there exists a $\cA$-walk of length $m$.
As in Lemma \ref{lf23}, we can easily prove that there exists a clopen set
$V\subset U$ such that if $x\in V$, then $\alpha_W(i)(x)\notin V$ for any
$1\leq i \leq l$.
By minimality, there exists some $l<k$ such
that for all $x\in V$, $\alpha_W(i)(x)\in V$ for some
$l<i\leq k$. Pick
$$c:=(\hat{q}_{-r},\hat{q}_{-r+1},\dots,\hat{q}_0,\hat{q}_1,\dots,\hat{q}_r)\in
\cA\,.$$
\noindent
For any integer $l<m\leq k$ such that $n_\cA\mid m$, pick an $\cA$-walk $wk(m)$ of length $m$, that
starts and ends with $c$. That is, for each such $m$ we have a sequence
of elements from $\cA$,
$$(wk^m_0, wk^m_1,\dots, wk^m_{m-1})\,,$$
\noindent
such that $wk^m_0=wk^m_{m-1}=c$,
and a corresponding sequence of elements of $Q$,
$$(q^m_0, q^m_2,\dots, q^m_{m-1})\,,$$
\noindent
where $q^m_i=wk^m_i(0)$.
\noindent
Now we define the continuous maps 
$$\phi_1:\cC\backslash V\to \{1,2,\dots, k\}\,\,\mbox{and} \,\,
\phi_2:\cC\backslash V\to \{1,2,\dots, k\}$$
\noindent
in the following way.
For $y\in \cC\backslash V$, there exist  unique positive integers $i,m$, $i<m\leq k$, $n_{\cA}\mid m$ and
$x_1\in V$, $x_2\in V$
such that
\begin{itemize}
\item $\alpha_W(i)(x_1)=y, \alpha_W(m)(x_1)=x_2,$
\item $\alpha_W(j)(x_1)\notin V$ for $1\leq j \leq m-1$.
\end{itemize}
\noindent
Then, let $\phi_1(y)=m, \phi_2(y)=i\,.$
Now, we can define the continuous map $\psi:\cC\to Q$ by
\begin{itemize}
\item $\psi(x)=\hat{q}_0$ for $x\in V$,
\item $\psi(y)=q_{\phi_2(y)}^{\phi_1(y)}$ for $y\in \cC\backslash V$.
\end{itemize}
\noindent
It is not hard to see that $\conf_{r,\psi}(\alpha_W)=\cA$.
Hence, the lemma follows. \qed
\begin{lemma}\label{l2f17}
Let $\beta$ be a minimal action such that $\cA\in\conf(\beta)$.
Then, $n_{\cA}\in PS(\beta)$.
\end{lemma}
\proof
Let $\phi:\cC\to Q$ be a continuous map such that
$\conf_{r,\phi}(\beta)=\cA$.
Then, we have a corresponding continuous map
$\Phi:\cC\to \cA$ such that for each $x\in \cC$
$$\Phi(x)=\left(\phi(\beta(-r)(x)),\phi(\beta(-r+1)(x)),\dots, \phi(\beta(r)(x))\right)\,.$$
\noindent
Let $c\in \cA$. By minimality, there exists some $t\geq 1$ such that
if $\Phi(x)=c$ $\Phi(\beta(i)(x))=c$ and $i\geq t$, then $n_\cA\mid i$ and the
 sequence
$$\left(\Phi(x),\Phi(\beta(1)(x)),\dots, \Phi(\beta(i)(x))\right)$$
defines an $\cA$-walk starting and ending at $c$.
Again, let $V\subset \Phi^{-1}(c)$ be a clopen set such that if $x\in V$ and
$j\leq t$, then $\beta(i)(x)\notin V$.
So, we can define the continuous map $\lambda:\cC\to\{0,1,\dots,n_\cA-1\}$
by setting
\begin{itemize}
\item $\lambda(x)=0$ if $x\in V$,
\item $\lambda(y)=i$, if $j$ is the smallest integer
such that $\beta(j)(x)=y$, for some $x\in V$ and $j\equiv i\, \mbox{mod} (n_\cA)$.
\end{itemize}
\noindent
Then, $\cup^{m-1}_{i=0} \beta(i)(U)$ is a clopen partition of $\cC$, where
$U=\lambda^{-1}(0).$ Hence, the lemma follows. \qed
\vskip 0.1in
\noindent
By the previous lemmas, $\cA\in \conf(\alpha_W)$ if and only if
$n_\cA\in W$, so our theorem follows. \qed
\subsection{The limits of countable graphs}
Let $\alpha:\Gamma\actson \cC$ be a stable action
of a finitely generated group $\Gamma$, with symmetric generating
set $S$. Then, for any $x\in\cC$, we can consider the simple, connected
orbit graph $G_x$, where
\begin{itemize}
\item $V(G_x)=\cup_{\gamma\in\Gamma} \alpha(\gamma)(x),$
\item the vertices $y\neq z\in V(G_x)$ are adjacent if
there exists $\sigma\in\Sigma$ such that $\alpha(\sigma)(y)=z$.
\end{itemize}
Let $\phi:\cC\to Q$ be a continuous map and let $d=\max_{x\in\cC} \deg_{G_x}(x)\,.$
Then $\conf_{r,\alpha}(\phi)\subset U^{r,Q}_d$ is the set of all rooted
$Q$-labeled balls of radius $r$, that occur in the
$Q$-labeled graph $(G_x,\phi)$ for some $x\in \cC$. Again, let
$\conf(\alpha)=\cup_{r,\phi} \conf_{r,\alpha}(\phi)\,.$
\begin{definition} \label{d2f20}
The action $\alpha:\Gamma\actson \cC$ is the limit of the
qualitatively convergent graph sequence $\{G_n\}^\infty_{n=1}\subset \grd$
if $\cA\in\conf(\alpha)$ if and only if $\cA\in\conf(G_n)$ for
large enough $n\geq 1$.
\end{definition}
\noindent
We have the following ``simple graph'' version of Theorem \ref{t2f8}.
\begin{proposition}
\label{p2f21}
For any qualitatively convergent graph sequence $\{G_n\}^\infty_{n=1}\subset \grd$,
there exists some finitely generated group $\Gamma$ and a stable
action $\alpha:\Gamma\actson\cC$ such that $\alpha$ is the limit
of $\{G_n\}^\infty_{n=1}.$
\end{proposition}
\proof
Before proving the proposition, let us make a short comment.
Let $Q$ be a finite set and $G\in\grd$ be a simple graph. Suppose that there
exists a labeling $\tau:E(G)\to Q$ such that adjacent edges have different
labels. 
Then, $\tau$ defines an action of the $Q$-fold free power $F_Q$ of the cyclic
group of two elements with generating system $Q$, where $q^2=e_{F_Q}$ for each
$q\in Q$. The underlying graph of the Schreier graph of the action $\tau$ is
just $G$. 
Suppose that  $\{G_n\}^\infty_{n=1}\subset \grd$ is a convergent graph
sequence, $\tau_n:E(G_n)\to Q$ is a sequence of maps as above and
$\{S_n\}^\infty_{n=1}$ is 
the associated sequence of Schreier graphs. If $S_n\stackrel{q}{\to} \beta$, where
$\beta$ is a stable action of $F_Q$, then $\beta$ is the qualitative limit of
our sequence $\{G_n\}^\infty_{n=1}$. It is not hard to show that one has maps
$\tau_n:E(G_n)\to Q$ such
that the associated Schreier graphs converge naively, but we do not know
whether one can achieve qualitative convergence.
Therefore, we need to pursue a somewhat different path towards the
proof. 
Let $\Gamma$ be a finitely generated group and $\Sigma$ be a symmetric
generating set, $|\Sigma|=d$. Then, we have a natural ``forgetting'' map
$\cal{F}:\cC\gsg\to \crg$, mapping each rooted Schreier graph into the underlying
rooted simple graph. Clearly, $\cal{F}$ is continuous and maps invariant sets
into invariant sets.
\vskip 0.1in
\noindent
Let $Y\subset \crg$ be a closed invariant subset. We say that $Y$ is
{\bf proper} if for any $r>0$, there exists $s_r>0$ such that
if $(G,\psi,x)\in Y$, $x,y\in V(G)$ and $0<d_G(x,y)\leq r$, then
$(\psi(x))_{[s_r]}\neq (\psi(y))_{[s_r]}$. Now, let $\cup_{i=1}^q W_i$
be a clopen partition of $Y$ such that
if $(G,\psi,x)\in W_i$, $0<d_G(x,y)\leq 2$, then
$(G,\psi,y)\neq W_i$.
Clearly, such partition exists by properness.
Let $\psi: Y\to \{1,2,\dots, q\}$ be the continuous function such that
$\psi((G,\psi,x))=i$, if $(G,\psi,x)\in W_i$. Let $Q=\{1,2,\dots,q\}$
and $\hat{Q}$ be the set of $2$-element subsets of $Q$. Let
$F_{\hat{Q}}$ be the $|\hat{Q}|$-fold free power of the cyclic group of two
elements
with generating system $\hat{Q}$, where $a^2=e_{F_{\hat{Q}}}$, if $a\in \hat{Q}$.
For $(G,\psi,x)\in Y$
and
$y,z\in V(G)$, $d_G(y,z)=1$, label the edge $(y,z)\in E(G)$ by $\{\psi((G,\psi,y)),
\psi((G,\psi,z))\}\in \hat{Q}$.
Observe that for each $(G,\psi,x)\in Y$ we obtain a
rooted $F_{\hat{Q}}$-Schreier graph $\overline{(G,\psi,x)}$. The following
lemma is easy to prove.
\begin{lemma} \label{szerda1}
The set $\overline{Y}=\bigcup_{(G,\psi,x)\in Y} \overline{(G,\psi,x)}$
is a closed invariant subset of $\cC\gsg$, where $\Gamma=F_{\hat{Q}}$ and
$\Sigma=\hat{Q}$. The action of $F_{\hat{Q}}$ on $\overline{Y}$ is stable
and the forgetting map $\cal{F}:\overline{Y}\to Y$ is a homeomorphism.
\end{lemma}
\noindent
Let $Y\subset \crg$ be a proper subset. We can define the configuration
spaces as in the proof of Theorem \ref{t2f8}. Let $\phi:Y\to Q$ be
a continuous map, where $Q$ is a finite set. Then, $\conf_{r,Y}(\phi)\subset U_d^{r,Q}$ is the
set of all rooted $Q$-labeled balls that occur in some labeled graph
$(G,\psi)$,
where $(G,\psi,x)\in Y$. Note that we do not use the function $\psi$ only the
graph structure. Again, $\conf(Y)$ is the set of all configurations.
We say that $Y$ is a limit of a qualitatively convergent
sequence $\{G_n\}^\infty_{n=1}\subset \grd$ if $\conf(Y)$ equals to
the set of all surviving configurations. By Lemma \ref{szerda1}, our
proposition
follows from the following lemma.
\begin{lemma}\label{szerda2}
For each qualitatively convergent sequence $\{G_n\}^\infty_{n=1}\subset \grd$
there exists a proper subset  $Y\subset \crg$ such that 
$G_n\stackrel{q}{\to}Y$.
\end{lemma}
\proof The proof of the lemma is almost identical to the one of Theorem
\ref{t2f8}. We glance through the proof for completeness.
It is enough to check convergence for finite sets $Q$, where $Q=\{0,1\}^k$.
Let $\{\cA_i\}^\infty_{i=1}$ be an enumeration
of the surviving configurations, where
$\cA_i\subseteq U_d^{r_i,\{0,1\}^{k_i}}$.
Again, we construct a sequence of Cantor labelings 
$\{\psi^n:V(G_n)\to \cC\}^\infty_{n=1}$.
Let $\prod^\infty_{j=1} \psi^n_j: V(G_n)\to \prod^\infty_{j=1} \{0,1\}^{k_j}$ be defined in
the following way.
\begin{itemize}
\item If $\cA_j\in\conf(G_n)$, then let $\psi^n_j:V(G_n)\to\{0,1\}^{k_j}$ be a
function such that $\conf_{r_j,G_n} (\psi^n_j)=\cA_j$,
\item if $\cA_j\notin\conf(G_n)$, then let $\psi^n_j(v)=\{0,0,\dots,0\}$ for all $v\in V(G_n)$.
\end{itemize}
\noindent
Now, let us define $Y\subset\cC \gsg$ as follows.
Let $y\in Y$, if $y=\lim_{k\to\infty} y_{n_k}$ for some convergent sequence
$\{y_{n_k}\}^\infty_{k=1}\subset \crg$, where $y_{n_k}\in \overline{O}((G_{n_k},\phi_{n,k})).$
The following lemma can be proven in the same way as Lemma \ref{l2f10}.
\begin{lemma} \label{szerda3}
The space $Y\subset\cC \gsg$ is proper.
\end{lemma}
\noindent
Now, we need to show that $\conf(Y)$ is equal to the set of surviving
configurations. This can be done exactly the same way as in the end of the
proof
of Theorem \ref{t2f8}. \qed

\section{Almost finiteness of graphs and stable actions}
The goal of this section is to introduce and study various notions of almost
finiteness
for classes of bounded degree graphs. Our definition of almost finiteness is
based on the notion of almost finiteness of \'etale groupoids introduced by
Matui \cite{Matui}.
\subsection{The geometric groupoid of a stable action}
Let $X$ be a compact, metrizable, totally disconnected continuum. Let 
$\Gamma$ be a countable group and let $\alpha:\Gamma\actson X$ be
a stable action. Now using the action $\alpha$, we define a locally compact, totally disconnected, Hausdorff, principal,
\'etale groupoid $\cG_\alpha$, the {\bf geometric groupoid} of $\alpha$.
One should note that $\cG_\alpha$ is isomorphic to the transformation
groupoid of the action $\alpha$ only if the action is free.
\vskip 0.1\in
\noindent
$(1)$ The elements of the groupoid $\cG_\alpha$
are the pairs $(x,y)$, where $y=\alpha(\gamma)(x)$
for some $\gamma\in\Gamma$. As usual, we have the range and source
maps $s(x,y)=x$, $r(x,y)=y$, and a groupoid multiplication 
$(x,y)(a,b)=(x,b)$, where
$a=z$.
\vskip 0.1in
\noindent
$(2)$ The basis for the topology on $\cG_\alpha$
is given in the following way. Let $\alpha(\gamma)(x)=y$. Then
by stability,
\begin{itemize}
\item there exist clopen sets $x\in U$ and $y\in V$ such that
$\alpha(\gamma):U\to V$ is a homeomorphism,
\item and if also $\alpha(\delta)(x)=y$, then there exists
a clopen set $x\in W$ such that
$$\alpha(\gamma)_{\mid W}=\alpha(\delta)_{\mid W}\,.$$
\end{itemize}
\noindent
The base neighbourhoods of the element $(x,y)$ are in
the form of $(U,x,y)$, where $\alpha(\gamma):U\to V$ as above,
and $(a,b)\in (U,x,y)$ if $a\in U$ and $\alpha(\gamma)(a)=b$.
\vskip 0.1in
\noindent
One can immediately see that
\begin{enumerate}
\item The base set $\cG^0_\alpha$ is homeomorphic to $X$. The multiplication,
source and range maps are continuous.
\item The range map is a local homeomorphism, so our groupoid is \'etale.
\item The isotropy groups of all $x\in X$ are trivial.
\end{enumerate}
\noindent
Hence, $\cG_\alpha$ is a locally compact, totally disconnected,
Hausdorff, principal, \'etale groupoid. The groupoid $\cG_\alpha$
is minimal if and only if the action $\alpha$ is minimal.
Now, let us suppose that $\Gamma$ is finitely generated and
$\Sigma$ is a symmetric generating system of $\Gamma$.
Then, for any $x\in X$, the orbit set of $x$
$$\{(x,y)\in \cG_\alpha\}$$
\noindent
is equipped with a bounded degree graph structure $\cG_\alpha(\Sigma)$. 
The element
$(x,a)$ is adjacent to $(x,b)$ if $\alpha(\sigma)(a)=b$ for some
generator $\sigma\in \Sigma$.
Let $\cH_{\alpha,\Sigma}$ be the components of $\cG_\alpha(\Sigma)$.
Now, we define the almost finiteness  of $\alpha$.
\begin{definition}\label{d2f26}
The action $\alpha:\Gamma\actson X$ is almost finite if for any $\eps>0$
there exist $K_\eps>0$, a finite set $Q$ and a continuous
map $\phi_\eps:X\to Q$ satisfying the following conditions.
\begin{enumerate}
\item If $x,y\in X$, $\alpha(\gamma)(x)=y$ for some $\gamma\in\Gamma$,
and $\phi_\eps(x)=\phi_\eps(y)$, then either  
$d_G(x,y)\leq K_\eps$ or $d_G(x,y)\geq 3K_\eps$, 
where $G\in \cH_{\alpha,\Sigma}$ and $x\in V(G)$.
\item If $x$ and $G$ are as above, $i_G(H_x)\leq \eps$, where
$$H_x=\{z\in V(G)\,\mid\, d_G(x,z)\leq K_\eps\,\mbox{and}\,\phi_\eps(x)= \phi_\eps(z)\}\,.$$
\end{enumerate}
Clearly, almost finiteness of $\alpha$ does not depend
on the choice of the generating system.
\end{definition}
\noindent
Note that $\phi_\eps$ defines a continuous field of almost finite partitions
on the orbit graphs in $\cH_{\alpha,\Sigma}$. 
One can easily check that the partitions above satisfies
the conditions given in \cite{Suzuki} Definition 3.6.
Hence, we have the following definition-proposition.
\begin{proposition}\label{suzuki}
The action $\alpha$ is almost finite in the sense of Definition \ref{d2f26}
if and only if the principal \'etale groupoid $\cG_\alpha$ is almost finite
in the sense of \cite{Suzuki} Definition 3.6.
\end{proposition}
\noindent
Using Definition \ref{partoracle} in Section \ref{CTDA} and Theorem \ref{fotetel},  we  have
the following proposition.
\begin{proposition}\label{csut2}
The action $\alpha$ is almost finite if the set of graphs
$\cH_{\alpha,\Sigma}$ is a distributed almost finite class. In particular, the
action $\alpha$ is almost finite if $\cH_{\alpha,\Sigma}$ is a $D$-doubling
family for some $D>0$.
\end{proposition}
In \cite{Elekurs}, we proved (see the remark after Proposition 5.2.) that
for any real number $t\geq 1$, there exists a minimal, stable action $\alpha_t$
of some finitely generated group and some constant $C_t$ such that
for any orbit graph $H$, $r\geq 1$ and $x\in V(H)$,
$$\frac{1}{C_t} r^t\leq |B_r(H,x)|\leq C_t r^t\,.$$
\noindent
So, there exists some $D>0$ such that
all the orbit graphs $H$ are $D$-doubling (see Section \ref{doublal}).
Hence by Proposition \ref{csut2},
the action $\alpha_t$ is a minimal, almost finite Cantor action. Thus, by
the Main Theorem of \cite{Suzuki}, we have the following corollary.
\begin{corollary}
For all $\alpha\geq 1$, the groupoid $C^*$-algebra of $\alpha_t$ is a 
simple $C^*$-algebra of stable rank one.
\end{corollary}
\subsection{Almost finiteness and convergence}
\begin{proposition} \label{limitalmost}
Let $\{G_n\}^\infty_{n=1}\subset \grd$ be an almost finite resp. a
strongly almost finite family and suppose that
$G_n\tart G$. Then, $G$ is an almost finite graph resp. a strongly almost
finite graph. 
\end{proposition}
\proof
First, we need a lemma.
\begin{lemma}\label{equialmost}
Let $G\in\grd$ be a countable infinite graph such that
for some $\eps$ and $K$, there exists an $(\eps,K)$-partition of $V(G)$. Let $H\in\grd$ be
equivalent to $G$. Then, there exists an $(\eps,K)$-partition of $V(H)$.
\end{lemma}
\proof
Let $\cal{S}=\{S_i\}^\infty_{i=1}$ be an $(\eps,K)$-partition of $V(G)$. We can define a graph structure
$G_{\cal{S}}$ on $\cal{S}$. Let $S_i\neq S_j$ be adjacent, if there exist
vertices $x\in S_i$ and $y\in S_j$ so that $d_G(x,y)\leq 3K$. Observe that for any $i\geq 1$, the degree
of $S_i$ in the graph $G_{\cal{S}}$ is less than $d^{5K+1}$. Indeed, let $x\in S_i$. Then any tile $S_j$ adjacent
to $S_i$ is contained in the ball $B_{5K}(G,x)$. Also,
$|B_{5K}(G,x)|<d^{5K+1}$. So if $|Q|=d^{5K+1}$,  we can  label 
$\cal{S}$ by elements of the set $Q$ in such a way that adjacent elements have different labels. Let us lift the labeling above to a labeling $\phi:V(G)\to Q$
in such a way that $\phi(x)$ is the label of the tile containing $x$. Hence, we obtain an element $(G,\phi)\in\grd^Q$.
\begin{lemma}\label{tech1}
There exists a labeling $\psi:V(H)\to Q$ such that for any labeled ball $B_R(H,z,\psi)$ there exists
a labeled ball $B_R(G,y,\phi)$ that is rooted-labeled isomorphic to $B_R(H,z,\psi)$.
\end{lemma}
\proof
Fix a vertex $w\in V(H)$. For any $n\geq 1$, we define $B_n(w,\psi_n)$ to be rooted-isomorphic to some ball
$B_n(y_n,\phi)$. Then we consider a convergent subsequence,
$$\{B_{n_k}(w,\psi_{n_k})\}^\infty_{k=1}\stackrel{\rg^Q}{\rightarrow} (H',\psi')\,,$$
\noindent
where $H'$ is the component of $H$ containing $w$.
Clearly, for any $z\in V(H')$ and $R>0$, there exists a labeled ball $B_R(G,y,\phi)$ that is
rooted-labeled isomorphic to $B_R(H',z,\psi')$. We can finish the proof of the lemma, by defining $\phi$ for
all components of $H$. \qed
\vskip 0.1in
\noindent
Now we can finish the proof of Lemma \ref{equialmost}.
We can define a partition of $V(H)$ using the labeling $\psi:V(H)\to Q$. 
Let $x\equiv_\psi y$ if $\psi(x)=\psi(y)$ and $d_H(x,y)\leq K$. It is easy to check that
$\equiv_\psi$ is in fact an equivalence relation. Also, by Lemma \ref{tech1},
the partition defined by $\equiv_H$ is an $(\eps,K)$-partition of $V(H)$. \qed
\vskip 0.1in
\noindent
Now we can conclude the proof of our proposition.
Since $\{G_n\}^\infty_{n=1}$ form an almost finite class, the
graph $\overline{G}$ that consists
of disjoint copies of the graphs $\{G_n\}^\infty_{n=1}$
is itself almost finite.
Let $\phi:V(\overline{G})\to Q$ encode an $(\eps,K)$-partition of 
$V(\overline{G})$ as above. Thus we have a sequence
$\{(G_n,\phi_n)\}^\infty_{n=1}\subset \grd^Q$ and we can consider
a convergent subsequence $(G_{n_k},\phi_{n_k})\to (H,\phi)$ by Proposition
\ref{compactpro}, where $H\in \grd$. By the argument at the end of Lemma \ref{equialmost}, $V(H)$ has
an $(\eps,K)$-partition. Since 
$\{G_{n_k}\}^\infty_{k=1}\stackrel{n}{\rightarrow} G$
and $(G_{n_k},\phi_{n_k})\to (H,\phi)$, the graphs $G$ and $H$ are
equivalent. Thus, by Lemma \ref{equialmost}, $V(G)$ has an $(\eps,K)$-partition
as well. Therefore, $G$ is almost finite. The strong almost finite version can
be proven in exactly the same way. \qed
\vskip 0.1in
\noindent
For qualitative convergence, we have the following proposition.
\begin{proposition}\label{csututolso}
Let $\{G_n\}^\infty_{n=1}\subset\grd$ be a set of
finite graphs.
Suppose that $G_n\stackrel{q}{\to} \alpha$, where $\alpha:\Gamma\actson X$ is
some action of a finitely generated group $\Gamma$.
finite. Then, the family $\{G_n\}^\infty_{n=1}$ is almost finite if and only if $\alpha$ is almost finite.
\end{proposition}
\proof First, suppose that $\alpha$ is almost finite and fix $\eps>0$.
Let $\psi_\eps:X\to Q$ be the mapping of Definition \ref{d2f26} and consider
the configuration $\cA=\conf_{5K_\eps,\alpha}(\psi_\eps)$.
Since $\cA$ is a surviving configuration, we
have $N>0$ such that for any $n\geq N$ there exists
$\phi_n:V(G_n)\to Q$ such that
$$\conf_{5K_\eps,G_n}(\phi_n)=\cA\,.$$
\noindent
Therefore, the family $\{G_n\}_{n\geq N}$ is $(\eps,K_\eps)$-almost finite.
Thus, $\{G_n\}^\infty_{n=1}$ is an $(\eps, L_\eps)$-almost finite family,
where 
$$L_\eps=\max (K_\eps, \max_{1\leq i\leq N} |V(G_i)|)\,.$$
\noindent
Now, suppose that $\{G_n\}^\infty_{n=1}$ is an almost finite family. Let $\eps>0$.
Let $\{(G_n,\phi_n)\}^\infty_{n=1}\subset \grd^Q$ be as in the proof of Proposition \ref{limitalmost}.
 For $n\geq 1$, let
$\cA_n=\conf_{5K_\eps,G_n}(\phi_n)$.
Then, there exists $\cA$ such that $\cA=\cA_n$ for infinitely many $n's$.
That is, $\cA$ is a surviving configuration. Hence, there exists $\psi:X\to Q$ such that
$\cA=\conf_{5K_\eps,\alpha}(\psi)$. So, by Definition \ref{d2f26}, the action $\alpha$ is almost finite. \qed

\subsection{The fractal construction}
Now, as a preparation for Theorem \ref{nonamena}, we construct an almost finite,
connected infinite graph $G$ of bounded vertex degrees, which
is minimal in the suitable sense and contains an expander sequence of finite induced subgraphs. A somewhat similar
construction has been used in \cite{Elekurs}.
Let $K$ be a finite graph.
We call a subset $A_i\subset V(K)$ an {\bf $i$-subset} if
\begin{itemize}
\item for any $x\neq y\in A_i$, $d(x,y)>3^i$,
\item $\cup_{p\in A_i} B_{10^i}(K,p)=V(K)$.
\end{itemize}
The following lemma is easy to prove.
\begin{lemma}\label{triv1}
For any finite $3$-regular graph $L$ such that $10^k\leq \diam(L) $,
we have disjoint subsets $A_1,A_2,\dots,A_k\subset V(L)$ such that $A_i$
is an $i$-subset and
$\cup_{i=1}^k A_i \neq V(L).$
\end{lemma}
\noindent
Now let  $\{K_j\}^\infty_{j=1}$ be
an expander sequence of finite $3$-regular graphs such that
\begin{itemize}
\item for any $j\geq 1$, $\diam(K_j)\geq 10^j$ and
\item $|V(K_j)| \geq j 3^{(10^j+1)}\,.$
\end{itemize}
\noindent
Using Lemma \ref{triv1}, we pick disjoint subsets $A_1^j, A_2^j,\dots A^j_j\subset V(K_j)$ and a vertex $p_j$
such that $A^j_i$ is an $i$-subset
and $p_j\notin \cup_{i=1}^j A^j_i$. We call the vertex $p_j$ the connecting vertex of $K_j$.
\vskip 0.05in
\noindent
{\bf Step 1.} Let $H_1=K_1$. The graph $H_2$ is constructed in the following
way. First, we consider our graph $K_2$. For each
$q\in A^2_1\subset V(K_2)$, we choose
a copy $H^q_1$ of the graph $H_1$. Also, we identify one of the $|A^2_1|$ copies with the graph $H_1$.
Then, we connect the graph $H^q_1$ with the vertex $q$ by adding an edge $e$ between $q$ and the
connecting vertex of $H^q_1$. The resulting graph will be denoted by $H_2$.
So, $|V(H_2)|=|V(K_2)|+|A^2_1| |V(H_1)|\,.$
From now on, $p_2\in V(K_2)$ will also be called the connecting vertex of the graph $H_2$.
\vskip 0.05in
\noindent{\bf Step 2.}
Consider the graph $K_3$. Again,
for each vertex $r\in A^3_1$, we pick a copy $H^r_1$ of the graph $H_1$. Then, we
connect $H^r_1$ with the vertex $r$ as above. Also, for each $s\in A^3_2$, we pick
a copy $H^s_2$ of the graph $H_2$ and connect $H^s_2$ to the vertex $s$ via the connecting vertex.
Again, we identify one of the copies $H^s_2$ with the graph $H_2$. So, the resulting graph $H_3$ contains
$H_2$ as an induced subgraph.
From now on, $p_3\in K_3$ will also be called the connecting vertex of $H_3$.
\vskip 0.05in
\noindent{\bf Step n.}
Suppose that the graphs $H_1\subset H_2\subset\dots \subset H_n$ have already been
constructed and each graph $H_i$ contains a connecting point of degree $3$.
Now, for any $1\leq i \leq n$ and for each $w\in A^{n+1}_i$, we pick
a copy $H_i^w$ of the graph $H_i$. Then we connect $H^w_i$ to the vertex $w$ by an edge between $w$
and the connecting vertex of $H^w_i$. Finally, we identify for some $z\in A^{n+1}_n$ the graph $H^z_n$ with
the graph $H_n$. Finally, the connecting vertex of $K_{n+1}$ will be called the connecting vertex of $H_{n+1}$.

\noindent
Now, let $H=\cup^\infty_{i=1} H_i$. Observe that $H$ is a connected, infinite graph with vertex degree bound $4$.
\vskip 0.05in
\noindent
For each vertex $x\in V(H)$, there is a unique integer $j_1(x)$ such that $x$ is a vertex of a certain copy $K_{j_1(x)}(x)$ of
the  graph $K_{j_1(x)}$. Then, we have an integer $j_2(x)$ and a copy  $K_{j_2(x)}(x)$ of the
graph $K_{j_2(x)}$  such that the subgraphs $K_{j_1(x)}(x)$ and $K_{j_2(x)}(x)$ are
connected by an edge in the above process of building the graph
$H$. Inductively, we have an infinite sequence
of integers $j_1(x)<j_2(x)<\dots$ and disjoint induced subgraphs $K_{j_1(x)}(x), K_{j_2(x)}(x),\dots$ so that
the graphs $K_{j_n(x)}(x)$ and $K_{j_{n+1}(x)}(x)$ are connected to each other by an edge. We call the integer $j_1(x)$ the type
of $x$.
\begin{lemma} \label{csut1}
For all $x\in V(H)$ and $l\geq 1$, there exists $y\in V(H)$ such that
\begin{itemize} 
\item The type of $y$ equals to $l$,
\item $d_H(x,y)\leq \diam(H_l)+10^l+2$.
\end{itemize}
\end{lemma}
\proof
Suppose that $j_1(x)=m>l$. Then, $x\in K_m(x)$ and there exists
an element $z$ of the distinguished $l$-subset of the copy $K_m(x)$ such that
$d_H(x,z)\leq 10^l$. By our construction, $z$ is adjacent of a vertex $y$ of type $l$.
Now, assume that $j_1(x)\leq l$, Let $k\geq 1$ be the integer such that $j_k(x)\leq l$, $j_{k+1}(x)>l$.
By our construction, there exists a vertex $w$ of type $j_{k+1}(x)$ such that $d_H(x,w)\leq \diam(H_l)+1$. Hence,
by our previous observation, the lemma follows. \qed
\begin{proposition}
The graph $H$ is almost finite.
\end{proposition}
\proof
Let $l\geq 1$ and let $S_l$ be the set of all vertices $x\in V(H)$ such that $j_1(x)>l$.
Then, $S_l$ is the disjoint union of subsets $P$, where the induced subgraph $\overline{P}$ on $P$ is isomorphic to
$K_i$ for some integer $i>l$.
For each set $P$, let $A^P_l=\{\alpha^P_1,\alpha^P_2,\dots,\alpha^P_{k^P_l}\}$ be the
distinguished $l$-subset in $\overline{P}$. 
So, we can partition $P$ into subsets $\cup_{j=1}^{k^P_l} U_{p,j}$, where for any $1\leq j \leq k^P_L$,
$$\alpha^P_j\in U_{P,j}\subset B_{10^l}(\overline{P},\alpha^P_j)\,.$$
\noindent
Now we define the subset $V_{P,j}\subset V(H)$ in the following way. Let $x\in V_{P,j}$
if either $x\in U_{P,j}$ or $x$ is a vertex of a copy of $H_i$, $1\leq i\leq l$ attached to a vertex in the set $U_{P,j}\cap A^P_i$
in the construction. Observe that $\cup_P \cup _{1\leq j \leq k^P_l} V_{P,j}$ is a partition of $V(H)$.
Also,
\begin{equation}\label{egy131}
3^{10^l+1} l\leq |K_l|\leq |H_l|\leq |V_{P,j}|.
\end{equation}
\noindent
and
\begin{equation}\label{egy132}
|\partial(V_{P,j})|\leq |U_{P,j}|\leq  3^{10^l+1}\,.
\end{equation}
\noindent
Therefore, for all $P$ and integer $1\leq j \leq k^P_l$ we have that
\begin{equation} \label{egy133}
\frac{|\partial(V_{P,j})|}{|V_{P,j}|}\leq \frac{1}{l}\,.
\end{equation}
\noindent
Thus, by \eqref{egy131} and \eqref{egy133}, we can immediately see that the graph $H$ is almost finite. \qed

\subsection{A non-amenable almost finite groupoid}
Now we prove the main result of this section. It answers a query of
Suzuki (Remark 3.7 \cite{Suzuki}).
\begin{theorem} \label{nonamena}
There exists a stable minimal action $\alpha:\Gamma\actson\cC$ of finitely
generated group such that the associated minimal geometric groupoid $\cG_\alpha$
is almost finite but non-amenable.
\end{theorem}
\proof
Consider the connected graph $H$ constructed in the previous subsection.
Let $Q$ be a finite set and $\phi:E(H)\to Q$ be 
a labeling such that
adjacent edges have different labels.
Again, we consider the $|Q|$-fold free power $F_Q$ of the cyclic group
of two elements with symmetric generator set $Q$.
The labeling induces a transitive action
of the group $F_Q$ on the set $V(H)$ such that $H$ is the underlying
graph of the associated Schreier graph $S_H$.
For each $n\geq 1$, consider
an $(1/n, K_n)$-almost finite partition of $H$. Pick a
large enough integer $k_n$ and a labeling $\tau_n:V(H)\to\{0,1\}^{k_n}$
that encodes the partition. That is, if $x$ and $y$ belong to the same part
then $\tau_n(x)=\tau_n(y)$, and if $\tau_n(x)=\tau_n(y)$ for some $x,y$ which
are not belonging to the same part, then $d_H(x,y)>3K_n$.
Also, for each $n\geq 1$, pick a large enough integer $l_n$ and
a map $\kappa_n:V(H)\to \{0,1\}^{l_n}$ such that if $0<d_H(x,y)\leq n$, then 
$\kappa_n(x)\neq \kappa_n(y)$.
Finally, for $n\geq 1$, let $\mu_{2n}=\kappa_n$ and
$\mu_{2n-1}=\tau_n$. Consider
$h\in H$ and the $\cC$-labeling $\mu=\prod^\infty_{n=1}\mu_n:V(H)\to
\{0,1\}^\N$
and the element $y=(H,h,\mu)\in\cC \gsg$, where $\Gamma=F_\Q$ and $\Sigma=Q$.
Let $Z$ be a closed, minimal, invariant set in the orbit closure
$\overline{O}(y)$.
Our theorem immediately follows from the proposition below.
\begin{proposition} \label{szerdapro}
The restricted action $\alpha:\Gamma\actson Z$ is stable, almost finite and
the associated geometric groupoid $\cG_\alpha$ is minimal.
\end{proposition}
\proof
By the choice of the mappings $\kappa_n$, there exist continuous maps $\rho_n:Z\to \{0,1\}^{l_n}$
such that if $x$ and $y$ are vertices of the same graph $G$ in $Z$ and
$0<d_G(x,y)\leq n$, then $\rho_n(x)\neq \rho_n(y)$. Therefore, the 
action of $F_Q$ on $Z$ is stable.
Also, by the choice of the mappings $\tau_n$, there
exist continuous maps $\lambda_n:Z\to \{0,1\}^{l_n}$ that define
$(1/n,K_n)$-almost finite partitions of the elements of $Z$ as required in
Definition \ref{d2f26}. We only need to show that the geometric groupoid
$\cG_\alpha$
is non-amenable. By the definition of an amenable groupoid \cite{Anant}, if
$\cG_\alpha$ is amenable then all the orbit graphs of $\cG_\alpha$ possess
Yu's Property A. It is well-known that if an infinite graph $G$ contains
a growing sequence of expanders as induced subgraphs, then $G$ does not
have Property A. Our graph $H$ was constructed in such a way
that it contains a sequence of expanders as induced subgraphs.
By Lemma \ref{csut1}, for any $n\geq 1$, there exists
$t_n>0$ such that if $x\in V(H)$, then the ball $B_{t_n}(H,x)$ contains
a copy of $K_n$ as induced subgraph. Hence, if $G$ is in the orbit closure
of $H$ and $y\in V(G)$, then  the ball $B_{t_n}(G,y)$ contains an induced
copy of $K_n$. Therefore, the orbit graphs in $\cG_\alpha$ do not have
Property $A$. Hence our proposition (and so the theorem) follows. \qed
\subsection{Fractionally almost finite graphs}
The notion of fractional almost finiteness (see Definition \ref{weaklyalm}) is
closely
related to the notion of fractional hyperfiniteness introduced by Lov\'asz
\cite{Lovhyp}.
Observe that if $G\in\grd$ is a fractionally almost finite graph and $H\subset G$ is
a subgraph, then $H$ is fractionally almost finite as well. On the other hand,
almost finiteness of $G$ does not necessarily imply the almost finiteness of
$H$, since the Cayley graph of some amenable groups contain copies of the
$3$-regular tree as a subgraph.
Lov\'asz proved that fractional hyperfiniteness implies hyperfinitess for
graphing. Nevertheless, we have the following proposition.
\begin{proposition}
The $k$-regular tree $T_k$ is fractionally almost finite (and, of course,
$T_k$ is not almost finite).
\end{proposition}
\proof Fix an infinite ray
$\{p_1,p_2,\dots\}$
in $T_k$. Let $l\geq 1$ be an integer. Now, for each integer
$1\leq i \leq l$ we construct a partition of $V(T_k)$. For $x\in V(T_k)$, let
$P_i(x)$ be the shortest path from $x$ to $p_i$. Let $x\equiv_{l,i} y$ if
there exists sonme $z\in V(T_k)$ such that
\begin{itemize}
\item both $P(x)$ and $P(y)$ contains $z$,
\item $d_{T_k}(x,z)<l$, $d_{T_k}(y,z)<l$ and
\item $d_{T_k}(z,p_i)$ is divisible by $l$.
\end{itemize}
\noindent
Then,
\begin{enumerate}
\item $\equiv_{l,i} $ is an equivalence relation.
\item The equivalence classes have bounded diameter.
\item The vertex $x$ is on the boundary of its class if and only if \begin{equation}\label{penteq}
d_{T_k}(x,p_i)\equiv 0\,\mbox{or}\, l-1\,(\mbox{mod}\, l)\,.\end{equation}
\end{enumerate}
Observe that for any $x\in V(G)$, there exists at most $2$ elements $i$
of the set $\{1,2,\dots, l\}$ such that \eqref{penteq} holds.
Hence, if $l>\frac{2}{\eps}$, then \eqref{este1} is satisfied. \qed
\begin{remark}
Using Proposition 2.10 of \cite{DG}, it is not hard to see that all
fractionally almost finite graphs have Property A. It would be interesting to
construct a Property A graph $G$ that is not fractionally almost finite.
\end{remark}

\section{The spectra of graphs}
The main goal of this section is to prove a spectral convergence result for
strongly almost finite graph classes.
\subsection{Uncountably many isospectral connected graphs}
Let $G\in\grd$ be an infinite graph and
$\cal{L}_G:l^2(V(G))\to l^2(V(G))$
\noindent
be the Laplacian operator on $G$.
That is, $$\cal{L}_G(f)(x)=\deg(x) f(x)-\sum_{x\sim y} f(y)\,.$$
\noindent
It is well-known that $\cl$ is a positive, self-adjoint operator and $\Spec(\cl_G)\subset [0,2d]$.
\begin{proposition}
If $G$ and $H$ are equivalent (see Section \ref{tgc}), then \\
$\Spec(\cl_G)=\Spec(\cl_H)$.
\end{proposition}
\proof
First, we need a lemma.
\begin{lemma} \label{spec1}
Let $P$ be real polynomial, then
$\|P(\cl_G)\|=\|P(\cl_H)\|$.
\end{lemma}
\proof Fix some $\eps>0$.
Let $f\in l^2(V(G)$ such that $\|f\|=1$ and
$\|P(\cl_G)(f)\|\geq (1-\eps) \|P(\cl_G)\|$. We can assume that $f$ is supported
on a ball $B_s(G,x)$ for some $s>0$ and $x\in V(G)$. Let $t$ be the degree of $P$.
Then, $P(\cl_G)(f)$ is supported in the ball $B_{s+t}(G,x)$.
Since $G$ and $H$ are equivalent, there exists $y\in V(H)$ such that the ball
$B_{s+t}(G,x)$ is rooted-isomorphic to the ball $B_{s+t}(H,y)$ under the rooted-isomorphism $j$.
Then, $\|j_*(f)\|=1$ and $\|P(\cl_G)(f)\|= \|P(\cl_H)(j_*(f))\|$, where
$j_*(f)(z)=f(j^{-1}(z)).$
Therefore, $\|P(\cl_H)\|\geq (1-\eps)\|P(\cl_G)\|$ holds for any $\eps>0$. Consequently,
 $\|P(\cl_H)\|\geq \|P(\cl_G)\|$. Similarly,  $\|P(\cl_G)\|\geq \|P(\cl_H)\|$, thus our lemma follows. \qed
\vskip 0.1in
\noindent
By functional calculus, we have that
\begin{equation} \label{eqspec}
\|f(\cl_G)\|=\|f(\cl_H)\|
\end{equation}
\noindent
holds for any real continuous function.
Observe that $\lambda\in\Spec(\cl_G)$ if and only if
for any $n\geq 1$ $\|f_n^\lambda(\cl_G)\|\neq 0$, where $f^\lambda_n$ is a piecewise linear, continuous, non-negative
function such that
\begin{itemize}
\item $f^\lambda_n(x)=1$ if $\lambda-\frac{1}{n}\leq x \leq \lambda + \frac{1}{n}\,.$
\item $f_n^\lambda(x)=0$ if $x\geq \lambda +\frac{2}{n}$ or $x\leq \lambda-\frac{2}{n}\,.$
\end{itemize}
\noindent
Therefore, by \eqref{eqspec} our proposition follows. \qed
\vskip 0.1in
\noindent
It is well-known that many isospectral finite graphs exist. It is less-known
(but certainly known for experts) that many isospectral infinite connected
graphs exist, so the following proposition might be interesting on its own.
\begin{proposition}
For any $d>3$, there exist uncountable many pairwise non-isomorphic trees
in $\grd$ which possess the same spectra.
\end{proposition}
\proof
Let $K_d$ be the set of all finite trees in $\grd$. 
Let $T$ be an infinite tree of vertex degree bound $d-1$.
We construct a tree $\hat{T}$ in the following way.
Let $\{t_i\}^\infty_{i=1}$ be an enumeration of the vertices of $T$ and $\{T_i\}^\infty_{i=1}$ be an enumeration
of the elements of $K_d$. For each $i\geq 1$ let us connect the tree $T_i$ with $T$ by adding an edges
between $t_i$ and a vertex $p_i\in T_i$ such that $\deg(p)<d$.
\begin{lemma} \label{spec2lemma}
The set of connected, finite induced subgraphs of $\hat{T}$ (up to isomorphisms) is exactly $K_d$.
\end{lemma}
\proof
By our construction, if $G\in K_d$, then $G$ is an induced subgraph of $\hat{T}$.
Conversely, let $H$ be a connected, finite induced subgraph of
$\hat{T}$. Clearly, $H$ is a tree of vertex degree bound $d$. \qed
\vskip 0.1in
\noindent
By our lemma, for any $T,S$ with vertex degree bound $d-1$, we have that $\hat{T}$ is equivalent to $\hat{S}$.
It is easy to check that if $T$ and $S$ are infinite non-isomorphic trees without leaves, then $\hat{T}$ and
$\hat{S}$  are not isomorphic. Since there are uncountably many such trees,
our proposition follows. \qed
\begin{remark}
Note that the end space of the tree $T$ and $\hat{T}$ are the same, so one can
actually has uncountably many trees with the same spectrum and pairwise
different end spaces.
\end{remark}
\subsection{Spectral convergence for strongly almost finite graphs}
The main goal of this section is to prove the following theorem.
\begin{theorem} \label{spectralconv}
Let $\{G_n\}^\infty_{n=1}\subset \grd$ be a strongly almost finite family of
graphs such that $G_n\stackrel{n}{\to} G$. Then $\Spec{G_n}\to \Spec(G)$ in the
Hausdorff topology.
\end{theorem}
\proof
Let $G$ be a strongly almost finite infinite graph. Then, it is not hard to
see
that for any $\delta>0$ and integer $l>1$ there exist constants
$M_{\delta,l}>0, L_{\delta,l}>0$, and finite partitions
$$\{H^i_1, H^i_2,\dots\}_{i=1}^{M_{\delta,l}}$$
\noindent
such that
\begin{itemize}
\item the diameters of the classes $\{H^i_j\}$ are bounded by $L_{\delta,l}$,
\item for each class $H^i_j$, $\frac{|\partial_l (H^i_j)|}{|H^i_j|}<\delta$,
\item for each $x\in V(G)$, 
$$\frac{|\{i\,\mid\, x\in \partial_l (H^i_j)|\,\mbox{
  for some}\, j\geq 1\}|}{M_{\delta,l}}<\delta.$$
\end{itemize}
\noindent
We call such a partition family $\cal{E}$ a $(\delta, l, M_{\delta,l} L_{\delta,l})$-almost finite partition family.
Now, let $f:V(G)\to\R$ be a function such that $\sum_{v\in V(G)} f^2(v)=1$.
Let $\cal{E}$ be a $(\delta, l, M_{\delta,l} L_{\delta,l})$-almost finite partition family and for $1\leq i \leq M_{\delta,l}$ let
$f_i:V(G)\to\R$
be defined by setting
\begin{itemize}
\item $f_i(x)=f(x)$ if $x\in H^i_j\backslash \partial_l (H^i_j)$ for some
  $j\geq 1$.
\item $f_i(x)=0$, otherwise.
\end{itemize}
\begin{lemma}\label{szombat1}
We have the following inequality.
$$|\{i\,\mid\, \|f_i\|^2< (1-\sqrt{\delta})\}|<\sqrt{\delta} M_{\delta,l}.$$
\end{lemma}
\proof
By our condition,
$$\sum_{x\in V(G)} \sum_{1\leq i \leq M_{\delta,l}} f^2_i(x)\geq  
\sum_{x\in V(G)} (1-\delta)  M_{\delta,l}f^2(x).$$
\noindent
That is,
$$\sum_{1\leq i \leq M_{\delta,l}} \|f_i\|^2\geq (1-\delta) M_{\delta,l}.$$
\noindent
Let $A$ be the set of integers $1\leq i \leq M_{l,\delta}$ for which
$\|f_i\|^2< (1-\sqrt{\delta}).$
Then,
$$|A|(1-\sqrt{\delta})+ ( M_{l,\delta}-|A|)\geq (1-\delta)  M_{l,\delta}\,.$$
\noindent
That is, 
$$\frac{|A|}{M_{\delta,l}} (1-\sqrt{\delta})+ (1-\frac{|A|}{M_{\delta,l}}
)\geq (1-\delta)\,,$$
Hence, $\frac{|A|}{M_{\delta,l}}<\sqrt{\delta}$. \qed \vi
Now, let $\cal{E}$ as above and $P$ be a polynomial of degree less than $l$.
Define $$\|P(\Delta)\|_{\cal{E}}:=\sup_g \frac{\|P(\Delta)g\|}{\|g\|}\,,$$
\noindent
where the supremum is taken for
all nonzero functions $g$ which are supported on $H^i_j\backslash\partial_l( H^i_j)$ for
some class $H^i_j$.
Clearly, $\|P(\Delta)\|_{\cal{E}}\leq \|P(\Delta)\|.$
Also, let 
$$\|P(\Delta)\|_{\hat{\cal{E}}}:=\sup_g \frac{\|P(\Delta)g\|}{\|g\|}\,,$$
\noindent
where the supremum is taken for 
all $L^2$-functions $g$ such that there exists $i$ for which $g(x)=0$ if $x\in \partial( H^i_j)$ for
all
classes $H^i_j$. By the degree condition on $P$, $\|P(\Delta)\|_{\cal{E}}=
\|P(\Delta)\|_{\hat{\cal{E}}}$.
\begin{lemma} \label{szomb2}
Suppose that $\max_{0\leq t \leq 2d} |P(t)|\leq 2$. Then,
$$\|P(\Delta)\|-\|P(\Delta)\|_{\cal{E}}< 3\left(1-\sqrt{1-\sqrt{\delta}}\right).$$
\end{lemma}
\proof
Let $f:V(G)\to \R$ such that
$\|f\|=1$ and $\|P(\Delta)f\|\geq \|P(\Delta\|-(1-\sqrt{1-\sqrt{\delta}})\,.$
Let $f_i$ be as above such that
$\|f_i\|>(1-\sqrt{\delta})$. Then,
$$\|P(\Delta)f_i\|\geq \|P(\Delta) f\|-2(1-\|f_i\|)\geq
\|P(\Delta)\|-3\left(1-\sqrt{1-\sqrt{\delta}}\right).$$
\noindent
Since $\|f_i\|\leq 1$, we have that
$$\|P(\Delta)\|_{\hat{\cal{E}}}\geq
\|P(\Delta)\|-3\left(1-\sqrt{1-\sqrt{\delta}}\right).$$
\noindent
Thus, our lemma follows. \qed
\vi
For any $m\geq 1$, we have polynomials, $\{P^m_i\}_{i=1}^{m2d}$ such that
\begin{itemize}
\item for any $0\leq t \leq 2d$, $0< P^m_i(t)\leq 1+\frac{1}{m},$
\item $0< P^m_i(t)<\frac{1}{m}$ if $t\leq \frac{i}{m2d}-\frac{1}{m}$ or $t\geq
  \frac{i}{m2d}+\frac{1}{m}$,
\item $1-\frac{1}{m}< P^m_i(t)<1 +\frac{1}{m}$, if $\frac{i}{m2d}-\frac{1}{2m}<t
  < \frac{i}{m2d}+\frac{1}{2m}.$
\end{itemize}
\vi
Now, let $G_n\stackrel{n}{\to} G$, where
$\{G_n\}^\infty_{n=1}$ is a strongly almost finite class. Thus, by
Proposition \ref{limitalmost} $\{\|G_n\|^\infty_{n=1},G\}$ is a strongly
almost finite class as well. By the previous lemma, for any $\delta>0$ and
$l>0$, we have $M_{\delta,l}>0, L_{\delta,l}>0$ and
$(\delta,l,M_{\delta,l},L_{\delta,l})$-almost finite partition families
$\cal{E}_n$ on $G_n$
and a $(\delta,l,M_{\delta,l},L_{\delta,l})$-almost finite partition family
$\cal{E}$ on $G$ such that
for all $n\geq 1$ and $1\leq i \leq m2d$,
$$\|P^m_i(\Delta_{G_n})\|-\|P(\Delta_{G_n})\|_{\cal{E}_n}\leq
3\left(1-\sqrt{1-\sqrt{\delta}}\right),$$
$$\|P^m_i(\Delta_{G})\|-\|P(\Delta_{G})\|_{\cal{E}}\leq
3\left(1-\sqrt{1-\sqrt{\delta}}\right).$$
\noindent
Since $G_n\tart G$, for large enough $n$ the set of all induced subgraphs of
diameter at most $L_{\delta,l}$ in $G_n$ equals to the set of
all induced subgraphs of
diameter at most $L_{\delta,l}$ in $G$.
Therefore, for large enough $n$, we have that for all $1\leq i \leq m2d$,
\begin{equation} \label{szombegyenlet}
\left| \|P^m_i(\Delta_{G_n})\|-\|P^m_i(\Delta_G)\|\right|<
6\left(1-\sqrt{1-\sqrt{\delta}}\right).
\end{equation}
\noindent
Suppose that $\lambda\in\Spec(\Delta_G)$. Then by \eqref{szombegyenlet}, for any $\zeta>0$
there exists $N_\zeta>0$ such that if $n>N_\zeta$, then we have some
$\lambda_n\in\Spec(G_n)$, $|\lambda_n-\lambda|<\zeta.$
Also, there exists $M_\zeta>0$ such that if $n> M_\zeta$ and
$\mu_n\in\Spec(G_n)$, then there exists $\mu\in\Spec(G)$ such that
$|\mu_n-\mu|<\zeta$. Hence, our theorem follows. \qed
\begin{remark}
Let $G_n\stackrel{n}{\to} G$, where $\{G_n\}^\infty_{n=1}$ is a large girth
sequence
of $3$-regular graphs and $G$ is the $3$-regular tree. Then for all $n\geq 1$,
$0\in\Spec{G_n}$ and $0\notin\Spec{G}$.
\end{remark}
\section{Constant-time distributed algorithms}\label{CTDA}
\noindent The theory of distributed graph algorithms is a vast subject
developed in the last thirty years (see the monograph \cite{Lynch}). Using the
ideas of our paper we introduce the notion of constant-time distributed
algorithms,
the qualitative analogue of constant-time randomized local algorithms.
\subsection{Algorithms and oracles}
First, let us give a formal definition
for an $r$-round distributed algorithm on a graph class $\cG\subseteq \grd$.
For a finite set $Q$ and a graph $G$, we will denote by
$L_Q(G)$ the set of all functions $\phi:V(G)\to Q$.
An $r$-round distributed algorithm can be represented by an operator
$$O:L_{Q_1}(G)\to L_{Q_2}(G)$$
\noindent
for some finite sets $Q_1,Q_2$.
However, the value $O(\phi)(x)$ depends only on
the values of $\phi$ on the ball $B_r(G,x)$. Therefore,
the operator $O$ can be described by a map $\Theta:U^{r,Q_1}_d\to Q_2$, where
$$O(\phi)(x)=\Theta(B_1(G,x,\phi))\,.$$
\noindent
We call $\Theta$ an {\bf oracle function} and we use the
notation $O_\Theta$ for the associated algorithm operator.
The following examples formalize how to combine simple algorithms into more complicated ones.
\begin{example} \label{ex3f1}
Let $\Theta_1: U^{r,Q_1}_d\to Q_2$ and $\Theta_2: U^{r,Q_2}_d\to Q_3$ be oracles.
Then, there exists a unique composition oracle $\Theta_3: U^{r+s,Q_1}\to Q_3$
such that
for any $\phi\in L_{Q_1}(G)$,
$$O_{\Theta_3}(\phi)=O_{\Theta_2}(O_{\Theta_1}(\phi))\,.$$
\end{example}
\begin{example} \label{ex3f2}
Let $\{\Theta_i:U_d^{r_i,M_i}\to Q_i\}^k_{i=1}$ be oracles.
Also, for the functions $\{\phi_i: L_{M_i}(G)\}^k_{i=1}$
let $\oplus^k_{i=1} \phi_i\in L_{\oplus_{i=1}^k M_i}(G)$ be the function
such that 
$$(\oplus^k_{i=1} \phi_i)(x)=\oplus^k_{i=1} (\phi_i(x)).$$
\noindent
Let $\Theta:U_d^{s,\oplus^k_{i=1}Q_i}\to N$ be an oracle and $r=\max_{1\leq i \leq k} r_i$.
Then, there exists a unique oracle $\Theta':U_d^{r+s, \oplus_{i=1}^k M_i}\to N$, such that
$$\Theta'(\oplus^k_{i=1} \phi_i)=\Theta(\oplus^k_{i=1} \Theta_i(\phi_i))\,.$$
\end{example}
\noindent
Now, let $k\geq 1$ be an integer and $Q$ be a finite set such that
$|Q|\geq |V(B)|$ (we say that $Q$ is $(r,d)$-large) for all $B\in U^k_d$. Let $L_Q^{(k)}(G)$ denote the set
of all $k$-separating functions $\phi:V(G)\to Q$. That is, functions $\phi$ for which
$\phi(x)\neq \phi(y)$ if $0<d_G(x,y)\leq k$.
As in the Introduction, a constant-time distributed algorithm (CTDA) on a graph class $\cG\subset \grd$,
takes a $k$-separating $Q_1$-labeling on $G\in\cG$ and produces some $Q_2$-labeling.
So, the algorithm operator is given by an oracle,
$$\Theta:U^{r,k,Q_1}_d\to Q_2\,,$$
\noindent
where $U^{r,k,Q_1}_d$ denote the set of all $k$-separating $Q_1$-labelings of the balls in $U^r_d$. 
Before getting further, let us consider a simple, but important example.
\begin{proposition}\label{p3f3}
We have a CTDA, which produces a maximal independent subset for any graph $G\in\grd$.
\end{proposition}
\proof
First, let us make the statement of our proposition completely precise.
Let $M\geq d+1$. Then, there exists an oracle
$\Theta:U^{|M|,1,M}_d\to \{a,b\}$
such that
for any graph $G\in\grd$ and $\phi\in L_M^{(1)}(G)$, the set
$(O_\Theta (\phi))^{-1}(a)$ is a maximal independent subset in $G$.
We will give a combinatorial description of $\Theta$ as
the composition of $|M|$ $1$-round oracles.
Let $\phi\in L_M^{(1)}(G)$. In the first round, we relabel
the vertices $x\in V(G)$ for which $\phi(x)=1$ by $a$ to obtain the labeling
$\phi_1:V(G)\to M\cup a$. Note that $\phi_1(y)=\phi(y)$ if $\phi(y)\neq 1$.
In the second round, we relabel the vertices $x$ for which $\phi_1(x)=2$.
Let $\phi_2(x)=a$ if $\phi_1(y)\neq a$ for all vertices $y$ adjacent to $x$.
Otherwise, let $\phi_2(x)=b.$ So, we obtain the labeling
$\phi_2:V(G)\to M\cup\{a,b\}$.
Inductively, in the $k$-th round we start with a labeling
$\phi_{k-1}:V(G)\to M\cup\{a,b\}$. Then, we relabel the vertices $x$ for which
$\phi_{k-1}(x)=k.$ Again, let $\phi_k(x)=a$ if $\phi_{k-1}(y)\neq a$ for all vertices $y$ adjacent to $x$.
Otherwise, let $\phi_k(x)=b.$ So, after the $|M|$-th round we obtain the labeling $\phi_{|M|}$ which maps
all vertices of $G$ into the set $\{a,b\}$ and $\{(\phi_{|M|})^{-1}(a)$ is a maximal independent set. Indeed,
if $\phi_{|M|}(x)=b$, then $\phi_{|M|}(y)=a$ for at least one vertex $y$ adjacent to $x$. \qed
\vskip 0.1in
Certain algorithms, e.g. the one described in Example \ref{ex3f2}, requires elements of the set
$L^{(r)}_M(G)$ as input. However, the input functions given are from the much 
larger set $L^{(r)}_N(G)$, where $|N|>|M|$. The following lemma shows that
we can convert $r$-separating $N$-valued functions into $r$-separating $M$-valued functions
using a simple oracle.
\begin{lemma} \label{convert}
Let $M=\{1,2,\dots,m\}, N=\{1,2,\dots,n\}$ be finite sets such that $n>m$ and $M$ is $(r,d)$-large. 
Then, we have an oracle
$$\Theta: U^{r(n-m),r,N}_d \to M\,,$$
\noindent
such that $O_\Theta$ maps $L_N^{(r)}(G)$ into $L_M^{(r)}(G)$.
\end{lemma}
\proof Let $\phi\in L_N^{(r)}(G).$ In the first round we relabel the elements
$x\in V(G)$, for which $\phi(x)=m+1$. Let $\phi_1(x)=i$, where $1\leq i \leq m$ is
the smallest integer such that $\phi(y)\neq i$, provided that $0\leq d_G(x,y)\leq r$.
By inductive relabelings, we construct the sequence $\phi_1,\phi_2,\dots,\phi_{n-m}$ in $n-m$ rounds.
Then, $\phi_{n-m}$ will be an $r$-separating $M$-labeling. \qed

\subsection{Distributed graph partitioning and almost finiteness}
We start with the notion of a  distributed graph partitioning oracle, the qualitative
analogue of the randomized partioning oracles introduced by Hassidim et. al. \cite{HKNO}.
Let $\cG\subseteq \grd$ be a class of graphs, $1\leq s \leq n$ and
$\Theta:U_d^{n,s,M}\to Q$ be an oracle, where $M$ is some $(s,d)$-large set.
Let $\eps>0$ be a real number and $K_\eps>0$ be an integer.
\begin{definition}\label{partoracle}
$\Theta$ is a {\bf distributed $(\eps, K_\eps)$-partitioning oracle} for the class $\cG$, if
for all $G\in\cG$ and $\phi\in L_M^{(s)}(G)$ the following two conditions are satisfied.
\begin{enumerate}
\item For any $x,y\in V(G)$, if $O_\Theta(\phi)(x)=O_\Theta(\phi)(y)$, we have that
either $d_G(x,y)\leq K_\eps$ or $d_G(x,y)\geq 3 K_\eps$.
\item For any $x\in V(G)$, $i_G(H_x)\leq \eps$, where
$$H_x=\{z\in V(G)\,\mid\, d_G(x,z)\leq K_\eps\,\mbox{and}\,O_\Theta(\phi)(x)= O_\Theta(\phi)(z)\}\,.$$
\end{enumerate}
\end{definition}
\noindent
Let $x\equiv_{\langle\Theta,\phi\rangle} y$ if $y\in H_x$. It is easy to see
that the relation $\equiv_{\langle\Theta,\phi\rangle}$ is, in fact, an equivalence relation. So, $\Theta$ computes
an $(\eps,K_\eps)$-partition, indeed.
\begin{definition}\label{distalm}
A family $\cG\subset \grd$ is distributed almost finite if
for any $\eps>0$, there exist integers $1 \leq s_\eps\leq n_\eps $, and
$K_\eps\geq 1$, finite sets $M_\eps, Q_\eps$ such that $M_\eps$ is
$(s_\eps,d)$-large and we also have an oracle
function $\Theta_\eps:U^{n_\eps,s_\eps,M_\eps}_d\to Q_\eps$ such that $\Theta_{\eps}$ is a $(\eps, K_\eps)$-partitioning oracle for $\cG$.
\end{definition}
\noindent
Let $\cG\in\grd$ be a  class of finite graphs, $1\leq s\leq n$ be integers and $M$ be an $(s,d)$-large set.  Again, let
$Q$ be another finite set and $*$ be an extra symbol.
Let $\Theta:U^{n,s,M}_d\to \{Q,*\}$ be an oracle function. We say that
$\Theta$ is an {\bf $(\eps,K_\eps)$-hyperfinite partitioning oracle} for the class $\cG$ if the following three conditions are satisfied.
\begin{enumerate}
\item For any $x,y\in V(G)$ and $\phi\in L^{s}_M(G)$ if $*\neq O_\Theta(\phi)(x)=O_\Theta(\phi)(y)$, we have that
either $d_G(x,y)\leq K_\eps$ or $d_G(x,y)\geq 3 K_\eps$.
\item For any $x\in V(G)$, such that  $O_\Theta(\phi)(x)\neq *$,  $i_G(H_x)\leq \eps$, where 
$$H_x=\{z\in V(G)\,\mid\, d_G(x,z)\leq K_\eps\,\mbox{and}\,O_\Theta(\phi)(x)= O_\Theta(\phi)(z)\}\,.$$
\item $|\{z\,\mid\,O_\Theta(\phi)(z) =*\}|\leq \eps |V(G)|\}\,.$
\end{enumerate}
\begin{definition}
A family of finite graphs $\cG\subset \grd$ is distributed hyperfinite if
for any $\eps>0$, there exist integers $1 \leq s_\eps \leq n_\eps$ and $K_\eps\geq 1$, finite sets $M_\eps, Q_\eps$ such that $M_\eps$ is $(s_\eps,d)$-large and an oracle
function $\Theta_\eps:U^{n_\eps,s_\eps,M_\eps}_d\to \{Q_\eps,*\}$ such that $\Theta_{\eps}$ is a $(\eps, K_\eps)$-hyperfinite partitioning oracle for $\cG$.
\end{definition}
\noindent
Finally, we define distributed strong almost finiteness.
\begin{definition}
A family $\cG\subset \grd$ is distributed strongly almost finite if
for any $\eps>0$, there exist integers $1\leq s_\eps\leq n_\eps $, and $N_\eps,K_\eps\geq 1$, finite sets $M_\eps, Q_\eps$  such that $M_\eps$ is $(s_\eps,d)$-large , and  oracle
functions $\{\Theta_i:U^{n_\eps, s_\eps,M_\eps}_d\to Q_\eps\}^{N_\eps}_{i=1}$ such that
\begin{enumerate}
\item 
$\{\Theta_i\}^{N_\eps}_{i=1}$ are $(\eps,K_\eps)$- partitioning oracles for $\cG$
\item for all $G\in\cG$, $\phi\in L^{s_\eps}_M(G)$  and $x\in V(G)$,
$$\frac{|\{i\,\mid x\in \partial_G(H_x^i)\}|}{|{N_\eps}|}\leq \eps\,,$$
\noindent
where $H^i_x$ is the class defined by the labeling $O_{\Theta_i}(\phi)$.
\end{enumerate}
\end{definition}
\noindent
The following theorem is the qualitative analogue of the main result of \cite{HKNO}.
\begin{theorem}\label{otostetel}
Let $\cG\in \grd$ be a hyperfinite family of finite graphs. Then, $\cG$ is
distributed hyperfinite as well.
\end{theorem}
\proof First, we need two lemmas.
\begin{lemma}\label{l201}
Let $G\in\grd$ be a hyperfinite family of finite graphs. Assume that for any
$0<\eps<1$ and $G\in\cG$, $G$ is $(\eps,K_\eps)$-hyperfinite. Then, provided
that $G\in\grd$ and $H\subset G$ is an induced subgraph such that
$|V(H)|\geq\eps |V(G)|$, one can remove $(d+1)\eps |V(H)|$ vertices from $H$ (with all
the adjacent edges) so that for all the remaining components $T$, $|T|\leq K_{\eps^3}$ and $i(T)\leq\eps$.
\end{lemma}
\proof Observe, that $H$ is $(\eps^2,K_{\eps^3})$-hyperfinite. So, let us
remove $\eps^2 |V(H)|$ vertices from $H$ in such a way that
all the remaining components have size at most $K_{\eps^3}$. Let $\cA$ be the
subset of the remaining components $M$ such that $i_H(M)>\eps$. Then, we have
that
$$\eps\sum_{M\in\cA} |M|\leq d\eps^2 |V(H)|\,.$$
\noindent
Hence, $\sum_{M\in\cA} |M|\leq d\eps |V(H)|$. Thus, we can remove
$(d\eps+\eps^2) |V(H)|$ vertices from $V(H)$ in such a way that all the
remaining components have size at most $K_{\eps^3}$ and have isoperimetric constant
at most $\eps$. \qed
\begin{lemma}\label{l202}
Let $(1-(d+1)\eps)>1/2$. Let $H\subset G$ be graphs as above. Suppose that
$\cM$ is a maximal system of connected induced subgraphs in $H$ such that
if $C\neq D\in \cM$, then $d_H(V(C),V(D))\geq 2$ and if $C\in \cM$, then
$|V(C)|\leq K_{\eps^3}$ and $i(V(C))\leq \eps$. Then,
$$\sum_{C\in\cM} |V(C)|\geq \frac{\eps}{4d^2 K^2_{\eps^3}}|V(G)|\,.$$
\end{lemma}
\proof
By the previous lemma, we have a system of induced subgraphs $\cN$ such that
\begin{itemize}
\item $\sum_{A\in\cN} |V(A)|\geq (1-(d+1)\eps)|V(H)|,$
\item if $A\neq B\in \cN$, then $d_H(V(A),V(B))\geq 2$,
\item if $A\in\cN$, then $|V(A)|\leq K_{\eps^3}$ and $i(V(A))\leq \eps$.
\end{itemize}
\noindent
Notice that if $A\in\cN$, then there exists an element $C\in\cM$ such
that $V(A)\cap B_2(H,V(C))\neq\emptyset$, that is, the $2$-neighbourhood of the set $V(C)$ in $H$
intersects $V(A)$. Indeed, if such subgraph $C$ did not exist, then $\cM$
could not be a maximal system.
Since $|B_2(H,V(C))|\leq 2d^2 K_{\eps^3}$, we have that
$$ |\cM|\geq \frac{1}{2d^2 K_{\eps^3}} |\cN|\,.$$
\noindent
Also, 
$$\sum_{A\in\cN} |V(A)|\geq (1-(d+1)\eps)|V(H)|\geq \frac{1}{2} |V(H)|.$$
\noindent
Therefore $|\cN|\geq \frac{1}{2 K_{\eps^3}}\eps |V(G))|.$ Thus,
$$\sum_{C\in\cM} |V(C)|\geq |\cM|\geq \frac{\eps}{4d^2 K^2_{\eps^3}}
  |V(G)|\quad\qed $$
\noindent
Now, we can prove our theorem. Repeating the argument of Proposition \ref{p3f3}, we can
devise an oracle $\Theta: U_d^{n,M\times \{a,b\}}\to \{c,d\}$ such that for any pair $\phi\in L^{(n)}_M(G)$ and $\psi\in L_{\{a,b\}}(G)$,
$O_{\Theta}(\phi\oplus\psi)=\rho\in L_{\{c,d\}}(G)$ and $\cM$ is a maximal system in the induced subgraph $H$ satisfying the following two conditions.
\begin{enumerate}
\item If $C\neq D\in\cM$, then $d_H(V(C),V(D))\geq 2.$
\item For any $C\in\cM$, $|V(C)|\leq K_{\eps^3}$ and $i(V(C))\leq\eps$,
\end{enumerate}
\noindent
where $H$ is the induced subgraph on the set $\psi^{-1}(a)$ and $\cM$ is the induced subgraph on the set $\rho^{-1}(c)$.
Now, we start with the graph $G$ and set $H_0=G$. We apply the 
distributed algorithm $O_\Theta$
to compute the maximal system $\cM_0$. Then, we remove the vertices covered by
the elements of $\cM_0$ and all the neighbouring vertices to obtain the
induced subgraph $H_1$. Now, we remove the vertices of the maximal system
$\cM_1$ in $H_1$ together with the neighbouring vertices to obtain $H_2$ and
so on. Hence, we compute a sequence of induced subgraphs $H_0\supset H_1\supset\dots \supset H_j$, 
where $j$ is an integer larger than $\frac{4d^2  K^2_{\eps^3}}{\eps}.$   By Lemma \ref{l202}, we have that 
$|V(H_j)|\leq \eps |V(G)|.$  Now, let $X$ be the union of all vertices
computed in the process as vertices covered by some maximal system. Let $Y$ be
the neighbours of the elements of $X$.
Then, 
$$V(G)=X \cup Y \cup V(H_j)\,.$$
Also, $|Y|\leq d\eps |X|$, since the elements of the maximal systems computed
in the process are graphs $C$ such that $i(V(C))\leq \eps$.
Therefore, 
$$|Y\cup V(H_j)|\leq (d+1) \eps |V(G)|\,.$$
\noindent
Hence, we have a CTDA that computes a set $Z$ of size
at most $(d+1)\eps |V(G)|$ in the graph $G$ such that if we remove
$Z$ from $V(G)$, all the remaining components have size at most
$K_{\eps^3}$. Thus,
$\cG$ is a distributed hyperfinite family of finite graphs. \qed

\subsection{Approximating maximum independent subsets}
\begin{definition}
Let $\cG\subset \grd$ be a class of finite graphs. We say that there exist CTDA's for the approximated maximum independent subset problem in $\cG$,
if for any $\eps>0$, there exist integers $1\leq s_\eps\leq n_\eps$, a finite $(s_\eps,d)$-large subset $M_\eps$ and an oracle $\Theta:U^{n_\eps,s_\eps,M}_d\to \{a,b\}$
satisfying the following properties.
\begin{enumerate}
\item For any $G\in\cG$ and $\phi\in L_M^{(s)}(G)$, $(O_\Theta (\phi))^{-1}(a)$ is an independent
subset in $G$.
\item $\frac{|(O_\Theta (\phi))^{-1}(a)|}{|V(G)|} \geq \frac{|I_G|}{|V(G)|}-\eps$, where $I_G$ is a maximum size independent subset in $G$.
\end{enumerate}
\end{definition}
\begin{proposition} \label{hypind}
Let $\cG\subset \grd$ be a hyperfinite class of finite graphs. Then, there exist CTDA's for approximated maximum independent sets in $\cG$.
\end{proposition}
\proof
Let $G\in\grd$, $Q$ be a finite set and $K\geq 1$ be an integer.
Then, we call a function $\rho:V(G)\to \{Q,*\}$ a $K$-subset function if the 
following two conditions are satisfied.
\begin{enumerate}
\item If $\rho(x)=\rho(y)$, then either $d_G(x,y)\leq K$ or $d_G(x,y)\geq 3K$.
\item If $\rho(x)\neq *$, $\rho(z)\neq *$ and $\rho(x)\neq \rho(y)$, then
$d_G(x,y)\geq 2$.\
\end{enumerate}
\noindent
For such a $K$-subset function $\rho$ and $x\in V(G)$ and $\rho(x)\neq *$,  we will denote by $T^\rho_x$
the induced subgraph on the set
$$\{z\,\mid\, d_G(x,y)\leq K \,\mbox{and}\, \rho(z)=\rho(x)\}.$$
\begin{lemma} \label{l3f7}
Let $K>0$ and $Q$ be a finite set.  Then, there exists a $(K,d)$-large set $M=\{1,2,\dots,m\}$ and an oracle $\Theta: U^{K,M\oplus\{Q,*\}}_d\to\{a,b\}$, such that
for all $G\in\grd$, $\phi\in L_M^{(K)}(G)$ and $K$-subset function $\rho\in L_{\{Q,*\}}(G)$, the function
$f=O_\Theta(\phi\oplus\rho)$  has the following properties.
\begin{itemize}
\item $f^{-1}(a)\subset \rho^{-1}(Q)$.
\item For any $x\in V(G)$, such that $\rho(x)\neq *$
the set $f^{-1}(a) \cap V(T^\rho_x)$ is a maximal
independent subset of the graph $T^\rho_x$.
\end{itemize}
\end{lemma}
\proof Let $\Lambda$ be the finite set of not necessarily connected graphs (up to isomorphism) that occur as induced subgraph
in some ball $B\in U^K_d$.
Let $\hat{\Lambda}$ be the
set of all graphs with linearly ordered vertices such that the underlying
graph is an element of $\Lambda$. Finally, for any $\hat{J}\in \hat{\Lambda}$
fix a maximal independent subset $I_{\hat{J}}\subset V(\hat{J})$.
We describe $\Theta$ by the corresponding algorithm operator $O_\Theta$ in the following way.
If $\rho(x)\neq *$, let $O_\Theta(\phi\oplus \rho)(x)=a$ if in the ordered
graph $(T^\rho_x,\phi)$, $x\in I_{(T^{\rho}_x,\phi)}$. Note that
the relation $x\in I_{(T^{\rho}_x,\phi)}$ is well-defined. It is not hard to see
that $\Theta$ satisfies the conditions of our lemma. \qed
\vskip 0.1in
\noindent
Now let $\cG$ be a hyperfinite class of finite graphs.
Combining Theorem \ref{otostetel} and Lemma \ref{l3f7}, we can immediately
see that for any $\eps>0$, there exists a CTDA which for any $G\in\cG$ produces an induced subgraph
$H$ and a maximal independent subset $I_H\subset V(H)$ such that $(1-\eps) |V(G)|\leq |V(H)|$.
Let $J$ denote the restriction of $I_G$ onto $I_H$. Then, $|J|\leq |I_H|$ and $|I_G|\leq |J|+\eps |V(G)|$ hold.
Consequently, 
$$\frac{|I_H|}{|V(G)|}\geq \frac{|I_G|}{|V(G)|}-\eps,$$
\noindent
hence our proposition follows. \qed
\begin{remark} \label{subr}
The proof of Proposition \ref{hypind} illustrates the basic subroutines we use
to build CTDA's.
\begin{description}
\item[A] The subroutine finds a maximal $r$-separating system in the graphs $G\in\cG$.
\item[B1] For some $s\geq 1$,$ 2s<r$, the subroutine takes a symmetry breaking
  function $\phi\in L^{2r}_M(G)$ and a maximal $r$-separating system $T$ as
  inputs.
Then, independently label the balls $B_r(G,x)$, $x\in T$ by some set $Q$
using local computations in the balls.
\item[B2] A slight modification of the previous one. The subroutine uses the
  function $\phi\in L^{2r}_M(G)$, a maximal $r$-separating system $T$ and
a previously constructed labeling $\psi:V(G)\to P$ to label \\ the balls
$B_r(G,x), x\in T$.
\end{description}
\end{remark}
\begin{remark}
One can apply algorithm oracles to Cantor subshifts $Z\in \cC\gsg$ as
well. It is important to note that $O_\Theta(Z)$ is always qualitatively
weakly contained in $Z$.
\end{remark}

\subsection{Approximated 
maximum
matchings}

The goal of this section is to prove the following qualitative
analogue of the main results of \cite{ELip} and \cite{NO}.
A similar result using a different concept of local algorithms was proved by
Even, Medina and Ron \cite{EMR} (see also, \cite{Astr}).
\begin{theorem} \label{matching}
There exist CTDA's for the approximated maximum matching for finite
graphs in $\grd$.
\end{theorem}
\proof
Again, before getting into details, let us explain
the precise meaning of the theorem.
For  a $(5,d)$-large set $Q$, we call the function
$\rho:V(G)\to Q$ an matching function on $G$ if the
following conditions hold.
\begin{enumerate}
\item For any $a\in V(G)$, there exists at most one $b\in V(G)$
such that $0<d_G(a,b)\leq 5$ and $\rho(a)=\rho(b)$.
\item If  $0<d_G(a,b)\leq 5$ and $\rho(a)=\rho(b)$, then $a$ and $b$
are adjacent vertices.
\end{enumerate}
\noindent
Clearly, the set of adjacent pairs $(a,b)$ for which $\rho(a)=\rho(b)$ form
a matching $M_\rho$ of $G$.
The existence of $CTDA's$ for the approximated maximum matching
in finite graphs means that for
any $\eps>0$, there exist integers $1\leq n \leq l$,
an $(n,d)$-large set $N$, an $(5,d)$-large set $Q$ and an
oracle $\Theta: U^{l,n,N}_d\to Q$ such that for any finite
graph $G\in\grd$ and $\phi\in L^{(n)}_N(G)$,
\begin{itemize}
\item $\rho=O_\Theta(\phi)$ is an matching function, and
\item $\frac{|M_\rho|}{|V(G)|}\geq \frac{|M_G|}{|V(G)|}-\eps$, where
$M_G$ is the maximum sized matching in $G$.
\end{itemize}
\noindent
We will closely follow the proof of the main result in \cite{ELip}.
Informally speaking, the distributed algorithm work as follows.
First, we build a ``local improvement'' algorithm (see also \cite{NO})
which takes a matching $M$ as an input (together with the usual symmetry
breaking auxilliary  function) and
produces a new matching $M'$ such that $|M'|> |M|$ provided that
$M$ has an augmenting path of length shorter than $T$.
By Lemma 2.1 of \cite{ELip}, if a matching $M$ has  at most $\frac{\eps}{2} |V(G)|$
vertices from which an augmenting path shorter than $T$ starts, then \begin{equation} \label{pentek1}
\frac{|M_G|}{|V(G)|}\leq \frac{|M|}{|V(G)|}\frac{T+1}{T}+\frac{\eps}{2}\leq 
\frac{|M|}{|V(G)|}+\frac{1}{T}+\frac{\eps}{2}\,.
\end{equation}
\noindent
Thus, if $T>\frac{2}{\eps}$, the repeated applications of such a local improvement
algorithm lead to the required $(1+\eps)$-approximation of our theorem.
The construction of such algorithm is not very hard, however, we need to
show that the number of repetitions needed is bounded for the class of
finite graphs in $\grd$. So, let us write down the algorithm as a crude
pseudo-code, where each step requires a simple basic subroutine algorithm
as described in Remark \ref{subr}. Let $T >\frac{2}{\eps}$ be an integer.
\vi
{\bf 10} For each finite graph $G\in\grd$ we set up a starting
matching function $\rho:V(G)\to Q$.
\vi
{\bf 20} We construct a finite family $J_1, J_2,\dots, J_t$ of maximal
$10T$-separating
systems in $V(G)$ such that $\cup^t_{i=1} J_i=V(G)$. It is easy to see that
for large enough $t$ such algorithm exists for graphs with vertex degree bound
$d$.
\vi
{\bf 30} LET $j=1$.
\vi
{\bf 40} LET $i=1$. 
\vi
{\bf 50} IF $j=\mbox{``BOUND''}$ THEN GO TO {\bf 90}.
\vi
{\bf 60} IF $i=t+1$, LET $j=j+1$ and GO TO {\bf 40}.
\vi
{\bf 70} Consider the vertices $x\in J_i$, if there exists an augmenting path
starting from $x$ the subroutine makes the improvement inside the ball $B_{4T}(G,x)$ to
obtain a new matching function $\rho:V(G)\to Q$ representing more edges. If
there
is no such augmenting path the algorithm does not change the matching inside
the ball.
\vi
{\bf 80} LET i=i+1 and GO TO {\bf 50}.
\vi
{\bf 90} STOP.
\vi
Our theorem follows from the proposition below.
\begin{proposition} \label{pentekpro} 
If ``BOUND'' is a large enough integer, then for any finite
graph $G\in\grd$, when the algorithm stops we end up with a matching satisfying \eqref{pentek1}.
\end{proposition}
\proof
We apply an infinite-to-finite argument motivated by a similar proof in
\cite{ELip}.
Let $s_1<s_2<\dots$, $m_1<m_2<\dots$ and $k_1<k_2<\dots$ be positive integers
such that $M_n=\{0,1\}^{k_n}$ is an $(s_n,d)$-large set
and $\Theta_n: U_d^{m_n,s_n,M_n}\to Q$ is an oracle
that takes an element of $L_{M_n}^{(s_n)}(G)$, $G\in\grd$ as
input, and construct a matching as in the pseudo-code above until the variable
$j$ reaches $n$. Suppose that the statement of our
proposition does not hold. Then, there exists a sequence of finite
graphs $\{G_n\}^\infty_{n=1}$ such that for any $n\geq 1$ as the variable
``BOUND'' equals to $n$, when the algorithm stops the number $l_n$ of
vertices $x\in V(G_n)$ for which
an augmenting path shorter than $T$ starts at $x$, is greater
than $\frac{\eps}{2}|V(G_n)|$.
Thus, we have $\phi_n\in L_{M_n}^{(s_n)}(G_n)$ such
that for the matching function $O_\Theta(\phi_n)$ we have that
$$\frac{l_n}{|V(G_n)|}> \frac{\eps}{2}|V(G_n)|.$$
\noindent
Note that for any $n\geq 1$ we can regard the labeled graph $(G_n,\phi_n)$ as
an element of $\cC \grd$ by extending the labeling from
$\{0,1\}^{k_n}$ to $\{0,1\}^\N$ as zero for all the coordinates larger than
$k_n$.
Now, we recall the $\cC$-labeled version of the Benjamini-Schramm convergence
(Section 3.2. \cite{Elekhyper}).
Let $\{H_n,\psi_n\}^\infty_{n=1}\subset \cC \grd$ be a sequence of finite
graphs.
For any $B\in U^{r,\{0,1\}^r}_d$ let $T_{G,\phi}(B)\subset V(H,\psi)$ be the set of
vertices $x$ in $(H,\psi)\in\cC \grd$ such that the rooted-labeled ball of
radius $r$ around $x$ is rooted-labeled isomorphic to $B$. We say that
$\{H_n,\psi_n\}^\infty_{n=1}\subset \cC \grd$ is convergent in the sense of
Benjamini and Schramm if for any $r\geq 1$ and $B\in U^{r,\{0,1\}^r}_d$
$$\lim_{n\to\infty} \frac{|T_{H_n,\psi_n}(\cal{B})|}{|V(H_n)|}$$
\noindent
exists.
Clearly, one can pick a convergent graph sequence from any sequence of finite
graphs in $\grd$, so we can suppose that our counterexample sequence
$\{G_n,\phi_n\}^\infty_{n=1}$ is convergent in the sense of Benjamini and
Schramm. We can also suppose that $\{G_n,\phi_n\}^\infty_{n=1}$ is convergent
in $\cC \rg$.
Now, let us consider the infinite version of our proposition. 
Let $Z\subset\cC\grd$ be a proper subset such that for any $n\geq 1$, if
$(H,x,\psi)\in Z$ and $y,w\in V(H)$ for which $0<d_H(y,w)\leq s_n$, then
we have that $(\psi(y))_{[k_n]}\neq (\psi(w))_{[k_n]}$, where
$\{k_n\}^\infty_{n=1}$, $\{s_n\}^\infty_{n=1}$ are the sequences as above. Note that $Z$ can be regarded
as a Borel graphing (see \cite{ELip}). Observe that
the operators $O_{\Theta_n}(x)$ defines a Borel matching $\cM_n$ in $Z$.
\begin{lemma} \label{kedd12}
If $\mu$ is an invariant probability measure on the Borel graphing $Z$, then
there
exists an integer $m_Z>1$ such that 
the $\mu$-measure of elements $z\in Z$, for which an augmenting path (with
respect to the matching $\cM_{m_Z} $) shorter
than $T$ starts at $z$ is less than $\frac{\eps}{3}$.
\end{lemma}
\proof
The lemma follows from Proposition 1.1 \cite{ELip}, nevertheless we give
a short proof for completeness. Let $e=(y,w)$ be an edge of $(H,\psi,x)\in Z$.
In the matchings $\{\cM_n\}^\infty_{n=1}$ for some $n$'s $e$ belongs
to matching $\cM_n$ for some  $n$'s $e$ does not  belong
to the matching $\cM_n$. However, the membership of the edge $e$ 
stabilizes. Indeed, any change of the membership of the edge $e$
increases the number of matched vertices in the ball $B_{5T}(H,y)$.
Hence, we have a well-defined limit matching $\cM$ for which there
is no element $z\in Z$ with augmenting path shorter than $T$ starting at $z$.
Let $S_n\subset Z$ be the set of vertices $z\in Z$ (do not forget that $z$ is
a rooted-labeled infinite graph $(H,x,\phi)$), for which all edges adjacent to $z$ in the
graphing structure are stabilized after the $n$-th step. 
Since $\mu(S_n)\to 1$, our lemma follows. \qed
\vi
Now let $(G,\phi)\in\cC \grd$ be a limit
of the sequence $\{G_n,\phi_n\}^\infty_{n=1}\subset\cC\grd$.
Let $Z\subset \cC\grd$
be the orbit closure of $(G,\phi)$ in
$\cC\grd$ as in Section \ref{tgc}.
Then, we have an invariant probability measure $\mu$ on $Z$, for which
$$\lim_{n\to\infty}
\frac{|T_{G_n,\phi_n}(\cal{B})|}{|V(G_n)|}=\mu(T_Z(\cal{B}))$$
\noindent
holds for all $r\geq 1$ and $\cB\in U^{r,\{0,1\}^r}_d$.
Here, $T_Z(\cal{B})$ is the clopen set
of elements $(H,x,\psi)\in Z$ such that the rooted-labeled ball
$B_r(H,x,\psi_{[r]})$ is rooted-labeled isomorphic to $\cB$.
That is, the measured graphing $(Z,\mu)$ is the Benjamini-Schramm limit
of the sequence $\{(G_n,\phi_n)\}^\infty_{n=1}$. Also, for any ball $\cB$
the set $T_Z(\cB)$ is nonempty if and only if
the sets $\{T_{G_n,\phi_n}(\cal{B})\}^\infty_{n=1}$ are nonempty for all but finitely
many values of $n$.
Now let $m_Z$ be the constant in Lemma \ref{kedd12}.
Notice that the sequence of $Q$-labeled graphs $\{O_{\Theta_{m_Z}}(\phi_n)\}^\infty_{n=1}$
is convergent in the compact space $\grd^Q$.
Also, it is not hard to see that for all $r\geq 1$ and $\cB\in U^{r,Q}_d$, we
have that
$$\lim_{n\to\infty}
\frac{|T_{G_n,\phi_n}(\cal{B})|}{|V(G_n)|}=\mu(T_Z(\cB)).$$
\noindent
Therefore,
$$\lim_{n\to\infty}
\frac{|A_{G_n,\phi_n}(\cal{B})|}{|V(G_n)|}=\mu(A_Z(\cB)),$$
\noindent
where $A_{G_n,\phi_n}(\cal{B})$ is the set of vertices $x$ in $G_n$
for which an augmenting path
shorter than $T$ starts at $x$ in the matching
defined by the matching function 
$O_{\Theta_{m_Z}}(\phi_n)$.
Similarly, $A_Z(\cB)$ is
the clopen set of vertices $z\in Z$ for which an augmenting path
of $\cM_{m_Z}$ shorter than $T$ starts at $z$.
Since by Lemma \ref{kedd12}, $\mu(A)\leq \frac{\eps}{3}$ and
$$\limsup_{n\to\infty} \frac{|A_{G_n,\phi_n}(\cal{B})|}{|V(G_n)|}\geq
\frac{\eps}{2}\,,$$
we obtain a contradiction. Hence, our proposition follows and so does
our theorem. \qed

\subsection{The unrestricted weighted independent subset problem}
Let $G\in\grd$ be a finite graph.
Let $w\in V(G)\to \N\cup 0$ be an arbitrary
function The maximum $w$-weighted independent subset in $G$ is an independent
set $J\subset V(G)$ such that $\sum_{x\in J} w(x)$ is maximal among all
independent subsets of $G$. 
The $(1+\eps)$-approximated $w$-weighted independent
problem is to find an independent set $J\subset V(G)$ such that
\begin{equation} \label{vasar1}
\frac{\sum_{x\in J} w(x)}{\sum_{x\in V(G)} w(x)}\geq \frac{\sum_{y\in I_{G,w}}
    w(x)}{\sum_{x\in V(G)} w(x)}-\eps\,,
\end{equation}
\noindent
where $I_{G,w}$ is a maximum $w$-weight independent subset in $G$.
A local distributed algorithm algorithm for the approximated
weighted independent subset problem must deal
with arbitrarily large integers, hence it must use the
full power of the $\mathcal{LOCAL}$-model. The messages between 
the processors as well as the local computations are supposed to 
be unbounded. The deterministic-random local distributed algorithm
for the  approximated
weighted independent subset problem
for any $\eps>0$ takes a symmetry breaking function 
$\phi\in L_{M_\eps}^{(r_\eps)}(G)$ as an input and for any
$w:V(G)\to \N\cup 0$ in $r_\eps$-rounds it produces independent
subsets $J_1, J_2,\dots,J_{T_\eps}$ in $G$ in such a way that
if we randomly pick one of the $J_i's$ then the
probability of picking an independent subset satisfying \eqref{vasar1}
is larger than $(1-\eps)$.
\begin{proposition} \label{detran}
Let $\cG\subset\grd$ be a distributed strongly almost finite graph class
(e.g. $D$-doubling graphs by Theorem \ref{fotetel}) then
there exist deterministic-random local distributed algorithms for the
approximated weighted independent subset problem.
\end{proposition}
\proof
Let $G\in\cG$ and $\cal{E}=\cup^t_{j=1} H_j$ be a partition
of $V(G)$ such that
\begin{equation} \label{vasar3}
\sum_{x\in V^{\cal{E}}} w(x)\geq (1-\eps)\sum_{x\in V(G)} w(x)\,
\end{equation}
\noindent
where $V^{\cal{E}}=\cup^t_{j=1} (H_j\backslash \partial(H_j))\,.$
For each $j\geq 1$, pick a maximum $w$-weighted independent set
$I^{\cal{E}}_j$ in the graph induced on $H_j\backslash \partial(H_j)$ and
let $J^{\cal{E}}=\cup^t_{j=1} I^{\cal{E}}_j$. Then, we have that
$$\sum_{x\in J^{\cal{E}}} w(x)+\eps \sum_{x\in V(G)} w(x)\geq \sum_{x\in
  I_{G,w}} w(x)\,.$$
\noindent
That is,

$$\frac{\sum_{x\in J^{\cal{E}}} w(x)}{\sum_{x\in V(G)} w(x)}\geq \frac{\sum_{y\in I_{G,w}}
    w(x)}{\sum_{x\in V(G)} w(x)}-\eps\,.$$
\noindent
Repeating the argument of Lemma \ref{szombat1}, we can immediately see
that we have a CTDA for the class $\cG$ that produces, for all $G\in\cG$,
partitions
$$\{\cal{E_i}\}_{i=1}^{T_\eps}=\{H^i_1,H^i_2,\dots\}_{i=1}^{T_\eps}$$
\noindent
such that
the diameter of each class $H^i_j$ is bounded by some constant $L_\eps$,
independently of $G$, and the number of $i'$s for which \eqref{vasar3}
holds is larger than $(1-\eps) T_\eps$. Now our proposition follows,
since simple local algorithms can find maximum size $w$-independent sets in
the bounded diameter parts, by checking all independent subsets( using the
full power of the $\mathcal{LOCAL}$-model). \qed

\subsection{Distributed parameter testing}
Our notion of distributed parameter testing
can be viewed as the qualitative version of the randomized parameter testing
due to Goldreich and Ron \cite{GR}. 
First, recall the randomized parameter testing model for bounded
degree graphs. Let $p:\cG\to [a,b]$ or, in general, $p:\cG\to K$ be
a function such that $K$ is a compact metrizable space.
We say that $p$ is testable or estimable in the class of finite graphs
$\cG\in\grd$  if we have the following algorithm. 
\begin{enumerate}
\item First, we pick $r_\eps$ vertices of the graph $G\in\cG$ uniformly
  randomly and explore the $s_\eps$-neighbourhood of the picked vertices.
\item Then, the algorithm makes a ``guess'' $\hat{p}(G)$ in such a way that
$$\mbox{PROB}(d_K(\hat{p}(G),p(G))>\eps)<\eps\,.$$
\noindent
\end{enumerate}
\noindent
Now we define a {\bf distributed parameter testing} algorithm for the class 
of finite graphs $\cG\in\grd$.
\begin{enumerate}
\item First, the algorithm learn $\cB_{s_\eps}(G)$, the set
of all rooted balls of radius $s_r$ that occur in $G$.
Note that the algorithm will not know the probability distribution on 
 $\cB_{s_\eps}(G)$, as opposed to the case of randomized testability,
where we have a very good estimate on the distribution with
high probability by the Law of Large Numbers.
\item Based on the knowledge of $\cB_{s_\eps}(G)$, the
algorithm makes a guess $\hat{p}(G)$  in such a way that
the inequality 
$d_K(\hat{p}(G),p(G))\leq\eps$ always holds.
\end{enumerate}
\noindent
It is well-known that for a class $\cG\in\grd$ the parameter is testable in
the randomized setting if and only if for all Benjamini-Schramm convergent sequences
$\{G_n\}^\infty_{n=1}\subset \cG$ the limit $\lim_{n\to\infty} p(G_n)$ exists.
We have the following proposition for the distributed setting.
\begin{proposition} \label{testing}
The parameter $p:\cG\to K$ is testable in the distributed sense if and only if
for all naively convergent sequence $\{G_n\}^\infty_{n=1}\subset \cG$ the 
limit $\lim_{n\to\infty} p(G_n)$ exists.
\end{proposition}
\proof
Suppose that for all naively convergent sequence
$\{G_n\}^\infty_{n=1}\subset \cG$ the 
limit $\lim_{n\to\infty} p(G_n)$ exists.
Then, for any $\eps>0$ there exists $s_\eps>0$ such that
if $d_K(p(G),p(H))<\eps$, then $\cB_{s_\eps}(G)=\cB_{s_\eps}(H)$.
Then, we have the following testing algorithm.
For all family $\cB\subset U^{s_\eps}_d$ which equals to $\cB_{s_\eps}(G)$
for some graph $G\in \grd$, we pick a representative $G_{\cB}\in \grd$. 
The algorithm will make the guess $\hat{p}(G)=p(G_{\cB_{s_\eps}(G)})$.
\noindent
Now, suppose that for some naively convergent sequence
$\{G_n\}^\infty_{n=1}\subset \cG$
the limit $\lim_{n\to\infty} p(G_n)$ does not exist.
Then, there exists some $\eps>0$ such that for some sequences
$\{H_n\}^\infty_{n=1}\subset \cG$ and $\{J_n\}^\infty_{n=1}\subset \cG$,
$\cB_n(H_n)=\cB_n(J_n)$ and $d_K(p(H_n),p(J_n))>\eps.$ Therefore, it is
impossible to guess the value $p(G)$ based on the knowledge of $\cB_k(G)$ for
any $k\geq 1$. \qed

\noindent
By Theorem \ref{spectralconv}, we have the following
corollary.
\begin{corollary}
Let $\cG\subset\grd$ be a strongly almost finite class
of finite graphs. Then the parameter $p_{norm}:\cG\to [0,2d]$, $p_{norm}(G)=\|\Delta_G\|$
or the parameter $p_{spec}:\cG\to Cl([0,2d])$, $p_{spec}(G)=\Spec(\Delta_G)$
are testable parameters, where $Cl([0,2d])$ is
the compact set of all closed subsets of the interval $[0,2d]$ with the Hausdorff metric.
\end{corollary}
\section{Doubling and almost finiteness}\label{doublal}
The goal of this section is to prove the following theorem.
\begin{theorem} \label{fotetel}
For any $D>1$, the class of $D$-doubling graphs is distributed strongly almost finite.
\end{theorem}
\subsection{Doubling graphs}
Let $D$ be a positive integer. A graph $G$ of bounded vertex degree is called {\bf $D$-doubling} if 
for any $x\in V(G)$ and integer $s\geq 1$, $|B_{2s}(G,x)|\leq D|B_s(G,x)|$. 
Any doubling graph $G$ has polynomial growth of order $\log_2(D)$, that is, there exists
$C>0$ such that
$$|B_r(G,x)|\leq C r^{\log_2(D)}$$
holds for all $x\in V(G)$ and $r\geq 1$. In fact, the constant $C$ depends only on $D$ and the vertex degree bound of $G$.
Although it is not true that all the graphs of polynomial growth are doubling, graphs
of strict polynomial growth are always doubling. Recall that a graph $G$ is of
strict polynomial growth if there exists $\alpha\geq 1$, $0<C_1<C_2$ such that
$$C_1 r^\alpha\leq B_r(G,x) \leq C_2 r^\alpha$$
holds for all $x\in V(G)$ and $r\geq 1$.
The following lemma is well-known and we prove it only for completeness.
\begin{lemma} \label{measurecovering}
If $G$ is $D$-doubling then any ball $B_{2s}(G,x)$ can be covered by at most $D^4$-balls of radius $s$.
\end{lemma}
\proof
Let $z\in V(G)$ and $s>0$. We need to show that the ball $B_{2s}(G,z)$ can be covered by $D^4$ balls of radius $s$. Let $X\subset B_{2s}(G,z)$
be a maximal set of vertices such that
if $p\neq q\in X$, then $d_G(p,q)>s$. Hence, $B_{2s}(G,x)\subset \cup_{p\in X} B_s(G,p)$. So, $B_{2s}(G,x)$ can be covered by at most $|X|$ balls of radius $s$.
Also, if $p\in X$ then
\begin{equation}
\label{becsles}
B_{4s}(G,z)\subset B_{8s}(G,p)\,.
\end{equation}

\noindent
Furthermore, the balls $\{B_{s/2}(G,p)\}_{p\in X}$
are disjoint and contained in $B_{4s}(G,z)$.
By the $D$-doubling property, for any $p\in X$
$$|B_{8s}(G,p)|\leq D^4 |B_{s/2}(G,p)|\,.$$
\noindent
Hence by  \eqref{becsles} if $p\in X$, then 
$|B_{4s}(G,z)|\leq D^4 |B_{s/2}(G,p)|\,.$
Therefore, $|X|\leq D^4\,.$  \qed
\vskip 0.1in
\noindent
The following lemma is due to Lang and Schlichenmaier \cite{LS}.
\begin{lemma} \label{lang}
Let $G$ be a $D$-doubling graph.
Then for all integers $s,n\geq 1$ one has a system
$\cal{B}=\cup^{D^{(n+3)4}}_{i=1} \cal{B}_i$ such that
\begin{enumerate}
\item $\cal{B}_i$ is the union of balls of radius $s$ and if
$B_s(G,z)$ and $B_s(G,z')$ are two elements of  $\cal{B}_i$ so that $z\neq z'$, then
$B_{2^ns}(G,z)\cap B_{2^ns}(G,z')=\emptyset$.
\item The elements of $\cal{B}$ covers all the vertices of $G$.
\end{enumerate} \end{lemma}
\proof
Let $Z\subset V(G)$ be a maximal set of vertices such 
that $d_G(z,z')>s$ if $z,z'\in Z$ and $z\neq z'$.
Then the family $\cal{B}:=\{B_s(G,z)\}_{z\in Z}$ covers $V(G)$.
By the previous lemma, for each $z\in Z$ the ball $B_{2^{n+2}s}(G,z)$
can be covered by $D^{(n+3)4}$ balls of diameter less or equal $s$.
Each of these balls contains at most one element of $Z$.
Therefore, any ball $B_{2^{n+2}s}(G,z)$, $z\in Z$ contains at most $D^{(n+3)4}$ elements
of $Z$. Hence there exists a colouring
$\chi:Z\to\{1,2,\dots, D^{(n+3)4}\}$ such that
$\chi(z)\neq \chi(z')$, whenever $z,z'\in Z$ and $d_G(z,z')\leq 2^{n+2}s$.
That is, the family of balls $\cal{B}_i=\{B_s(G,z), \chi(z)=i\}$
satisfies the conditions of our lemma. \qed
\vskip 0.15in
\noindent
Let $H$ be a finite subset of the graph $G$. Then $\partial_K(H)$ denotes the set of vertices $x$ of $H$ for which
there exists a vertex $y\notin H$ such that $d_G(x,y)\leq K$. Also, $B_K(H)$ denotes the set of vertices $y$ of $G$ for which
there exists a vertex $x\in H$ such that $d_G(x,y)\leq K$.
The following lemma is a straightforward consequence of the definitions.
\begin{lemma}\label{trivial}
Let $G$ be a $D$-doubling graph of vertex degree bound $d$, $K>0$, $0<p<1$ and $\delta>0$.
Then, there exists an integer $M=M_{D,d,K,p,\delta}$ such that for any $x\in V(G)$ and $s\geq M$,
$$|\{ t\, \mid\, s\leq t < 2s\,, \frac {|B_{t+K}(G,x)|}{|B_{t-K} (G,x)|}<1+\delta\}| > ps\,.$$
\end{lemma}
\noindent
That is, if $s$ is large enough, then most of the balls in the form $B_t(G,x)$, $s\leq t < 2s$ has small boundary. 
\subsection{The Basic Algorithm}\label{sub31}
\noindent
Let $0<\eps<\frac{1}{2}$ be a real constant and $D$ be a positive integer. We call an $N$-tuple of positive integers \\
$S_1> S_2 >S_3>\dots>S_N$ $(D,\eps)$-good if for any $D$-doubling graph $G$, any integer $1\leq i \leq N$ and
$q\in V(G)$, there exists an integer $S_i\leq r_i(q) < 2S_i$ such that
\begin{equation} \label{eq1}
\frac{|B_{r_i(q)+16N S_{i+1}}(q)|}{|B_{r_i(q)}(q)|}<1+\frac{\eps}{10D^3}\,\,,\mbox{if $1\leq i \leq N-1$}. \end{equation}
\begin{equation} \label{eq1pont5}
i_G(B_{r_i(q)}(q))<\eps\,\,,\mbox{if $1\leq i \leq N$}.  \end{equation}
\begin{equation} \label{eq2}  \mbox{For any}\,\, 1\leq i\leq N-1, S_i>4S_{i+1}+4S_{i+2}+\dots +4S_N\,. \end{equation}
\begin{equation} \label{eq3}   (1-\frac{1}{4D^3})^N<\eps\,. \end{equation}
\noindent
The existence of such $(D,\eps)$-good tuples easily follows from Lemma
\ref{trivial}. Note that since all the balls in this section are in our graph
$G$, we will use the notation $B_r(x)$ instead of $B_r(G,x)$.
\vskip 0.2in
\noindent
{\bf The preliminary round.} Let $S_1> S_2 >S_3>\dots>S_N$ be a fixed $(D,\eps)$-good $N$-tuple and $G$ be
a $D$-doubling graph. For each integer $1\leq i \leq N$ pick a maximal system of vertices $\{q^i_\alpha\}_{\alpha\in J_i}$
in $G$ in such a way that if $\alpha\neq\beta$ then $d_G(q^i_\alpha,q^i_\beta)> 8 S_i\,.$ Using our definition
of $(D,\eps)$-goodness, for each chosen vertex $q^i_\alpha$ we pick an integer $S_i\leq r_i(q^i_\alpha) < 2S_i$ such that
\begin{equation}
\label{hatar} \frac{|B_{r_i(q^i_\alpha)+16NS_{i+1}}(q^i_\alpha)|}{|B_{r_i(q^i_\alpha)}(q^i_\alpha)|}<
1+\frac{\eps}{10D^3}\,,\mbox{if $1\leq i \leq N-1$}\,.
\end{equation}
and 
\begin{equation}
\label{hatar2}i_G(B_{r_i(q^i_\alpha)}(q^i_\alpha))<\eps \,,\mbox{if $1\leq i \leq N$}\,.
\end{equation}
For $1\leq i \leq N$ we  call the chosen balls $B_{r_i(q^i_\alpha)}(q^i_\alpha)$ balls of type-$i$.
\vskip 0.1in
\noindent
{\bf The construction round.} 
First, we discard all balls of type-$2$ that are intersecting a chosen ball of type-$1$. Inductively,
we discard all chosen balls of type-$i$ that are intersecting a chosen ball $B$ of type-$j$, $j<i$ such that
$B$ has not been previously discarded. Finishing the process we obtain a disjoint system of balls $B^1, B^2,\dots$, 
which we
call {\bf nice} balls.
Our main technical proposition goes as follows.
\begin{proposition} \label{fopropo}
Let $A\subset G$ be a finite subset such that
\begin{equation} \label{hatarA} 
\frac{|\partial_{16N S_1}(A)|}{|A|}<\frac{\eps}{10D^3},
\end{equation}
\noindent
then we have that
$$|\cup_{B^i\subset A} B^i|\geq (1-\eps) |A|\,.$$
\end{proposition}
\proof
First we need a lemma.
\begin{lemma}
\label{lemma1}
Let $1\leq i \leq N$ and $Q\subset G$ be a finite subset such that
$$\frac{|\partial_{16S_i}(Q)|}{|Q|}<\frac{1}{2},$$
then
$$\sum_{q^i_\alpha, q^i_\alpha \in Q\backslash \partial_{8S_i}(Q)} |B_{r_i(q^i_\alpha)}(q^i_\alpha)|> \frac{1}{2D^3}
|Q|\,.$$
\end{lemma}
\proof
Observe that
\begin{equation}
\label{cover}
\cup_{q^i_\alpha, q^i_\alpha \in Q\backslash \partial_{8S_i}(Q)}B_{8S_i}(q^i_\alpha) \supset Q\backslash \partial_{16S_i}(Q)\,.
\end{equation}
\noindent
Indeed, let $x$ be a vertex in $Q\backslash \partial_{16S_i}(Q)$. Then, by the maximality of the system
$\{q^i_\alpha\}_{\alpha\in J_i}$, there exists $\beta\in J_i$ such that
$x\in B_{8S_i}(q^i_\beta)$. Clearly, $q^i_\beta\in Q\backslash \partial_{8S_i}(Q)\,.$
By \eqref{cover},
$$D^3 \sum_{q^i_\alpha, q^i_\alpha \in Q\backslash \partial_{8S_i}(Q)} |B_{r_i(q^i_\alpha)} (q^i_\alpha)|\geq | Q\backslash \partial_{16S_i}(Q)|\,,$$
\noindent
hence, our lemma follows. \qed
\vskip 0.2in
\noindent
Now we turn to the proof of our proposition.
Let $A\subset G$ be a finite subset satisfying
the inequality \eqref{hatarA}. We define a process during which we inductively pick nice balls inside $A$ in such a way that
eventually the picked balls will cover at least $(1-\eps)|A|$ vertices of $A$.
Before starting our construction we need two technical lemmas.
\begin{lemma}
\label{lemmakicsi}
Let $A\subset G$ be a finite subset such that
$$\frac{|\partial_{16N S_1}(A)|}{|A|}<\frac{\eps}{10D^3}\,.$$
For $1\leq j \leq i$, let $C_j\subset A$ be the union of some nice balls of type-$j$.
Suppose that $|A\backslash \cup^i_{j=1} C_j|\geq \eps |A|\,.$
Then
\begin{equation}
\label{eqkicsi}
\frac{|\partial_{16 S_{i+1}}(A\backslash \cup^i_{j=1} C_j)|}{|A\backslash \cup^i_{j=1} C_j|}<\frac{1}{2}\,.
\end{equation}
\end{lemma}
\proof
Let $x\in \partial_{16 S_{i+1}}(A\backslash \cup^i_{j=1} C_j)$. Then at least one of the two conditions below are
satisfied.
\begin{itemize}
\item $x\in \partial_{16S_{i+1}}(A)\,.$
\item For some $1\leq j\leq i$ there exists $y\in C_j$, such that $d_G(y, C_j)\leq 16 S_{i+1}\,.$ That is
$x\in B_{16 S_{i+1}}(C_j)\backslash C_j\,.$
\end{itemize}
Therefore by \eqref{hatar},
$$|\partial_{16 S_{i+1}(A\backslash \cup^i_{j=1} C_j)}|\leq \frac{\eps}{10}|A|+\frac{\eps}{10}(\sum_{j=1}^i |C_j|)\leq \frac{\eps}{5}|A|\,.$$
\noindent
Hence, if $|A\backslash \cup^i_{j=1} C_j|\geq \eps |A|$, then \eqref{eqkicsi} holds. \qed
\begin{lemma}
\label{lemmanagy}
Let $A$ and $\{C_j\}^i_{j=1}$ be as in Lemma \ref{lemmakicsi}.
Let $$B=\partial_{16NS_{i+1}}(A)\cup \bigcup_{j=1}^i (B_{16NS_{i+1}}(C_j)\backslash C_j)\,.$$
\noindent
Then
\begin{equation}
\label{eqnagy}
|B|\leq \frac{\eps}{4D^3}|A|\,.
\end{equation}
\end{lemma}
\proof
Observe that by \eqref{hatarA},
$$|\partial_{16NS_{i+1}}(A)|\leq \frac{\eps}{10D^3} |A|\,.$$
\noindent
Also, if $L$ is a ball of type-$j$ $j\leq i$, then by \eqref{eq1}
$$|B_{16NS_{i+1}}(L)\backslash L |\leq \frac{\eps}{10D^3}|L|\,.$$
\noindent
Therefore
$$|\bigcup_{j=1}^i B_{16NS_{i+1}}(C_j)\backslash C_j|\leq \frac{\eps}{10D^3}|A|\,, $$
\noindent
hence \eqref{eqnagy} follows. \qed

\noindent
Now we define our covering process. First, we cover $A$ with balls of type-$1$, then
we inductively cover the rest of $A$ by smaller and smaller chosen balls. The trick is that
 that we use only smaller balls which are further and further away from the balls we used
previously. In this way, we can assure that all the balls we use in the covering process are nice.

\noindent
So, let $Q_0=A$.

\noindent
{\bf Step 1.} We pick all the chosen balls of type-$1$ in $A$ that does not intersect
$\partial_{4S_1}(A)$. These balls will be called $A$-covering balls.
\begin{itemize}
\item The union of $A$-covering ball of type-$1$ will be denoted by $C_1$.
\item We set $Q_1:=Q_0\backslash C_1$.
\end{itemize}

\noindent
{\bf Step 2.} We continue the covering process by picking all the
chosen balls $L$ of type-$2$ inside the set $A$ that does not intersect the set $\partial_{8S_1}(A)$ 
and for any $A$-covering ball $D$ of type-$1$, $L$ does not intersect $B_{4S_2}(D)$ either.
\begin{itemize}
\item Again, we will call the balls picked above $A$-covering balls and denote their union by $C_2$.
\item We set $Q_2:=Q_1\backslash C_2.$
\end{itemize}

\noindent
We will see in a moment that all the $A$-covering balls of type-$2$ are nice.

\noindent
{\bf Step (i+1).} In the first $i$-steps we have already defined disjoint 
sets \\
$C_1, C_2,\dots, C_i$ inside the set $A$, where for any $1\leq j \leq i$, the set $C_j$ is the union
of nice balls picked at the step $j$. The $A$-covering balls have the following properties:
\begin{itemize}
\item If $L$ is an $A$-covering ball of type-$j$ and $M$ is an $A$-covering ball of type-$k$, $k<j$
then $B_{4(j-k)S_{k+1}}(M)\cap L= \emptyset$.
\item If $L$ is an $A$-covering ball of type-$j$, then
$\partial_{4jS_1}(A)\cap L=\emptyset\,.$
\item $Q_j=Q_{j-1}\backslash C_j$.
\end{itemize}

\noindent
Now we continue our covering process. We pick all the chosen balls $P$ of type-$(i+1)$ contained in $A$ for
which both conditions below are satisfied:
\begin{enumerate}
\item $\partial_{4(i+1)S_1}(A)\cap P=\emptyset.$
\item For any $1\leq j \leq i$ and $A$-covering ball $M$ of type-$j$, $B_{4(i+1-j) S_{j+1}}(M)\cap P=\emptyset\,.$
\end{enumerate}
\noindent
Finally, we set $Q_{i+1}:=Q_i\backslash C_{i+1}$. Our crucial observation is formulated in the following lemma.
\begin{lemma} \label{covnice} All the $A$-covering balls of type-$(i+1)$ are nice.
\end{lemma}
\proof
Let $L$ be an $A$-covering ball of type-$(i+1)$. Suppose that $L$ is not nice.
Then there exists a nice ball $M$ of type-$j$, $j\leq i$ such that $M\cap L\neq \emptyset$.
In our construction, $A$-covering balls of type-$j$ cannot intersect
an $A$-covering ball  of type-$(i+1)$. That is, $M$ cannot be an $A$-covering ball.
Since $L\subset A$ and $\partial_{4(i+1)S_1}(A)\cap L=\emptyset$, the ball $M$ is a subset of $A$ as well.
Hence, the reason that $M$ has not picked at step $j$ was 
\begin{itemize}
\item either that $M\cap \partial_{4jS_1}(A)\neq \emptyset$
\item or that $M\cap B_{4(j-k)S_{k+1}}(D)\neq \emptyset$ for some $A$-covering ball $D$ of type-$k$.
\end{itemize}
\noindent
{\bf Case 1.}
\begin{enumerate}
\item There exists $x\in M$, $x\in\partial_{4jS_1}(A)$ and
\item there exists $y\in L$ such that $y\notin \partial_{4(i+1)S_1}(A)$ and $y\in M$.
\end{enumerate}
\noindent
That is, $d_G(x,y)\geq 4(i+1-j)S_1$ and $d_G(x,y)\leq \diam (M) \leq 4 S_j$ leading to a contradiction.
\vskip 0.1in
\noindent
{\bf Case 2.}
\begin{enumerate}
\item There exists an $A$-covering ball $D$ of type-$k$, $k<j$ and $x\in M$ such that $x\in B_{4(j-k)S_{k+1}}(D)$
\item there exists $y\in L\cap M$, $y\notin \partial_{4(i+1-k)S_{k+1}}(D).$
\end{enumerate}
\noindent
Then $d_G(x,y)> 4(i+1-j)S_{k+1}$ and
$d_G(x,y)\leq \diam (M) \leq 4S_j$
leading again to a contradiction.

\begin{lemma} \label{hanyados}
If $|Q_i|>\eps |A|$ then $|Q_{i+1}|<(1-\frac{1}{4D^3}) |Q_i|$.
\end{lemma}
\proof
Observe that $Q_i=(A\backslash \cup^i_{j=1} C_j).$ By Lemma \ref{lemmakicsi},
$$\frac{|\partial_{16 S_{i+1}} Q_i|}{|Q_i|}<\frac{1}{2}\,.$$
Hence by Lemma \ref{lemma1} we have that
$$\sum_{q^{i+1}_\alpha, q^{i+1}_\alpha\in Q_i\backslash \partial_{8 S_{i+1}}(Q_i)}
|B_{r_{i+1}(q^{i+1}_\alpha)}(q^{i+1}_\alpha)|>\frac{1}{2D^3} |Q_i|\,.$$
\noindent
Let $q^{i+1}_\alpha\in  Q_i\backslash \partial_{8 S_{i+1}}Q_i$
such that the ball $L=B_{r_{i+1}(q^{i+1}_\alpha)}(q^{i+1}_\alpha)$ does not
intersect the set
$$B'=\partial_{4(i+1)S_1}(A)\cup \bigcup_{j=1}^i (B_{4(j-i)S_{i+1}}(C_j)\backslash C_j)\,.$$
\noindent
Then, $L$ is an $A$-covering ball of type-$(i+1)$.
On the other hand, if $L$ does intersect $B'$, then
$$L\subset \partial_{16NS_1}(A)\cup \bigcup_{j=1}^i (B_{16NS_{i+1}}(C_j)\backslash C_j)\,.$$
\noindent
Therefore by Lemma \ref{lemmanagy},
$$|C_{i+1}|>\frac{1}{2D^3}|Q_i|-\frac{\eps}{4D^3} |A|>\frac{1}{4D^3} |Q_i|\,.$$
\noindent
That is,
$$|Q_{i+1}|=|Q_i\backslash C_{i+1}|< (1-\frac{1}{4D^3} )|Q_i|\,\quad\qed$$
\noindent
By Lemma \ref{hanyados} and \eqref{eq3}, Proposition \eqref{fopropo} immediately follows. \qed
\subsection{The class of D-doubling graphs is almost finite}
The goal of this subsection is to prove the following proposition.
\begin{proposition}\label{doubalm}
The class of $D$-doubling graphs is almost finite.
\end{proposition}
\proof
Let $G$ be a $D$-doubling countable graph with vertex degree bound $d$.
Fix $\eps>0$ and let $\eps'=\frac{\eps}{4D^{16}}$. Let
$S_1>S_2>\dots>S_N$ be a $(D,\eps')$-good $N$-tuple of integers. Let $\cal{N}_i$ be the set
of nice balls of type-$i$ obtained in the construction round of our Basic Algorithm of  Subsection
\ref{sub31}. Let $U\subset V(G)$ be the set of vertices $x$ such that
$x$ is not contained in any nice ball $L\in \cal{N}_i$, $1\leq i \leq N$.
By \eqref{eq1pont5}, $i_G(L)<\eps'$ if $L$ is a nice ball.
Following \cite{DHZ}, we will construct an injective map $\Psi:U\to V(G)\backslash U$ in
such a way that for any nice ball $L$,
$$i_G(L\ \cup \Psi^{-1}(L))<\eps\,.$$
\noindent
Also, we will assure that
$$\sup_{x\in V(G)} d_G(\Psi(x),x)<\infty\,.$$
\noindent
Thus, we have a covering
$V(G)=\bigcup_{1\leq i \leq N} \bigcup_{L\in\cal{N}_i} (L\cup \Psi^{-1}(L))$
of the vertices of $G$ with disjoint sets of bounded diameter with isoperimetric constant less or equal than $\eps$.
This will show that $G$ is almost finite.

\noindent
Now, let us construct $\Psi$. First of all, 
by Lemma \ref{lang} and Lemma \ref{trivial}, we have an
integer $s>0$ and  a system $\cal{B}=\cup_{i=1}^{D^{16}} \cal{B}_i$ so that
\begin{itemize}
\item The elements of $\cal{B}_i$ are disjoint balls $B$ such that
$s\leq \mbox{radius of}\, B \leq 2s$ and $|\partial_{16 N S_1}(B)|< \frac{\eps'}{10D^3} |B|\,.$ Hence
by Proposition \ref{fopropo}, for any $B\in \cal{B}_i$, we have that $|U\cap B|<\eps'|B|.$
\item $\cal{B}$ covers $V(G)$.
\end{itemize}
\noindent
For a nice ball $L$ and $1\leq i \leq D^{16}$, let $L_1, L_2,\dots, L_{D^{16}}$
be disjoint
subsets of $L$ such that
$$2\eps' |L|< |L_i|<  3\eps' |L|\,.$$
\noindent
Thus for any $1\leq i \leq D^{16}$ and $B\in\cal{B_i}$
$$|U\cap B|< |\cup_{L\subset B} L_i|\,.$$
\noindent
Let $\phi_i:U\cap(\bigcup_{B\in \cal{B_i}} B)\to \bigcup_{B\in \cal{B_i}}
\bigcup_{L\subset B} L_i$ be an arbitrary injective map such that if $x\in B_i$
then $\phi_i(x)\in B_i$.
Our injective map $\Psi:U\to V(G)\backslash U$ is defined as follows.
If $x\in U$, let $\Psi(x)=\phi_i(x)$ if
$i$ is the smallest integer such that $x\in U\cap B$ for some ball $B\in\cal{B_i}$.
Hence, if $L$ is a nice ball then
$$\partial_G(L\cup \Psi^{-1}(L)|\leq |\partial_G(L)| + |\Psi^{-1}(L)|\leq $$
$$\leq |\partial_G(L)|+ 4\eps' D^{16} |L| \leq \eps |L|\,.$$
\noindent
Therefore our proposition follows. \qed
\vskip 0.1in
\begin{remark}
Recently, Downarowicz and Zhang \cite{DZ} proved that Cayley graphs of groups
of subexponential growth are distributed almost finite. Nevertheless, they proof used 
the transitivity of the graphs in a significant way. 
\end{remark}
\subsection{The class of D-doubling graphs is strongly almost finite}
In this subsection we go one step further and strengthen Proposition \ref{doubalm}.
\begin{proposition}\label{doubstrongalm}
For any $D\geq 1$, the class of $D$-doubling graphs is even strongly almost finite.
\end{proposition}
\proof
We need to modify the Basic Algorithm  of Subsection \ref{sub31} to fit our purposes.
Let $0<\delta,\eps<1$ be real constants and $D$ be a positive integer. Also, let $N>0$ be
so large that $(1-\frac{1}{4D^3})^N<\eps$, $(1-\frac{1}{D^{20}})^N<\frac{\delta}{2}$.
We call an $N$-tuple of positive integers $S_1>S_2>\dots> S_N$ $(D,\eps,\delta)$-good if there
exist positive integers $\{R_i\}^N_{i=1}$ such that
for any $D$-doubling graph $G$, for any integer $1\leq i \leq N$ and $q\in V(G)$, there exists an
integer $S_i\leq r_i(q)< 2S_i$, such that $r_i(q)+R_i< 2S_i$ and
\begin{equation} \label{neq1}
\frac{|B_{r_i(q)+j+20NS_{i+1}}(q)|}{|B_{r_i(q)_j}(q)|}<1+\frac{\eps}{10D^3},\quad\mbox{if $1\leq i \leq N-1$},\, 1\leq j\leq R_i.
 \end{equation}
\begin{equation} \label{neq1pont5}
i_G(B_{r_i(q)_j}(q))<\eps,\quad\mbox{if $1\leq i \leq N$},\, 1\leq j\leq R_i.  \end{equation}
\begin{equation} \label{neq2}  \mbox{For any}\,\, 1\leq i\leq N-1,\,\, S_i>4S_{i+1}+4S_{i+2}+\dots +4S_N\,. \end{equation}
\begin{equation} \label{neq3}   \mbox{For any}\,\, 1\leq i\leq N-1,\,\,\frac{10 S_{i+1}}{R_i}< \frac{\delta}{2N}\,.  \end{equation}
\noindent
Again, the existence of such $(D,\eps,\delta)$-good tuples follows from Lemma
\ref{trivial}.
Now, we fix a  $(D,\eps,\delta)$-good tuple $S_1>S_2>\dots> S_N$ and a system of integers $\{R_i\}^N_i$ satisfying the
conditions above. Using Lemma \ref{lang}, for each $1\leq i \leq N$ we pick $D^{20}$ maximal systems of vertices 
$\{\{q^{i,t}_\alpha\}_{\alpha\in J_{i,t}}\}_{t=1}^{D^{20}}$ such that
\begin{itemize}
\item for any $1\leq t \leq D^{20}$ and $\alpha\neq \beta$, $d_G(q^{i,t}_\alpha, q^{i,t}_\beta)>8S_i\,,$
\item the balls $\{\{B_{S_i}(q^{i,t}_\alpha)\}_{\alpha\in J_{i,t}}\}^{D^{20}}_{t=1}$ cover $V(G)$.
\end{itemize}
\noindent
Then, using the definition of $(D,\eps,\delta)$-goodness, for each $1\leq i \leq N$, $1\leq t\leq D^{20}$ and
vertex $q^{i,t}_\alpha$, we pick an integer $S_i\leq r_i(q^{i,t}_\alpha)< 2S_i$ such that
$r_i(q^{i,t}_\alpha)+R_i< 2S_i$ and
\begin{equation} \label{ujneq1}
\frac{|B_{r_i(q^{i,t}_\alpha)+j+20NS_{i+1}}(q)|}{|B_{r_i(q^{i,t}_\alpha)+j}(q)|}<1+\frac{\eps}{10D^3},\quad\mbox{if $1\leq i \leq N-1$},\, 1\leq j\leq R_i.
 \end{equation}
\begin{equation} \label{ujneq1pont5}
i_G(B_{r_i(q^{i,t}_\alpha)+j}(q))<\eps,\quad\mbox{if $1\leq i \leq N$},\, 1\leq j\leq R_i.  \end{equation}
\noindent
We call an $N$-tuple $\{(t_i,j_i)\}^N_{i=1}$ a {\bf code} if for any $1\leq i \leq N$ we have that $1\leq t_i\leq  D^{20}$ and
$1\leq j_i\leq R_i$. Thus, the set of codes $\cM$ has $D^{20N} \prod_{i=1}^N R_i$ elements.
For each code $\{(t_i,j_i)\}^N_{i=1}$ and $1\leq i \leq N$, we have a system of chosen balls
with centers $\{q^{i,t_i}_\alpha\}_{\alpha\in J_{i,t_i}}$ and radii $\{r_i(q^{i,t_i}_\alpha)+{j_i}\}_{\alpha\in J_{i,t_i}}$ and
an associated $(\eps, K_\eps)$-partition given by the Basic Algorithm for a certain integer $K_\eps>0$ that does not depend
on $G$. By the definition of strong almost finiteness, Proposition \ref{doubstrongalm} easily follows from the proposition 
below.
\begin{proposition} \label{segedprop} Let $x\in V(G)$ and $C_x\in \cM$
be the set of codes $\{(t_i,j_i)\}^N_{i=1}$ for which
there exists a nice ball $B$ given by the Basic Algorithm, such that $x\in B$ and $x\notin \partial_G(B)$.
Then, $|C_x|\geq (1-\delta)|\cM|$.
\end{proposition}
\proof
Let $D_x$ be the set of codes $\{(t_i,j_i)\}^N_{i=1}$ for which the chosen balls
given by the Basic Algorithm are not covering $x$.
Also, let $E_x$ be the set of codes  $\{(t_i,j_i)\}^N_{i=1}$ for which there
exists $1\leq i \leq N-1$ and
a chosen $i$-ball $B=B_{r_i(q^{i,t}_\alpha)+j_i}(q^{i,t}_\alpha)$ so
that
\begin{itemize}
\item $B\cap B_{5S_{i+1}}(x)\neq \emptyset$ and
\item $B$ does not contain the ball $B_{5S_{i+1}}(x)$.
\end{itemize}
\begin{lemma} \label{hl1}
If  $\{(t_i,j_i)\}^N_{i=1}\notin D_x\cup E_x$, then $\{(t_i,j_i)\}^N_{i=1}\in C_x$.
\end{lemma}
\proof
Suppose that $\{(t_i,j_i)\}^N_{i=1}\notin D_x\cup E_x$. Let $1\leq i \leq N$ be the smallest integer
such that there exists a chosen $i$-ball $B$ such that $x\in B$. Since our code is not in $E_x$, $x\notin \partial_G(B)$. Also,
there is no chosen $j$-ball, $j<i$ which intersects $B$. Hence, $B$ is a nice ball and therefore $\{(t_i,j_i)\}^N_{i=1}\in C_x$.
\begin{lemma} \label{hl2}
$|D_x|<\delta/2 |\cM|\,.$
\end{lemma}
\proof
For any $1\leq i \leq N$, there exists at least one $1\leq t \leq D^{20}$ such
that $x\in\cup_{\alpha\in J_{i,t}} B_{S_i}(q_\alpha^{i,t})$.
Hence, $|D_x|\leq \prod^N_{i=1} (D^{20}-1) R_i.$ Thus,
$$\frac{|D_x|}{|\cM|}<(1-\frac{1}{D^{20}})^N<\frac{\delta}{2}\,.\,\,\qed$$
\begin{lemma} \label{hl3}
$|E_x|<\delta/2 |\cM|\,.$
\end{lemma}
\proof Fix an $N$-tuple $\{t_i\}^N_{i=1}$, $1\leq t_i\leq D^{20}$.
Then, for any $1\leq i\leq N-1$, there exists
at most one $\alpha\in J_{i,t}$ such that both
\begin{equation}\label{eqma1}
B_{r_i(q^{i,t}_\alpha)+j}(q^{i,t}_\alpha)\cap B_{5S_{i+1}}(x)\neq\emptyset
\end{equation}
\begin{equation}\label{eqma2}
B_{r_i(q^{i,t}_\alpha)+j}(q^{i,t}_\alpha)\not\supset B_{5S_{i+1}}(x)
\end{equation}
hold for a certain element $1\leq j \leq R_i$.
Clearly, the number of $j$'s for which both \eqref{eqma1} and \eqref{eqma2} hold is
not greater than $10 S_{i+1}$.
Thus,
$$ |E_x|\leq D^{20N} (\sum_{j=1}^{N-1}\frac{10 S_{j+1}}{R_j})\prod^N_{i=1} R_i$$

\noindent Hence our lemma follows. \qed
\vskip 0.1in
\noindent
Now, Lemma \ref{hl1}, Lemma \ref{hl2} and Lemma \ref{hl3} immediately
imply our proposition. \qed

\subsection{Distributed strong almost finiteness}
Now, we finish the proof of Theorem \ref{fotetel}. It is not hard to see that 
all the constructions in Proposition \ref{doubalm} and Proposition
\ref{doubstrongalm} can be done locally. Nevertheless, we show step by step how
to build the partition families witnessing strong almost finiteness in a
distributed fashion, using the simple subroutines described in Remark
\ref{subr}.
\begin{enumerate}
\item We pick the $D^{20}$ maximal systems of vertices
  $\{\{q^{i,t}_\alpha\}_{\alpha\in J_{i,t}}\}^{D^{20}}_{t=1}$ using an
(A)-type local algorithm.
\item Using (B)-type local algorithms in the balls around the picked vertices
we construct the chosen balls.
\item Again, using (B)-type local algorithms we construct the system of nice
  balls.
\item Now, we construct the system of balls $\cup_{i=1}^{D^{16}} \cB_i$ as in
  the proof of Proposition \ref{doubalm} using (A)- and (B)- type local
  algorithms.
\item We construct the injective maps $\Psi$ as in the proof of Proposition
  \ref{doubalm} for each nice ball system using a (B)-type local algorithm.
\end{enumerate}
\vi
This finishes the proof of our theorem. \qed

\end{document}